\documentclass{cas-sc}
\pdfoutput=1

\usepackage{appendix}
\usepackage{graphbox}
\usepackage{mathtools}
\usepackage{dashrule}
\usepackage{enumitem}
\usepackage[numbers]{natbib}
\usepackage{lipsum}
\usepackage{caption}
\usepackage{subcaption}
\usepackage{cleveref}
\usepackage{tabularx}
\usepackage{amsthm}

\definecolor{stability_dark}{RGB}{102, 102, 102}
\definecolor{stability_light}{RGB}{204, 204, 204}

\definecolor{matlab_blue}{RGB}{0, 114, 189}
\definecolor{matlab_red}{RGB}{217, 83, 25}
\definecolor{matlab_yellow}{RGB}{237, 177, 32}
\definecolor{matlab_purple}{RGB}{126, 47, 142}
\definecolor{matlab_green}{RGB}{119, 172, 48}
\definecolor{matlab_azure}{RGB}{77, 190, 238}
\definecolor{matlab_dark_red}{RGB}{162, 20, 47}
\definecolor{nls_ic_plot_gray}{RGB}{180, 180, 180}

\definecolor{nls_fourier_gray}{RGB}{128, 128, 128}
\definecolor{nls_fourier_red}{RGB}{192, 57, 43}
\definecolor{nls_fourier_blue}{RGB}{51 153 255}
\definecolor{nls_fourier_black}{RGB}{30, 30, 30}

\setcitestyle{square}
\setcitestyle{citesep={,}}
\hypersetup{
 colorlinks=true,
 citecolor=Cerulean,
 linkcolor=Cerulean,
 urlcolor=Cerulean}
\crefname{equation}{}{} 
 
\newtheorem{remark}{Remark}

\begin{document}

\title[mode = title]{Exponential Runge-Kutta Parareal for Non-Diffusive Equations }
\tnotemark[1]
\tnotetext[1]{This work was funded by the National Science Foundation, Computational Mathematics Program DMS-2012875}

\shorttitle{Exponential Runge-Kutta Parareal}
\shortauthors{T Buvoli et~al.}

\author[1]{Tommaso Buvoli}
\cormark[1] 
\ead{tbuvoli@tulane.edu} 
\address[1]{Tulane University, New Orleans, LA 70118, USA}

\author[2]{Michael Minion}
\ead{mlminion@lbl.gov} 
\address[2]{Lawrence Berkeley National Lab,Berkeley, CA 94720, USA}
 
\cortext[cor1]{Corresponding author} 

\begin{abstract}
	Parareal is a well-known parallel-in-time algorithm that combines a coarse and fine propagator within a parallel iteration. It allows for large-scale parallelism that leads to significantly reduced computational time compared to serial time-stepping methods. However, like many parallel-in-time methods it can fail to converge when applied to non-diffusive equations such as hyperbolic systems or dispersive nonlinear wave equations. This paper explores the use of exponential integrators within the Parareal iteration. Exponential integrators are particularly interesting candidates for Parareal because of their ability to resolve fast-moving waves, even at the large stepsizes used by coarse propagators. This work begins with an introduction to exponential Parareal integrators followed by several motivating numerical experiments involving the nonlinear Schr\"{o}dinger equation. These experiments are then analyzed using linear analysis that approximates the stability and convergence properties of the exponential Parareal iteration on nonlinear problems. The paper concludes with two additional numerical experiments involving the dispersive Kadomtsev-Petviashvili equation and the hyperbolic Vlasov-Poisson equation. These experiments demonstrate that exponential Parareal methods offer improved time-to-solution compared to serial exponential integrators when solving certain non-diffusive equations. 
\end{abstract}

\maketitle

\def\Oof {{\mathcal{O}}}
\def\Fprop {{\mathcal{F}}}
\def\Gprop {{\mathcal{G}}}
\def\y {{y}}
\def\yn {\y_n}
\def\ynk {{\y_n^k}}
\def\ynkp {{\y_n^{k+1}}}
\def\ynp {\y_{n+1}}
\def\tn  {{t_n}}
\def\tnp  {{t_{n+1}}}
\def\Nstepstot {{N_s}}
\def\Nsteps {{N_T}}
\def\Nprocs {{N_p}}
\def\Nblocks {{N_b}}
\def\NFprop {{N_f}}
\def\NGprop {{N_g}}
\def\Tzero {{t_0}}
\def\Tfinal {{t_{\text{final}}}}
\def\Tblock {{\Delta T}}
\def\Tproc {{\Delta T_p}}
\def\CostG  {{C_\Gprop}}  %
\def\CostF  {{C_\Fprop}}  %
\def\CostB  {{C_B}}  %
\def\CostS  {{C_s}}  %
\def\CostP  {{C_p}}  %
\def\CostGstep  {{c_g}}  %
\def\CostFstep  {{c_f}}  %

\def\Speedup  {{S}}
\def\Eff  {{E}}
\def\u  {{u}}
\def\g  {{u}}
\def\uhat  {{\hat{\u}}}
\def\ghat  {{\hat{\g}}}
\def\uhatk  {{\hat{\u}_k}}
\def\ghatk  {{\hat{\g}_k}}
\def\Tfin  {{T}}
\def\cRate {\|\mathbf{E}\|}
\def\NIters {{K}}
\def\Dt  {{h}}
\def\DT  {{\Delta T}}
\def\DtF  {{\Delta t_\Fprop}}
\def\DtG  {{\Delta t_\Gprop}}
\def\K  {{K}}
\def\barK  {{\bar{\K}}}
\def\A {\mathcal{A}} %

\definecolor{plot_blue}{RGB}{41, 128, 185}
\definecolor{plot_orange}{RGB}{211, 84, 0}
\definecolor{plot_red}{RGB}{192, 57, 43}
\definecolor{plot_green}{RGB}{39, 174, 96}

\newcommand{\cdline}[1]{{\color{gray}  #1.} \hspace{0.5em}} 
\section{Introduction}
\label{sec:introduction}

Time integrators \cite{hairer1993solving,wanner1996solving,butcher2016numerical} are numerical methods that solve an initial value problem by sequentially advancing the solution via a series of discrete timesteps. For more than half a century, these methods have proven invaluable for modeling a range of dynamical processes appearing in both science and engineering. In a typical calculation one iteratively applies a time integrator to evolve a system over thousands or even millions of timesteps. Therefore, the total computational cost is not just that of a single timestep, but rather the combined cost of applying the method sequentially over the full temporal domain.

The sequential nature of classical time-stepping methods has come under increasing scrutiny in light of modern parallel hardware like multicore processors, massively parallel high performance computing systems, and specialized accelerators. For more than two decades, these advancements have spurred the development of new parallel-in-time (PinT) methods \cite{HortonVandewalle1995,LionsEtAl2001,EmmettMinion2012,FriedhoffEtAl2013,Gander2015_Review} that distribute the full temporal domain over a large number of computational nodes. Perhaps the most well-known PinT method is the Parareal algorithm \cite{LionsEtAl2001}. Parareal consists of a parallel iteration that combines a fine propagator (a computationally expensive and accurate integrator)  with a coarse propagator (a computationally cheap and less accurate integrator). The aim of Parareal is to obtain the solution of the fine propagator at a similar computational cost to that of running the coarse propagator as a serial one-step method.

Parareal has proven effective for accelerating the solution of diffusive equations \cite{FischerEtAl2005,Trindade2004,nielsen2012feasibility,KreienbuehlEtAl2015} and its theoretical convergence properties are well understood in the presence of diffusion \cite{GanderVandewalle2007_SISC}. In contrast, non-diffusive equations (e.g. hyperbolic systems or dispersive nonlinear wave equations) introduce significant numerical difficulties that lead to slow convergence or instabilities in the Parareal iteration \cite{Bal2005,StaffRonquist2005,GanderHalpern2017,Ruprecht2018,buvoli2021imexparareal}. Though numerous modifications have been proposed \cite{ChenEtAl2014,DaiEtAl2013,EghbalEtAl2016,FarhatEtAl2003,GanderPetcu2008,KooijEtAl2017}, the resulting methods introduce additional complexities that make them less applicable to all types of problems.

For unmodified Parareal, rapid convergence on non-diffusive equations requires that the coarse propagator accurately approximate the fine propagator \cite{Ruprecht2018,buvoli2021imexparareal}. However, classical integrators have difficulty resolving rapid oscillations at coarse stepsizes due to phase errors and numerical diffusion. For this reason, we are motivated to consider exponential integrators \cite{hochbruck2010exponentialreview} that treat a linear component exactly and possess the ability to accurately resolve fast moving waves even at coarse stepsizes.

In the past two decades, multiple families of exponential integrators have been developed to  efficiently solve both diffusive and non-diffusive equations \cite{beylkin1998ELP,cox2002ETDRK4,krogstad2005IF,hochbruckostermann2005ETDRKSTIFFB,luan2014explicit,buvoli2019esdc,ostermann2006general,buvoli2021epbm,crouseilles2020exponential, gaudreault2022high, hamon2020parallel, peixoto2019semi}. The majority of the works discussing the parallelization of exponential integrators focus on the computation of exponential matrix functions \cite{schreiber2019exponential,caliari2021accurate} and on parallel function evaluations or output computations \cite{buvoli2021epbm,luan2016parallel}. Within the context of PinT methods, asymptotic averaging techniques that are closely related to exponential integration have been proposed for oscillatory problems \cite{haut2014asymptotic,peddle2019parareal}, and exponential Krylov methods are used in the ParaExp alogorithm \cite{gander2013paraexp}.

In this work we combine exponential integrators \cite{hochbruck2010exponentialreview} with Parareal and demonstrate, both theoretically and experimentally, that this pairing can provide reduced time-to-solution on non-diffusive equations. We focus specifically on the non-diffusive, semilinear initial value problem 
	\begin{align}
		\mathbf{y}' = \mathbf{Ly} + N(t,\mathbf{y}), \quad \mathbf{y}(t_0) = \mathbf{y}_0
	\end{align}
where the eigenvalues of $\mathbf{L}$ are purely imaginary. Exponential integrators treat the linear component $\mathbf{L}$ exactly, granting them the ability to accurately resolve fast moving waves even at coarse stepsizes. However, special care must be taken when applying exponential integrators on non-diffusive equations since the methods are classically unstable  \cite{buvoli2022stability,crouseilles2020exponential}. In fact, we will demonstrate that the repartitioning strategy introduced in \cite{buvoli2022stability} is essential for obtaining stable exponential Parareal methods on stiff non-diffusive problems.

The organization of this paper is as follows. In \cref{sec:introduction-parareal} and \cref{sec:introduction-exponential-integrators} we respectively introduce Parareal and exponential integrators. In \cref{sec:motivating-nls-experiment} we motivate this paper by presenting several numerical experiments involving the one-dimensional nonlinear Schr\"{o}dinger equation. Then, in \cref{sec:linear-stability-convergence} we introduce analytical tools for understanding the convergence and stability properties of a Parareal configuration. Lastly, in \cref{sec:higher-dimensional-numerical-experiments}, we consider two-dimensional hyperbolic and dispersive wave equations and demonstrate that exponential Parareal can achieve reduced time-to-solution compared to serial exponential integrators.  
\section{Parareal Introduction}
\label{sec:introduction-parareal}

In this section we describe the Parareal algorithm \cite{LionsEtAl2001}, present a formula for parallel speedup, and provide a complete table of all Parareal parameters. We begin by supposing that one seeks an accurate, numerical solution to the initial value problem
\begin{align}
	\mathbf{y}'(t) = f(\mathbf{y}(t)), \quad \mathbf{y}(t_0) = \mathbf{y}_0.
	\label{eq:model-ode}
\end{align}
If computational cost can be neglected, then an accurate numerical integrator $\Fprop$, such as a high-order method with small timesteps, should be considered. However, this will not always be practical since the time to run the calculation can become prohibitive. Therefore, we often settle for a less accurate integrator $\Gprop$, such as a low-order method that is run using larger timesteps. Can the situation be improved with access to parallel computational hardware? 

The Parareal method is a parallel iteration that combines a {\em coarse propagator $\Gprop$} with a {\em fine propagator} $\Fprop$, and converges to the solution of $\Fprop$. Provided that the iteration can be efficiently parallelized and that the convergence rate is sufficiently high, then the computational time needed to run Parareal is similar to that of running the coarse propagator. In the following subsections, we explore the algorithm in more detail.

\subsection{Method definition}

Let $\mathcal{F}$ and $\mathcal{G}$ be two one-step methods that are respectively called the coarse and fine propagators; it is assumed that $\Fprop$ is more computationally expensive to apply than $\Gprop$. Next, suppose that we want to approximate \cref{eq:model-ode} at a discrete set of time points using the fine propagator, such that
	\begin{align}
		y_{n+1} = \Fprop(y_n), \quad n = 0, \ldots, N_p,
		\label{eq:fine_integrator_solution}
	\end{align}
	where $y_{n} \approx y(t_n)$. The Parareal algorithm converges to \cref{eq:fine_integrator_solution} by taking a provisional solution, $y^0_n \approx y(t_n)$, that is usually computed by a serial application of the coarse propagator $\Gprop$, and then correcting it via the iteration
	\begin{align}
		\ynp^{k+1} = \Gprop(\ynkp) + \Fprop(\ynk) - \Gprop(\ynk), \quad
		\left\{
			\begin{aligned}
				 n &= 0, \ldots, \Nprocs, \\
				 k &= 0, \ldots, K-1.	
			\end{aligned}
		\right.
		\label{eq:parareal_iteration}
	\end{align}
	In order to run the Parareal iteration, it is necessary to store and iteratively update the solution values along the entire time interval. The key property of the Parareal iteration is that the fine integrator $\Fprop$ can be applied in parallel on $\Nprocs$ processors. To further clarify this point, we show pseudocode for the Parareal iteration \cref{eq:parareal_iteration} in \cref{tab:parareal_pseudocode}.

\begin{table}

	\begin{center}
			
		\begin{tabularx}{0.45\linewidth}{|X|} \hline
		
			{\bf Parareal Pseudocode} \\ \hline

			\cdline{1} {\color{plot_green} \% provisional solution}
			
			\cdline{2} {\bf for} n = 0 : $\Nprocs - 1$
						
			\cdline{3} \hspace{1em} $y_{n+1}^0 = {\color{plot_blue} \mathcal{G}(y_n^0)}$
			
			\cdline{4} {\color{plot_green} \% Parareal iteration}
			
			\cdline{5} {\bf for} k = 0 : K - 1
			
			\cdline{6} \hspace{1em} {\bf parfor} j = 0 : $\Nprocs - 1$  {\color{plot_green} \% parallel loop}
			
			\cdline{7} \hspace{2em} $F_j = {\color{plot_red}\mathcal{F}(y_j^{k})}$
			
			\cdline{8} \hspace{1em} {\bf for} j = 0 : $\Nprocs - 1$
			
			\cdline{9} \hspace{2em} $y_{j+1}^{k+1} = {\color{plot_blue} \mathcal{G}(y_j^{k+1})} + F_j - {\color{plot_blue} \mathcal{G}(y_j^{k})}$
			
			\cdline{10}{\bf return} $y^{K}_{\Nprocs}$ \\ \\ \hline

		\end{tabularx}

	\end{center}
	
	\caption{Pseudocode for the Parareal iteration \cref{eq:parareal_iteration}. The fine integrator (colored in {\color{plot_red} red}) can run in a parallel loop, while the loops containing the coarse propagator (colored in {\color{plot_blue} blue}) are serial. Pseudocode for more efficient pipelined implementations are contained in \cite{Aubanel2011,ruprecht2017shared}.}
	\label{tab:parareal_pseudocode}
	
\end{table}

\subsection{Parallel speedup}
\label{subsec:parallel-speedup}

Parallel speedup is defined as the ratio between the computational time for running a serial algorithm and its parallel equivalent. For Parareal we compute speedup by dividing the computational cost of the sequential fine integrator \cref{eq:fine_integrator_solution} by the computational cost of the Parareal iteration \cref{eq:parareal_iteration} when run using $\Nprocs$ processors. Let the cost for a single step of the fine propagator $\Fprop$ and the coarse propagator $\Gprop$ be $\CostF$ and $\CostG$, respectively. The cost of performing $K$ Parareal iterations is the sum of the cost of the predictor, $\Nprocs  \CostG$, plus the additional cost of each iteration which, neglecting communication, is $K(\CostF + \CostG)$; see \cite{Aubanel2011}. In summary, the serial cost for computing  \cref{eq:fine_integrator_solution} is $C_s = \Nprocs \CostF$ and a parallel cost for \cref{eq:parareal_iteration} is $C_p = \Nprocs  \CostG + K(\CostF + \CostG)$.
If we let $\alpha=\CostG/\CostF$, then the parallel speedup is
\begin{align}
  \Speedup= \frac{\CostS}{\CostP} = \frac{\Nprocs}{\Nprocs \alpha + K (1+ \alpha )}. 
  \label{eq:parareal_speedup}
\end{align}
Lastly we remark the speedup formula will change if one considers more elaborate parallelization strategies such as those presented in \cite{Aubanel2011,BerryEtAl2012,ArteagaEtAl2015,Ruprecht2017_lncs}.

\subsection{Selecting the coarse and fine propagators}

\begin{figure}
	\centering
	\includegraphics[width=0.8\linewidth]{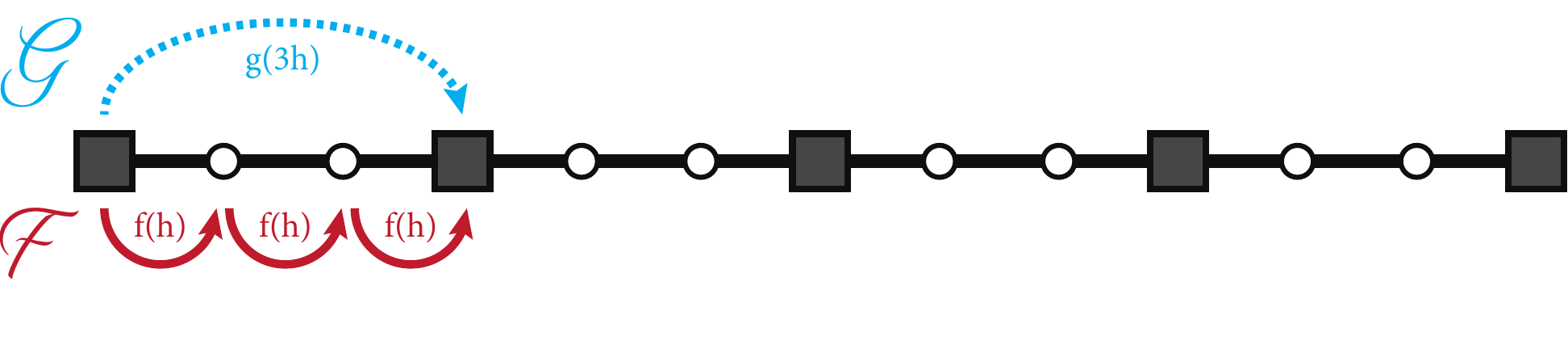}
	\caption{An illustration of the coarse and fine propagators (dashed blue arrow for coarse and solid red arrows for fine), the coarse temporal grid (large, black squares), and fine temporal grid (small, white circles). This illustration depicts the following parameters: the time interval has been divided into 12 total timesteps ($\Nstepstot = 12$), the coarse propagator $\Gprop$ that takes a single step of $g$ ($\NGprop = 1$), the fine propagator $\Fprop$ that takes three steps of $f$ ($\NFprop = 3$), and the resulting Parareal iteration requires $4$ processors ($\Nprocs = 4$).}
	\label{fig:parareal-coarse-fine-grid}
\end{figure}

A user can select any pair of one-step methods to be the coarse and fine propagators. A common approach, which we will use in this work, is to set the coarse propagator $\Gprop$ equal to $\NGprop$ steps of an inexpensive one-step method $g$ and the fine propagator $\Fprop$ equal to  $\NFprop$ steps of an expensive integrator $f$. Both $\Fprop$ and $\Gprop$ must advance the solution by the same amount; therefore, if we let $h$ be the stepsize of $f$, then the stepsize of $g$ must be $ \Dt \NFprop / \NGprop$. If we use the notation $M^\kappa(\eta)$ to denote $\kappa$ steps of a method $M$ run with stepsize $\eta$, then the fine and coarse propagators are
	\begin{align}
		\Fprop = f^\NFprop(\Dt) \quad \text{and} \quad \Gprop = g^\NGprop( \Dt\NFprop  / \NGprop).
		\label{eq:coarse_fine_propagators}
	\end{align}
	In \cref{fig:parareal-coarse-fine-grid} we illustrate the resulting coarse and fine grids for the coarse and fine propagators \cref{eq:coarse_fine_propagators}. Since the fine propagator $\Fprop$ is now $\NFprop$ steps of the method $f$, a Parareal method that converges to the solution of $\Fprop$ applied over $\Nprocs$ steps is also converging to the solution of $f$ applied over total of $\Nstepstot = \NFprop \Nprocs$ steps.  Throughout this work we will frequently characterize Parareal in terms of $(f, \NFprop)$, $(g, \NGprop)$, and $\Nstepstot$ instead of $\Fprop$, $\Gprop$, and $\Nprocs$. 

\subsection{A complete table of parameters}

As we have seen, the Parareal algorithm has a large number of free parameters. In \cref{tab:parareal-parameters} we make a complete list of parameters that are relevant to this work. As we note in the table, the integer variables $\Nprocs$, $\Nstepstot$, and $\NFprop$ are related by the equation $\Nstepstot = \Nprocs \NFprop$ and therefore the user can only select two of these variables with the third being automatically determined.

\begin{table}[h]
		
	{\bf Interdependent User Defined Parameters} \\ (User must select two in such a way that all three variables are integers) \\[0.5em]
	
	\begin{minipage}[b]{0.35\linewidth}
		\begin{tabular}{| c |  l | l |}
			\hline
			Variable & Meaning  \\ 
			\hline
			$\Nprocs$  & Number of processors \\
			$\Nstepstot$  & Total number of RK steps \\ 
			$\NFprop$  & Number of RK steps in $\mathcal{F}$ \\ \hline
		\end{tabular} \\[0.5em]	
	\end{minipage}
	\hspace{1em}
	\begin{minipage}[b]{0.15\linewidth}
		Interdependency
		
		$\Nstepstot = \Nprocs\NFprop$
		\vspace{2.5em}	
	\end{minipage}
	
	{\bf Independent User Defined Parameters} (User must select all) \\[0.5em]
	
	\begin{tabular}{| c |  l |}
		\hline 
		Variable & Meaning \\ \hline
		$f$ & RK method used in $\Fprop$ \\
		$g$ & RK method used in $\Gprop$ \\		
		$\NGprop$ & Number of RK steps in $\Gprop$ \\	
		$K$ & Number of Parareal iterations \\ \hline
	\end{tabular} \\[0.5em]

	{\bf Dependent Parameters} \\[0.5em]
	
	\begin{tabular}{| c |  l | l |}
		\hline
		Variable & Meaning & Definition  \\ \hline
		$\Dt$       & Timestep for serial method &  $\Tfinal/\Nstepstot$ \\  
		$\Gprop$ & Coarse propagator & $\NGprop$ steps of RK method $g$\\  
		$\Fprop$ & Fine propagator & $\NFprop$ steps of RK method $f$\\  
		$\Nsteps$ & Total number of fine steps &  $\Nstepstot$ \\  
		$\CostGstep$  & Cost of $g$ per step  & User defined   \\
		$\CostFstep$  & Cost of $f$ per step   & User defined   \\
		$\CostG$  & Cost of $\Gprop$ per coarse step  & $\NGprop \CostGstep$  \\
		$\CostF$  & Cost of $\Fprop$ per coarse step  & $\NFprop \CostFstep$   \\
		\hline
	\end{tabular}
	\caption{Parareal parameters names and definitions that are relevant to this work.}
	\label{tab:parareal-parameters}
	
\end{table}
\section{Exponential integrators}
\label{sec:introduction-exponential-integrators}

\label{sec:exponential_integrators}

The aim of this work is to study Parareal methods where the coarse and fine propagators are exponential Runge-Kutta methods. In this section we provide an introduction to exponential integrators and discuss their stability properties on non-diffusive equations. Exponential integrators \cite{hochbruck2010exponentialreview} are a class of numerical methods for solving the semilinear initial value problem
\begin{align}
	\mathbf{y}' = \mathbf{Ly} + N(t,\mathbf{y}), \quad \mathbf{y}(t_0) = \mathbf{y}_0.
	\label{eq:semilinear-ode}
\end{align}
In the past two decades, they have proven highly efficient for solving stiff systems and can offer certain advantages over both fully-implicit and linearly-implicit methods \cite{grooms2011IMEXETDCOMP, KassamTrefethen05ETDRK4, loffeld2013comparative, montanelli2016solving,hochbruck1997krylov, hochbruck1998exponential}. The main idea behind exponential integrators is to consider the exact solution to \cref{eq:semilinear-ode}, namely
\begin{align}
  \mathbf{y}(t_0 + \Dt) = e^{\Dt \mathbf{L}}\mathbf{y}_0 + \int^{t_0+\Dt}_{t_0}e^{(t_0 + h -\tau)\mathbf{L}}
  N(\tau,\mathbf{y}(\tau))d\tau,
  \label{eq:exp_derivation_formula}
\end{align}
and replace the nonlinear term $N(\tau,\mathbf{y}(\tau))$ with an explicit polynomial approximation in $\tau$. Such approximations are then used to compute the stages or outputs of exponential integrators families such as linear multistep methods \cite{beylkin1998ELP}, Runge-Kutta methods \cite{cox2002ETDRK4,krogstad2005IF,hochbruckostermann2005ETDRKSTIFFB,luan2014explicit,buvoli2019esdc}, and general linear methods \cite{ostermann2006general,buvoli2021epbm}.

A polynomial approximation of the nonlinear term implies that the formula of all exponential integrators can be expressed in terms of the exponential functions
\begin{align}
	\varphi_0(h\mathbf{L}) = e^{h\mathbf{L}} \quad \text{and} \quad 
	\varphi_j(h\mathbf{L}) = \int_{0}^{1} e^{(1-s)h\mathbf{L}}s^j ds \quad (j \ge 1);
 	\label{eq:phi-functions}
\end{align}
specifically, replacing  $N(\tau, y(\tau))$ in  \cref{eq:exp_derivation_formula} with $\sum_{i} \mathbf{c}_i \tau^i$ and defining $s=h\tau+t_0$ leads to a linear combination of functions $\varphi_j(hL)$. Therefore, at each timestep an exponential integrator requires matrix-vector products with the $\varphi$-functions of the linear operator $\mathbf{L}$. For many problems this can be done efficiently using a number of different algorithms \cite{ashi2009comparison,caliari2014comparison,higham2020catalogue}, including those based on squaring methods \cite{koikari2007error,al2011computing,al2010new}, contour integration \cite{KassamTrefethen05ETDRK4,trefethen2007,schreiber2019exponential}, Krylov-subspaces \cite{hochbruck1997krylov,hochbruck1998exponential,NiesenWright2011Krylov,NiesenWright2012Krylov,GAUDREAULT2018236}, and parallel rational approximations \cite{haut2015high,schreiber2019exponential,schreiber2019parallel}.

\subsection{Exponential Runge-Kutta methods}
\label{subsec:erk}

Exponential Runge-Kutta (ERK) methods are one-step methods that approximate the solution to \cref{eq:model-ode} by taking a linear combination of stage values at each timestep . The simplest ERK integrator is the exponential Euler method that is obtained by replacing $N(\tau,\mathbf{y}(\tau))$ in \cref{eq:exp_derivation_formula} with the constant approximation $N(t_n, \mathbf{y}_n)$, yielding
\begin{align}
	\mathbf{y}_{n+1} = \varphi_0(h\mathbf{L}) \mathbf{y}_n + \varphi_1(h \mathbf{L})hN(t_n, \mathbf{y}_n).
	\label{eq:exponential-euler}
\end{align}
More generally, an $s$-stage ERK method is
\begin{align}
	Y_i &= \varphi_0(h c_j\mathbf{L}) \mathbf{y}_n + \sum_{j=1}^{i-1} a_{ij}(h\mathbf{L}) N(c_j, Y_j), \quad i = 1, \ldots, s, \label{eq:erk_stages} \\
	\mathbf{y}_{n+1} &= \varphi_0(h \mathbf{L}) \mathbf{y}_n + \sum_{j=1}^{s} b_{j}(h \mathbf{L}) N(c_j, Y_j) \label{eq:erk_output}
\end{align}
where $a_{ij}(h\mathbf{L})$ and $b_j(h\mathbf{L})$ are functions that include linear combinations or products of the $\varphi$-functions \cref{eq:phi-functions}. By applying the identity $\varphi_0(h\mathbf{L})\mathbf{y}_n = \mathbf{y}_{n} + \varphi_1(h\mathbf{L}) \mathbf{L y}_n$ one can rewrite the equations \cref{eq:erk_stages,eq:erk_output} in terms of $\varphi_j(hK)$ for $j\ge 1$; this can be advantageous both for method analysis and implementation. Lastly, like classical RK methods, ERK methods can be represented using the Butcher tableau
\begin{center}
	\begin{tabular}{l|llll}
			$c_1$ 		& $0$ 			& 				& 				& 		\\ 
			$c_2$ 		& $a_{21}$ 		& 	$0$	 		& 				& 		\\
			$\vdots$	& $\vdots$ 		& 	$\ddots$ 	& $\ddots$ 		& 		\\
			$c_s$ 		& $a_{s,1}$ 		& 	$\hdots$ 	& $a_{s,s-1}$ 	& $0$ 	\\ \hline
						& $b_1$ 			& 	$\hdots$ 	& $b_{s-1}$ 		& $b_s$ 
	\end{tabular}
\end{center}
where the coefficients $a_{ij}$ and $b_j$ are now matrix functions of the linear operator $h\mathbf{L}$.

In this work, we will consider ERK methods of orders one to four from \cite{cox2002ETDRK4,krogstad2005IF}. We name these methods ERK1, ERK2, ERK3, and EKR4, and list their tableaux in \cref{app:coefficients}. 

\subsection{Stability and repartitioning for non-diffusive equations}

Since exponential integrators treat the linear operator $\mathbf{L}$ exactly, we would expect that they offer significantly improved stability properties compared to explicit integrators. While this is true for diffusive operators, the situation is more nuanced when $\mathbf{L}$ has purely imaginary eigenvalues \cite{buvoli2022stability,crouseilles2020exponential}. In particular, both exponential and explicit integrators have similarly sized stability regions, but the magnitude of the instabilities is often very small for exponential integrators. Therefore, unlike explicit methods, exponential integrators can still produce usable solutions on stiff non-diffusive equations so long as the total number of timesteps is not overly large \cite{buvoli2022stability}. 

In \cite{buvoli2022stability} we proposed a strategy that stabilizes exponential integrators by repartitioning the right-hand-side of \cref{eq:semilinear-ode} using perturbed linear and nonlinear operators $\widehat{\mathbf{L}}$ and $\widehat{N}$. This enables long-time simulations with exponential integrators and also removes instabilities when the underlying equation focuses energy into unstable modes. The perturbed operators are formed by respectively adding and subtracting a diffusive operator $\mathbf{D}$ such that
\begin{align}
	\widehat{\mathbf{L}} = \mathbf{L} + \epsilon	 \mathbf{D}
	\quad \text{and} \quad
	\widehat{N}(t,\mathbf{y}) = N(t,\mathbf{y}) - \epsilon \mathbf{D}.
	\label{eq:repartitioned-operators}
\end{align}
In short, we add damping to the linear operator $\widehat{\mathbf{L}}$ and excitation to the nonlinear operator $\widehat{N}(t,\mathbf{y})$. The differential equation \cref{eq:semilinear-ode} can then be written in terms of the perturbed operators as
	\begin{align}
		\mathbf{y}' &= \widehat{\mathbf{L}} \mathbf{y} + \widehat{N}(t,\mathbf{y}).
		\label{eq:semilinear-ode-epartitioned}
	\end{align}
Therefore, an exponential integrator that solves the repartitioned equation \cref{eq:semilinear-ode-epartitioned} is simultaneously solving \cref{eq:semilinear-ode}. However, instead of treating $\mathbf{L}$ exactly and approximating $N(t,\mathbf{y})$, a repartitioned integrator treats $\widehat{\mathbf{L}}$ exactly and approximates $\widehat{N}(t,\mathbf{y})$. The advantage of repartitioning is that the exponential integrator now possesses a large stability region for a continuous range of small $\epsilon$ values \cite{buvoli2022stability}. If $\mathbf{L}$ is diagonalizable, such that $\mathbf{L} = \mathbf{U \Lambda U}^{-1}$, and we select
\begin{align}
	\mathbf{D} = - \mathbf{U |\Lambda| U}^{-1}	\quad \text{and} \quad \epsilon = \frac{1}{\tan(\pi/2 - \rho)} \quad \text{for} \quad \rho \in (0, \pi/2),
	\label{eq:rho-repartitioning}
\end{align}
then we rotate all the eigenvalues of a non-diffusive linear operator $\mathbf{L}$ by $\rho$ degrees into the left-half plane. In other words, the eigenvalues of $\widehat{\mathbf{L}}$ all lie on the wedge $re^{i\theta}$ for $r\ge 0$ and $\theta \in \{\pi/2 + \rho, 3\pi/2 - \rho\}$. Many other choices for $\mathbf{D}$ are possible (e.g. low-order, even spatial derivatives), however in \cite{buvoli2022stability} we proposed \cref{eq:rho-repartitioning} because of its convenience when analyzing the stability effects of repartitioning.

In summary, exponential integrators exhibit mild instabilities for stiff non-diffusive equations that can be eliminated through re-partitioning. In \cref{sec:motivating-nls-experiment,sec:linear-stability-convergence} we will show that these instabilities are greatly exacerbated by the Parareal iteration, and that repartitioning is essential for obtaining stable, exponential Parareal methods for solving stiff non-diffusive equations.
\section{Motivating numerical experiments}
\label{sec:motivating-nls-experiment}

In this section we present three numerical experiments that highlight key properties of exponential Parareal integrators applied to non-diffusive equations. Later, in \cref{sec:linear-stability-convergence} we will see how linear stability analysis and linear convergence analysis can be used to more rigorously quantify our results. All three numerical experiments involve the one-dimensional nonlinear Schr\"{o}dinger (NLS) equation 
\begin{align}
	& iu_t + u_{xx} + 2|u|^2u = 0	\label{eq:nls}	
\end{align}
on the domain $x\in[-4\pi, 4\pi]$ with periodic boundary conditions. We discretize the equation in space using a 1024 point Fourier spectral method that is dealiased using the classical $3/2$ rule.  The equation is then integrated in Fourier space where the derivative operators are diagonal. This results in the semilinear equation \cref{eq:semilinear-ode} with
\begin{align}
	\mathbf{L} = \text{diag}(-i \mathbf{k}^2) \quad \text{and} \quad N(t,\mathbf{y}) = 2i\mathcal{F}(\mathcal{F}^{-1}(\mathbf{y}) \text{ * } \text{abs}(\mathcal{F}^{-1}(\mathbf{y}))),
	\label{eq:nls-system-discretized}
\end{align}
where $\mathbf{k}$ is a vector of Fourier wavenumbers, $\mathcal{F}$ denotes the discrete Fourier transform, and $*$ is an elementwise multiply (i.e. Hadamard product). To ensure the classical stability of exponential integrators we also consider the repartitioning \cref{eq:repartitioned-operators,eq:semilinear-ode-epartitioned} where the diffusive operator $\mathbf{D}$ and $\epsilon$ are selected according to \cref{eq:rho-repartitioning} such that
\begin{align}
	\mathbf{D} = \text{diag}(-\mathbf{k}^2),
	\quad \epsilon = \frac{1}{\tan(\pi/2 - \rho)}
	\quad \text{and} \quad
	\rho = \frac{\pi}{128}.
	\label{eq:repartitioning-nls}
\end{align}

In addition to investigating stability and convergence, we also compare the theoretical speedup of the Parareal iteration to its real-world performance on a distributed memory system. To do this, we implemented the exponential Parareal method as part of the open source package LibPFASST\footnote{\url{https://github.com/libpfasst/LibPFASST}} and performed the numerical experiment on the Cray XC40 {\em Cori} at the National Energy Research Scientific Computing Center.

Our first experiment uses the initial condition
\begin{align}
	u(x,t=0) = 1 + \tfrac{1}{100}	 \cos(x / 4), \label{eq:nls-low-frequency-ic}	
\end{align}
integrated out to time $t=14$. Repartitioning is not required for serial ERK methods on this short time-scale, and the convergence curves for classical and repartitioned exponential integrators look identical (see convergence plots in \cref{app:nls-serial-convergence}). We now consider two Parareal methods: one with classical ERK integrators, and the other with repartitioned ERK integrators. To obtain a high-accuracy solution, we select ERK4 as the fine integrator $f$, and $\Nstepstot = 2^{16}$ as the total number of fine steps; from the serial ERK convergence diagrams in \cref{fig:nls-serial-convergence-smooth} we see that a fully converged Parareal method will yield the solution with an error of $3\times 10^{-9}$. Next we must select a coarse integrator $g$ that is stable at large stepsizes. Though it may seem tempting to select ERK1 because it is the least expensive method per timestep, its poor stability leads us to choose the more stable ERK3 method. Lastly, we select $\NGprop = 1$ and $\NFprop = 32$; this implies that the coarse propagator $\Gprop$ consists of a single step of $g$, while the fine propagator $\Fprop$ consists of $32$ steps of $f$. With these parameters the serial coarse propagator has an accuracy of $2\times 10^{-1}$, which is approximately eight orders of magnitude less than the serial fine integrator (see the black crosses in \cref{fig:nls-serial-convergence-smooth}). A complete list of the Parareal parameters used for this experiment is contained in \cref{tab:nls-parareal-configuration}.  

\begin{table}
	\renewcommand*{\arraystretch}{1.25}
	\begin{tabular}{ll|ll|ll|ll}
		$\Nprocs$ 		& $2048$ 	& $\NFprop$ 		& $32$	& $f$ 	& ERK4   & $K$  & $1, \ldots, 6$\\ \hline
		$\Nstepstot$	& $2^{16}$ 	& $\NGprop$	 	& $1$ 	& $g$	& ERK3   & 
	\end{tabular}
	\caption{Parareal parameters used for the nonlinear Schrodinger equation \cref{eq:nls} with initial conditions \cref{eq:nls-low-frequency-ic,eq:nls-high-frequency-ic}.}
	\label{tab:nls-parareal-configuration}
\end{table}

\begin{figure}

	\begin{center}

		\begin{minipage}{0.32\textwidth}
			\centering
			\hspace{1em} {{\bf (a)} \footnotesize Error versus Iteration}
		\end{minipage}
		\begin{minipage}{0.32\textwidth}
			\centering
			{{\bf (b)} \footnotesize Parallel Speedup \cref{eq:parareal_speedup} versus Iteration}
		\end{minipage}
		\begin{minipage}{0.32\textwidth}
			\centering
			{{\bf (c)} \footnotesize NLS Solution}
		\end{minipage}

		\begin{minipage}{0.32\textwidth}
			\includegraphics[width=\linewidth]{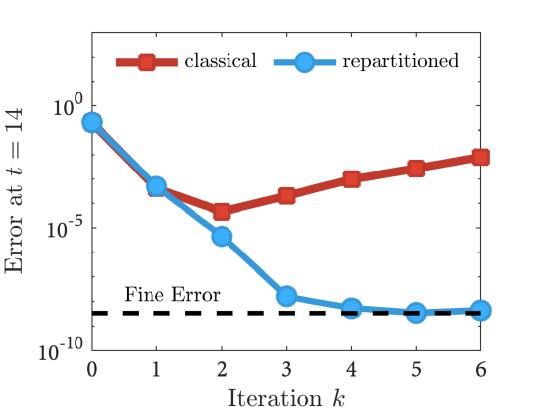}
		\end{minipage}
		\begin{minipage}{0.32\textwidth}
			\includegraphics[width=\textwidth]{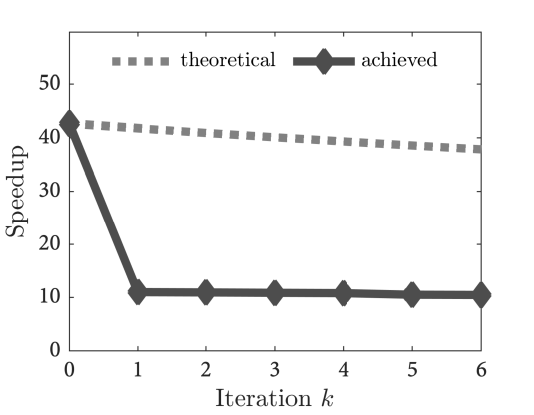}
		\end{minipage}
		\begin{minipage}{0.32\textwidth}
			\includegraphics[width=\linewidth,trim={0 0 0 1.5cm},clip]{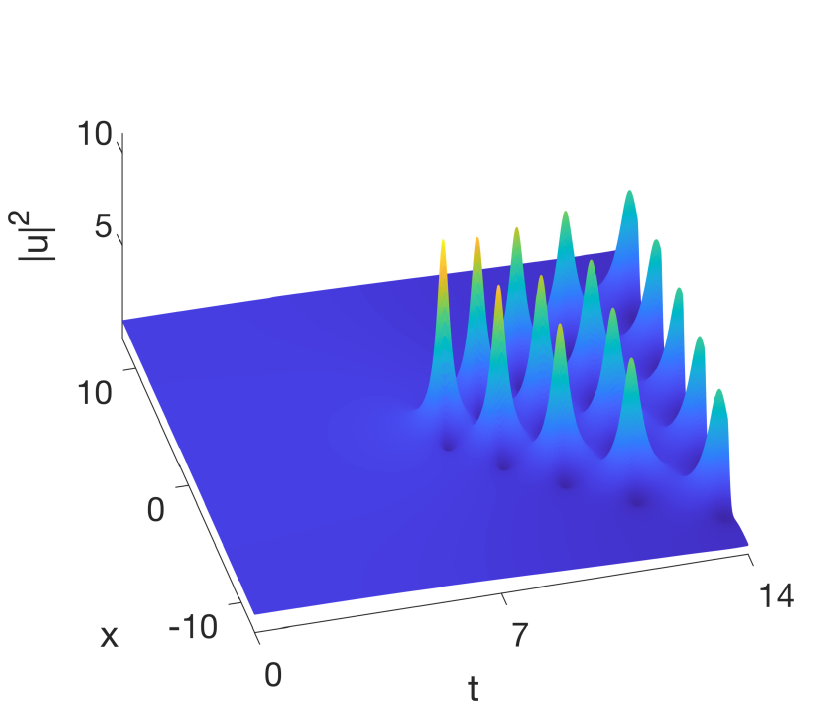}	
		\end{minipage}
		
	\end{center}
		
	\caption{Error and speedup of the Parareal configuration from \cref{tab:nls-parareal-configuration} applied to the NLS equation \cref{eq:nls} with the smooth initial condition \cref{eq:nls-low-frequency-ic}; error is measured with respect to our numerically computed reference solution, not the fine integrator solution. Subfigure {\bf (a)} shows the error at $\Tfinal=14$ as a function of the Parareal iteration $k$. Line color is used to distinguish classical exponential integrators from repartitioned exponential integrators. When $k=0$ the Parareal method is equivalent to running the coarse integrator with $\Nstepstot/\NFprop$ steps. Subfigure {\bf (b)} shows the parallel speedup as a function of the iteration $k$; this compares the running time of the Parareal iteration to that of taking $\Nstepstot$ serial steps with the fine integrator $f$. Note that speedup is identical for both classical and repartitioned Parareal. Lastly, figure {\bf 2(c)} shows the magnitude squared NLS solution $|u(x,t)|^2$. 
	}
	\label{fig:nls-experiment-smooth}
	
\end{figure}

 The computation was distributed on 64, 32-core Intel Haswell nodes that provided a total of 2048 compute cores. In \cref{fig:nls-experiment-smooth} we show how the error of the solution obtained by the Parareal iteration evolves as a function of the iteration number; error is measured with respect to our numerically computed reference solution, not the fine integrator solution. We also show the corresponding theoretical and achieved parallel speedup, along with a space-time plot of the NLS solution. The results demonstrate that repartitioning is essential for obtaining a convergent Parareal iteration; moreover, by using repartitioned ERK methods,  Parareal is able to obtain a high-accuracy solution up to 38.5 times faster than a serial ERK4 method. In contrast unmodified exponential integrators lead to a divergent Parareal iteration whose error increases monotonically for iteration number $k$ greater than two.   
 
Our timing results also reveal a practical challenge that can occur when applying a PinT method on a distributed memory system. Using the Parareal configuration from \cref{tab:nls-parareal-configuration}, we were only able to achieve a speedup factor of 10.4; approximately one quarter of the theoretical speedup factor of 38.5 predicted by \cref{eq:parareal_speedup}. This difference is due to unaccounted communication overhead. In fact, \cref{eq:parareal_speedup} will only be accurate if the time required to compute a single step of the propagator $\Gprop$ is significantly greater than the time for transferring the solution vector between two nodes. Although this does not hold true for the one-dimensional NLS equation, the exponential Parareal method is neverthless able to provide a high-accuracy solution an order of magnitude faster than the serial ERK4 method. Moreover, in \cref{sec:higher-dimensional-numerical-experiments} we will see that theoretical speedup very accurately predicts achievable speedup on more computationally expensive two-dimensional problems.
   
\begin{figure}

	\begin{center}
	
		\begin{minipage}[t]{0.35\linewidth}
			\centering
			\hspace{0.5em} {{\bf (a)} \footnotesize Error versus Iteration}
			
			\includegraphics[width=\linewidth]{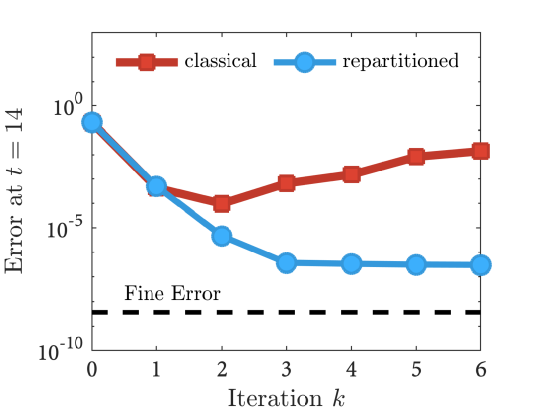}
		\end{minipage}		
		\begin{minipage}[t]{0.6\linewidth}
			\begin{minipage}[t]{0.5\linewidth}
				\centering
				\hspace{0.5em} {{\bf (b)} \footnotesize Solution at $t=14$}
			
				\includegraphics[width=\linewidth]{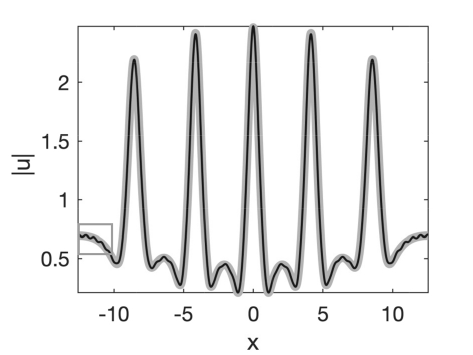}
			\end{minipage}
			\begin{minipage}[t]{0.5\linewidth}
				\centering
				\hspace{1em} {{\bf (c)} \footnotesize  Magnified Solution at $t=14$}
				
				\includegraphics[width=\linewidth]{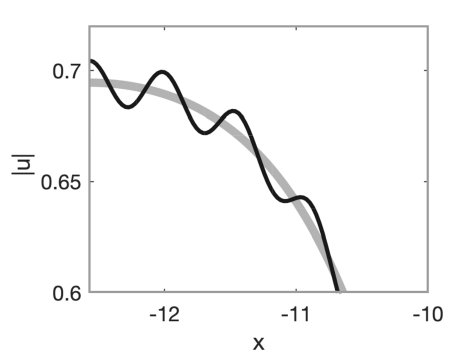}
			\end{minipage}
			
			\vspace{0.5em}
			\begin{minipage}{1\linewidth}
				\begin{center}				
					
					{\footnotesize \textcolor{nls_ic_plot_gray}{\hdashrule[0.2ex]{2em}{3pt}{}} initial condition is \cref{eq:nls-low-frequency-ic} \hspace{2em}}
					{\footnotesize \textcolor{black}{\hdashrule[0.2ex]{2em}{1pt}{}} initial condition is \cref{eq:nls-high-frequency-ic}}
				\end{center}
			\end{minipage}
		\end{minipage}

	\end{center}

	\caption{ Subfigure {\bf(a)} shows the error of the Parareal configuration from \cref{tab:nls-parareal-configuration} applied to the NLS equation \cref{eq:nls} with the oscillatory initial condition \cref{eq:nls-high-frequency-ic}.  Subfigure {\bf (b)} compares the NLS solution at time $t=14$ for the two initial conditions \cref{eq:nls-low-frequency-ic} and \cref{eq:nls-high-frequency-ic} that are drawn using a thick gray line and thin black line, respectively. The region enclosed by a blue square is magnified in  subfigure {\bf (c)} to highlight the small amplitude oscillation that arises from the initial condition \cref{eq:nls-high-frequency-ic}.}
	
	\label{fig:nls-experiment-oscillatory}
		
\end{figure}

For our second experiment we consider a modified initial condition that contains a high-frequency component, namely
\begin{align}
	u(x,t=0) = 1 + \tfrac{1}{100} \left[	 \cos(x / 4) + \cos(45x / 4) \right]. \label{eq:nls-high-frequency-ic}		
\end{align}
The newly added perturbation does not fundamentally change the solution, but rather introduces a low-amplitude oscillation that persists throughout the temporal integration window; see  \cref{fig:nls-experiment-oscillatory}(b)-(c). Moreover, the high-frequency mode does not make the computation more challenging for serial ERK methods, as evidenced by the convergence and efficiency plots that are nearly identical to those generated using the initial condition \cref{eq:nls-low-frequency-ic}; compare \cref{fig:nls-serial-convergence-smooth,fig:nls-serial-convergence-oscillatory}. We again consider the 
Parareal method from \cref{tab:nls-parareal-configuration} with either classical or repartitioned ERK integrators. In \cref{fig:nls-experiment-oscillatory}(a) we show error as a function of the iteration number $k$, and see that repartitioning is again required to prevent instabilities. However, the high-frequency component has now prevented the repartitioned Parareal method from fully converging; the method  remains stable as $k$ increases, however the error does not improve beyond $3\times 10^{-7}$. This is our first indication that highly-oscillatory solutions will cause convergence problems for Parareal. \Cref{app:nls-experiment-oscillatory-extended-k} contains a complementary error versus iteration graph for $\NIters$ up to 160. It further reveals that: (1) Parareal with classical ERK methods becomes completely unstable after 9 iterations, and (2) repartitioned Parareal requires additional iterations to converge, thus eliminating the hope for significantly reduced time-to-solution.

\begin{figure}

	\begin{center}
	
		\begin{minipage}{0.4\textwidth}
			\centering
			\hspace{0em} { \footnotesize {\bf (a)} Error versus Iteration}
		\end{minipage}
		\begin{minipage}{0.3\textwidth}
			\centering
			\hspace{4em} {\footnotesize {\bf (b)}  Parallel Speedup \cref{eq:parareal_speedup}}
		\end{minipage}
	
		\begin{minipage}{0.4\textwidth}
			\centering
			\includegraphics[height=0.75\linewidth]{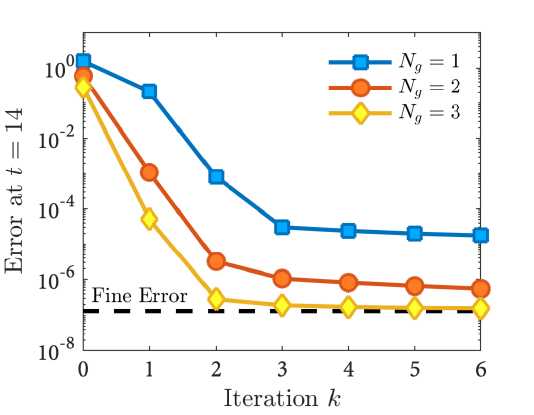}
		\end{minipage}
		\begin{minipage}{0.3\textwidth}
			\centering
			\includegraphics[height=\linewidth]{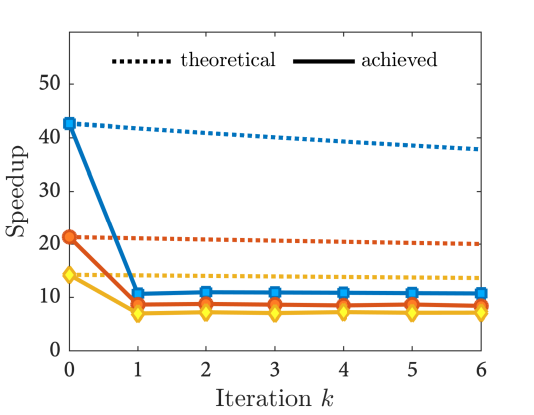}
		\end{minipage}
		
	\end{center}

	\caption{Parareal configuration with $\NGprop\in\{1,2,3\}$ (all other parameters are in \cref{tab:nls-parareal-configuration})  applied to the nonlinear Schr\"{o}dinger equation \cref{eq:nls} with initial conditions \cref{eq:nls-full-spectrum-ic}. Subfigure {\bf (a)} shows how the error at $y(\Tfinal=14)$ evolves as a function of the Parareal iteration $k$. Subfigure {(\bf b)} shows how much faster the Parareal algorithms are compared to taking $\Nstepstot$ serial steps with the fine integrator $f$. Increasing $\NGprop$ improves convergence but also decreases parallel speedup. We again see significant decrease in speedup due to communication overheads.
	}
	
	\label{fig:nls-ng-experiment}
	
\end{figure}

Thus far, we have seen that repartitioning is essential for preventing instabilities in the exponential Parareal iteration. Therefore, we will no longer consider Parareal with classically partitioned ERK methods in this section. However, we have also seen that repartitioning does not guarantee rapid convergence. In \cite{buvoli2021imexparareal,Ruprecht2018} it was shown that Parareal convergence on non-diffusive problems improves when the coarse propagator $\Gprop$ more closely approximates the fine propagator $\Fprop$. In our final motivating  experiment, we will explore this phenomenon using an even more challenging initial condition that contains $45$ spatial modes, namely
\begin{align}
	u(x,t=0) = 1 + \tfrac{1}{100}\ \sum_{k=1}^{45} \cos(kx/4).
	\label{eq:nls-full-spectrum-ic}
\end{align}
The NLS solution is now full of high-frequency information (see \cref{fig:nls-solution-visualization}) that causes even serial ERK integrators to achieve slightly diminished accuracy for the same number of steps; compare \cref{fig:nls-serial-convergence-full-spectrum} to \cref{fig:nls-serial-convergence-oscillatory}. Based on the previous experiment we expect that the high-frequency oscillations will prevent the Parareal configuration in \cref{tab:nls-parareal-configuration} from rapidly converging to the fine solution. We therefore change the accuracy of the coarse propagator $\Gprop$ by considering two additional choices for the number of coarse steps, namely $\NGprop \in \{1, 2, 3\}$. By increasing $\NGprop$ (the number of steps of $g$ in $\Gprop$) we can exchange parallel speedup for a more accurate coarse integrator. In \cref{fig:nls-ng-experiment} we show error and parallel speedup as a function of the iteration number $k$ for these three Parareal configurations. As expected, the Parareal configuration with $\NGprop=1$ does not converge to the fine solution within six iterations. However, as $\NGprop$ increases we see that convergence properties improve significantly at the cost of decreased parallel speedup. 

In summary, repartitioning is required to prevent instabilities and the accuracy of the coarse integrator must be increased if we want to resolve high-frequency modes. In the next section we will use linear stability analysis and linear convergence analysis to more carefully quantify these statements.

\section{Linear stability and convergence analysis}
\label{sec:linear-stability-convergence}

In this section we study the linear stability and convergence properties of exponential Parareal and provide a mathematical foundation for understanding the numerical experiments from \cref{sec:motivating-nls-experiment}. For classical Parareal methods, there are many existing works studying stability and convergence \cite{Bal2005,StaffRonquist2005,GanderVandewalle2007_SISC,Ruprecht2018} including several that develop rigorous mathematical convergence bounds for diffusive problems \cite{Gander2008,Southworth2019,Southworth2021TightTwoLevel}. Our analysis is based on the partitioned Dahlquist equation 
	\begin{align}
		y' = \lambda_1 y + \lambda_2 y, \quad y(0) = 1,
		\label{eq:partitioned-dahlquist}
	\end{align}
	and follows closely with our previous works \cite{buvoli2021imexparareal,buvoli2022stability} that respectively studied implicit-explicit Parareal and repartitioned exponential integrators. It is important to note that this equation is a considerable simplification of the nonlinear semilinear equation \cref{eq:semilinear-ode} and represents a scenario for which the linear and nonlinear term can be simultaneously diagonalized. Nevertheless, we will see that linear analysis accurately predicts the phenomena observed in \cref{sec:motivating-nls-experiment}, along with the higher dimensional problems in \cref{sec:higher-dimensional-numerical-experiments}. However, for general nonlinear systems there is no guarantee that the partitioned Dahlquist will always be an accurate model. 
	
Any one-step exponential integrator, including Parareal, applied to \cref{eq:partitioned-dahlquist} reduces to an iteration of the form 
\begin{align}
	y_{n+1} = R(z_1,z_2) y_n	 \quad \text{where} \quad z_1 = h\lambda_1,~ z_2 = h \lambda_2,
	\label{eq:one-step-dahlquist-iteration}
\end{align}
and $h$ is the method's stepsize. Consequently, \cref{eq:partitioned-dahlquist} is commonly used to study the stability of both exponential and implicit-explicit methods \cite{ascher1995implicit,cox2002ETDRK4,krogstad2005IF,izzo2017highly,buvoli2022stability};  in the case of exponential integrators the term $\lambda_1y$ is exponentiated while the term $\lambda_2y$ is treated explicitly. It should also be noted that \cref{eq:semilinear-ode} reduces to a system of decoupled, partitioned Dahlquist equations when the linear and nonlinear operator can be simultaneously diagonalized. Since we are only considering non-diffusive equations, we assume that $\lambda_1$ and $\lambda_2$ are purely imaginary. The following table summarizes the relevant equations.
\begin{center}
	\renewcommand*{\arraystretch}{1.5}
	\begin{tabular}{l|ll}
		Nonlinear system 			& $\mathbf{y}' = \mathbf{Ly} + N(\mathbf{y})$ 	&  $\text{eig}(\mathbf{L}), \text{eig}(\mathbf{\frac{\partial N}{\partial y}}) \in i \mathbb{R}$\\
		Partitioned Dahlquist 		& $y' = \lambda_1 y + \lambda_2 y$ 				&  $\lambda_1, \lambda_2 \in i\mathbb{R}$. \\
	\end{tabular}
\end{center}

To estimate stability and convergence properties for a specific nonlinear system, we consider a family of partitioned Dahlquist equations with continuous $\lambda_1$, $\lambda_2$ values that respectively enclose the spectrums of the linear operator $\mathbf{L}$ and the Jacobian of the nonlinear operator $\tfrac{\partial N}{\partial \mathbf{y}}$. This rectangular parameter space in the scaled coordinates $z_1$, $z_2$ is
	\begin{align}
		Z(h) = \left\{ z_1 \in h[-\bar{\lambda}_1,~ \bar{\lambda}_1],~ z_2 \in h[-\bar{\lambda}_2,~ \bar{\lambda}_2] \right\} \quad \bar{\lambda}_1 = i\rho(\mathbf{L}), \quad  \bar{\lambda}_2 = i\max_{t\in[t_0, \Tfinal]} \rho \left(\tfrac{\partial N}{\partial \mathbf{y}}(\mathbf{y}(t))\right)
		\label{eq:dahlquist-parameter-region}
	\end{align}
	where $\rho(\cdot)$ returns the spectral radius and $h$ is the stepsize required by the fine integrator to achieve a desired error tolerance. Ideally, we would like a Parareal configuration to be stable and rapidly convergent for any $(z_1, z_2) \in Z(h)$.  
	
Finally, due to the limited stability of exponential integrators on non-diffusive equations, we must consider repartitioning. If one applies the repartitioning \cref{eq:repartitioned-operators,eq:semilinear-ode-epartitioned,eq:rho-repartitioning}, the equations from the previous table have the following analogs.
\begin{center}
	\renewcommand*{\arraystretch}{1.5}
	\begin{tabular}{l|ll}
		Repartitioned nonlinear system 	&
		$\mathbf{u}' = \underbrace{\left(\mathbf{L} + \epsilon\mathbf{D}\right)}_{\hat{\mathbf{L}}}\mathbf{u} + \underbrace{(N(\mathbf{u}) - \epsilon \mathbf{Du})}_{\hat{N}(\mathbf{u})}$ &
		$\mathbf{L} = \mathbf{U \Lambda U^{-1}}$, $\mathbf{D} = -\mathbf{U|\Lambda|U^{-1}}$ \\
		Repartitioned Dahlquist 			& $y' = \underbrace{(\lambda_1 - \epsilon|\lambda_1|)}_{\hat{\lambda}_1} y + \underbrace{(\lambda_2 + \epsilon|\lambda_1|)}_{\hat{\lambda}_2}y$ \\
	\end{tabular}
\end{center}
Since repartitioning preserves linearity, the iteration \cref{eq:one-step-dahlquist-iteration} for a repartitioned integrator simply becomes
\begin{align}
	y_{n+1} = R(z_1 + \epsilon|z_1|, z_2 - \epsilon |z_1|) y_n.
	\label{eq:one-step-dahlquist-iteration-repartioned}
\end{align}

The remainder of this section is organized as follows. In \cref{subsec:nls-operators} we briefly quantify the parameter ranges that are pertinent for the discretized nonlinear Schr\"{o}dinger equation from \cref{sec:motivating-nls-experiment}. \Cref{subsec:parareal-linear-formulas} then contains simplified formulas for the Parareal method on the partitioned Dahlquist equation. In \cref{subsec:parareal-linear-stability,subsec:parareal-linear-convergence} we use linear analysis to study the stability and convergence properties of the exponential Parareal iteration. This allows us to quantify the stability effects of repartitioning, and to understand why high-frequency oscillations cause convergence problems for Parareal. In \cref{subsec:parareal-nls-validation} we then compare the predictions of linear analysis against the results of our nonlinear numerical experiments. In \cref{subsec:convergence-additional-configuration} we briefly analyze how certain Parareal parameters affect convergence. Lastly, we conclude with \cref{subsec:convergence-implications-pde} where we discuss the implications of convergence analysis for the solution of partial differential equations.

\subsection{Spectral radius of the nonlinear Schr\"{o}dinger operators}
\label{subsec:nls-operators}

To analyze the numerical experiments from \cref{sec:motivating-nls-experiment}, we first determine the parameters of the Dahlquist equation that most closely approximate the discretized nonlinear Schr\"{o}dinger equation \cref{eq:nls-system-discretized}. We proceed by bounding the spectral radius of the linear and nonlinear operators to estimate the rectangular parameter space $Z(h)$ defined in \cref{eq:dahlquist-parameter-region}:
\begin{itemize}
	\item {\em Linear operator}. The linear operator $\mathbf{L}$ for the discretized nonlinear Schr\"{o}dinger equation with an even number of spatial grid points $N_x$ is 
	\begin{align}
		\mathbf{L} = \text{diag}(-i\mathbf{k}^2) \quad \text{for} \quad \mathbf{k} = \tfrac{1}{4} \text{[0:$N_x/2 - 1$, $-N_x/2$:$-1$]}^T.
	\end{align}
	Using $N_x=1024$ and applying dealiasing we have $\bar{\lambda}_1 = \rho(\mathbf{L}) = (341/4)^2$; dealiasing removes the top one-third of the highest frequency modes so that only modes $-341,\ldots,341$ remain.
	
	\item {\em Nonlinear operator Jacobian}. Obtaining the exact spectral radius for the nonlinear Jacobian $\frac{\partial N}{\partial u}$ is more involved. Instead, we estimate its magnitude by assuming there is no coupling between Fourier modes. The continuous nonlinear operator in physical space is $2i|u|^2u$, which, when applied to a single mode $u(x) = a_k e^{ikx}$, leads to $2|a_k|^2 a_k e^{ikx}$. Ignoring mode coupling, the discretized nonlinearity in Fourier space acts like the diagonal operator diag($2i|{\mathbf{u}}|^2$). In each of the experiments from \cref{sec:motivating-nls-experiment}, the elements of $\mathbf{u}$ are all bounded above by 1.0001 throughout the temporal domain, so we estimate that $\bar{\lambda}_2 = \max_{t\in[0,14]} \rho\left(\frac{\partial N}{\partial \mathbf{y}}(\mathbf{y}(t))\right) \approx 2$. 
\end{itemize}
Lastly, all the experiments from \cref{sec:motivating-nls-experiment} use a fine stepsize of $h=14/2^{16}$. Therefore, using \cref{eq:dahlquist-parameter-region}, the nonlinear Schr\"{o}dinger equation can be approximately analyzed using the family of Dahlquist equations with scaled parameters
\begin{align}
	z_1 \in i [0, 1.6]
	\quad \text{and} \quad
	z_2 \in i  [-4.3, 4.3] \times 10^{-4}.
	\label{eq:nls-z1-z2-ranges}
\end{align}
To avoid imaginary numbers, it is convenient to consider the real-valued dimensions of this parameter region, namely
\begin{align}
	r_1 \in [0, 1.6]
	\quad \text{and} \quad
	r_2 \in [-4.3, 4.3] \times 10^{-4}.
	\label{eq:nls-r1-r2-ranges}
\end{align}
We will frequently refer back to these numbers as we study the stability and convergence properties of Parareal on the nonlinear Schr\"{o}dinger equation.

\subsection{Parareal for the partitioned Dahlquist equation}
\label{subsec:parareal-linear-formulas}

We now present several formulas that describe the Parareal iteration \cref{eq:parareal_iteration} applied to the  Dahlquist equation \cref{eq:partitioned-dahlquist}; these formulas were originally developed in \cite{Ruprecht2018} for unpartitioned linear problems. We begin by considering the coarse and fine propagators ($\mathcal{G}$, $\mathcal{F}$) and their underlying integrators ($g$, $f$). When applied to \cref{eq:partitioned-dahlquist} these methods reduce to the scalar iterations
\begin{align}
	\renewcommand*{\arraystretch}{1.5}
	\begin{array}{rllrl}
		\text{$g$:} & y_{n+1} = R_{g}(z_1, z_2)y_n & \hspace{3em} \vphantom{.} 	& \text{$\mathcal{G}$:} & y_{n+1} = R_\mathcal{G}(z_1, z_2)y_n = R_{g}(\delta z_1, \delta z_2)^{\NGprop}y_n \quad \text{for} \quad \delta = \frac{\NFprop}{\NGprop} \\
		\text{$f$:} & y_{n+1} = R_{f}(z_1, z_2)y_n & \hspace{3em} 				& \text{$\mathcal{F}$:} & y_{n+1} = R_\mathcal{F}(z_1, z_2)y_n = R_{f}(z_1, z_2)^{\NFprop}y_n
	\end{array}
	\label{eq:coarse-fine-stability-functions}
\end{align}
where $h$ is the stepsize and $z_1 = h \lambda_1$ $z_2 = h \lambda_2$. Note that for $g$ and $f$, $y_n$ corresponds to the $n$th fine step, while for $\mathcal{F}$ and $\mathcal{G}$, $y_n$ corresponds to the $n$th coarse step (see \cref{fig:parareal-coarse-fine-grid} for an illustration of coarse and fine steps). The Parareal iteration \cref{eq:parareal_iteration} then reduces to the matrix iteration
\begin{align}
	\mathbf{M}_{\mathcal{G}} \mathbf{y}^{k+1} &= (\mathbf{M}_{\mathcal{G}} - \mathbf{M}_\mathcal{F})\mathbf{y}^{k} + \mathbf{b}
	\label{eq:parareal-matrix-iteration}
\end{align}
where the vector $\mathbf{y}^k = [y_j^k] \in \mathbb{R}^{\Nprocs + 1}$ stores the solution at each coarse step, and
 the matrices $\mathbf{M}_\mathcal{F}, \mathbf{M}_\mathcal{F} \in \mathbb{R}^{\Nprocs+1, \Nprocs+1}$  and vector $\mathbf{b} \in \mathbb{R}^{\Nprocs+1}$ are 
	\begin{align}
		\mathbf{M}_{\mathcal{G}} = 	\left[
				\begin{array}{cccc}
				I \\
				-R_\mathcal{G} & I \\
				   & \ddots & \ddots \\
				   & & -R_\mathcal{G} & I
				\end{array}
			\right], \quad
		\mathbf{M}_{\mathcal{F}} = 
			\left[
				\begin{array}{cccc}
				I \\
				-R_\mathcal{F} & I \\
				   & \ddots & \ddots \\
				   & & -R_\mathcal{F} & I
				\end{array}
			\right], \quad
		\mathbf{b} = 
		\left[
			\begin{array}{c}
				y_0 \\
				0 \\
				\vdots \\
				0
			\end{array}
		\right].
		\label{eq:parareal-matrices}
	\end{align}
Note that the values $R_{\mathcal{G}}$ and $R_{\mathcal{F}}$ are the stability functions from \cref{eq:coarse-fine-stability-functions} that depend on $z_1$ and $z_2$. 

Next, solving the recurrence relation \cref{eq:parareal-matrix-iteration} yields
\begin{align}
	\mathbf{y}^{k+1} = \sum_{j=0}^k \mathbf{E}^j \mathbf{M}_{\mathcal{G}}^{-1} \mathbf{b}, \quad \text{for} \quad \mathbf{E} = \mathbf{I} - \mathbf{M}_G^{-1}\mathbf{M}_F^{-1},	
	\label{eq:parareal-iteration}
\end{align}
and we can now interpret the Parareal algorithm  as a fixed point iteration that converges to the fine solution 
\begin{align}
	\mathbf{y}_{\mathcal{F}} = y_0 \left[1, R_\mathcal{F}, R_\mathcal{F}^2, \ldots, R_\mathcal{F}^\Nprocs\right]^T  \in \mathbb{R}^{\Nprocs + 1}.
\end{align}
Lastly, if we define the error at the $k$th iteration as $\mathbf{e}^k = \mathbf{y}^k - \mathbf{y}_\mathcal{F}$ (i.e. the difference between the Parareal solution and the serial fine integrator solution), then the error at the $k$th iteration, evolves according to the matrix iteration
\begin{align}
	\mathbf{e}^k = \mathbf{E} \mathbf{e}^{k-1}.
	\label{eq:iteration_matrix_parareal_error}
\end{align}
To obtain \cref{eq:iteration_matrix_parareal_error}, we substitute $\mathbf{y}^k = \mathbf{e}^k + \mathbf{y}_{\mathcal{F}}$ into \cref{eq:parareal-matrix-iteration}, left multiply by $M_G^{-1}$, then simplify using $\mathbf{M}_{\mathcal{F}}\mathbf{y}_\mathcal{F} = \mathbf{b}$. In the following subsections we will use \cref{eq:parareal-iteration} to study stability and \cref{eq:iteration_matrix_parareal_error} to study convergence.

\subsection{Linear stability analysis}
\label{subsec:parareal-linear-stability}

Linear stability analysis \cite[IV.2]{wanner1996solving} is a well-known technique that is used to determine the types of equations for which a time integrator is stable (e.g. diffusive or advective). The analysis proceeds by considering the Dahlquist equation and determining the subset of parameters that lead to a stable iteration. All one-step exponential integrators applied to the partitioned Dahlquist equation \cref{eq:partitioned-dahlquist} reduce to the iteration \cref{eq:one-step-dahlquist-iteration}. The function $R(z_1,z_2)$ is the {\em stability function} of the method and its magnitude must be smaller than or equal to one to guarantee stability. The {\em stability region} of a method contains all the $(z_1,z_2)$ pairs for which this holds true and is formally defined as 
\begin{align}
	S = \left\{ (z_1, z_2) \in \mathbb{C}^2 : |R(z_1, z_2)| \le 1 \right\}.	
	\label{eq:stability-region}
\end{align}

For a fixed set of parameters, Parareal is a one-step method that advances the solution by $\Nprocs$ coarse timesteps, or equivalently, $\Nstepstot$ fine timesteps. If $y_j$ denotes the $j$th fine timestep, then Parareal applied to (\ref{eq:partitioned-dahlquist}) reduces to the iteration
\begin{align}
	y_{(\Nstepstot(n+1))} = R(\Nstepstot z_1,\Nstepstot z_2) y_{(\Nstepstot n)} \quad \text{where} \quad z_1 = h\lambda_1, ~ z_2 = h\lambda_2,
	\label{eq:one-step-dahlquist-iteration-parareal}
\end{align}
$h$ is the stepsize of the fine integrator $f$, and the stability function $R$ is
\begin{align}
	R(z_1, z_2) = \mathbf{c}_2 \left( \sum_{j=0}^k \mathbf{E}^j \right) \mathbf{M}_G^{-1} \mathbf{c}_1, && 
	\begin{aligned}
		\mathbf{c_1} &= \left[ 1, 0, \ldots, 0 \right]^T \in \mathbb{R}^{\Nprocs+1}, \\
		\mathbf{c_2} &= \left[ 0, \ldots, 0, 1 \right] \in \mathbb{R}^{\Nprocs+1}.
	\end{aligned}
\end{align}
This stability function follows directly from \cref{eq:parareal-iteration}; $\mathbf{c}_1$ is equivalent to $\mathbf{b}$ with $y_0 = 1$ and $\mathbf{c}_2$ extracts the solution at the final coarse step.

We now apply linear stability analysis to study Parareal methods with classical and repartitioned ERK integrators. Since we are interested in non-diffusive problems with $\lambda_1,\lambda_2 \in i\mathbb{R}$, we only consider the two-dimensional stability region
\begin{align}
		\widehat{S} = \left\{ (r_1, r_2) \in \mathbb{R}^2 : |R(ir_1, ir_2)| \le 1 \right\}.	
	\label{eq:stability-region-non-diffusive}		
\end{align}
Our aim is to determine if the Parareal method from \cref{tab:nls-parareal-configuration} is stable for the ($r_1$, $r_2$) region \cref{eq:nls-r1-r2-ranges} that encloses the eigenvalues of the discretized nonlinear Schr\"{o}dinger equation. 

\begin{figure}

	\centering
	
	Parareal Stability Regions and Instability Factors
	\rule{\textwidth}{0.5px} \\[0.5em]
	
	\begin{tabular}{cccccc}
	& $k=0$ & $k = 2$ & $k=4$ & $k=6$ \\ 
		\rotatebox{90}{\hspace{2em} Classical ERK} & 
		\includegraphics[height=.2\textwidth,align=b]{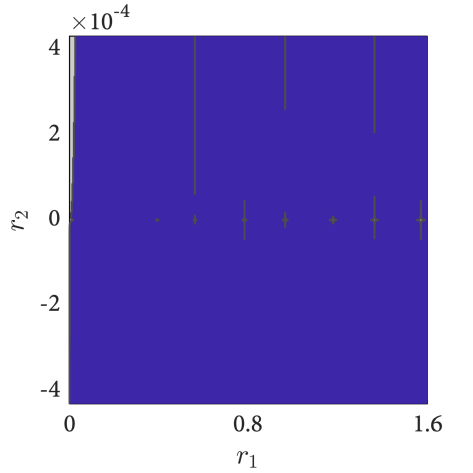} &
		\includegraphics[height=.2\textwidth,align=b]{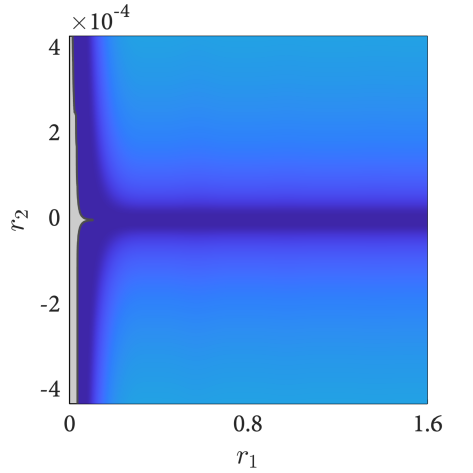} &
		\includegraphics[height=.2\textwidth,align=b]{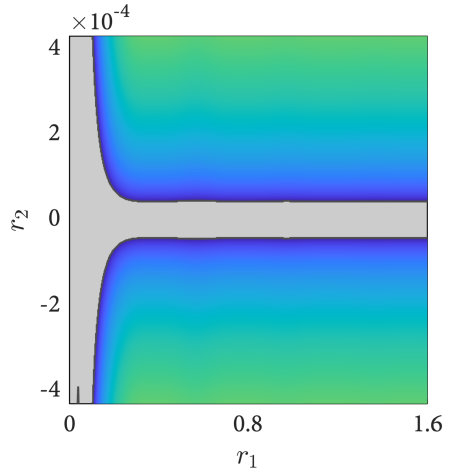} &
		\includegraphics[height=.2\textwidth,align=b]{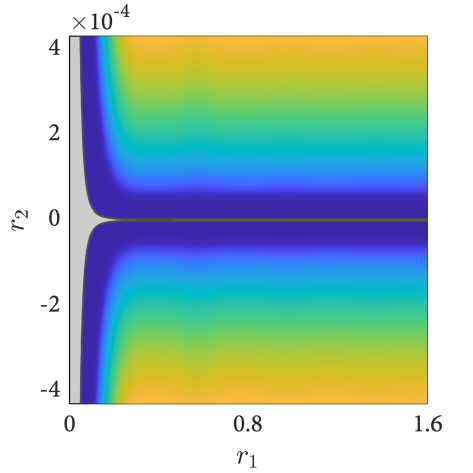} &
		\includegraphics[height=.2\textwidth,trim={6cm 0 0 0},clip,align=b]{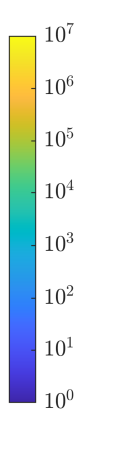} \\
		
		\rotatebox{90}{\hspace{1em} Repartitioned ERK} &
		\includegraphics[height=.2\textwidth,align=b]{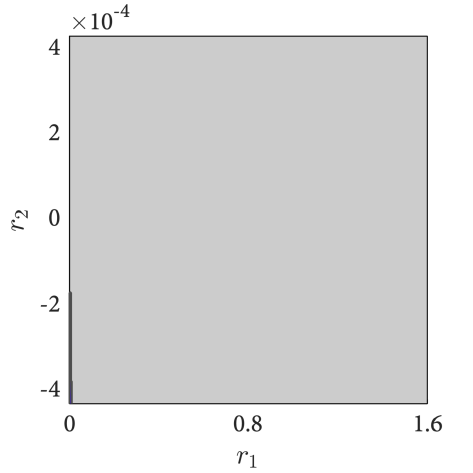} &
		\includegraphics[height=.2\textwidth,align=b]{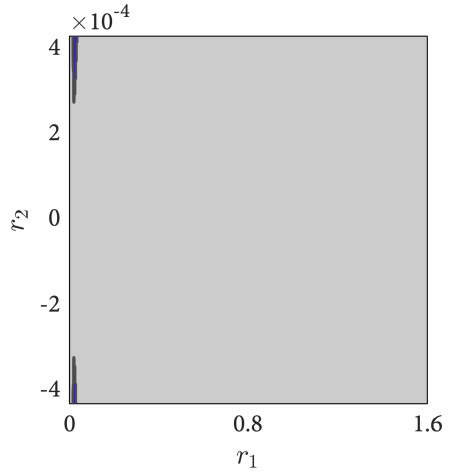} &
		\includegraphics[height=.2\textwidth,align=b]{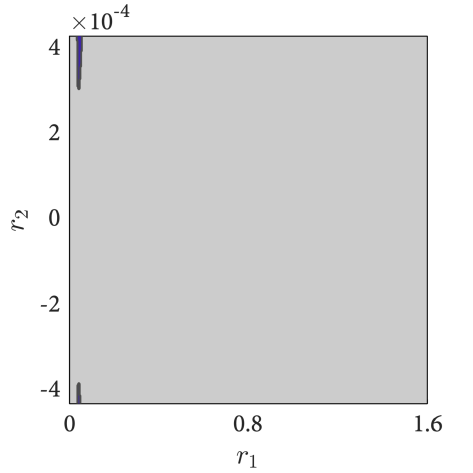} &
		\includegraphics[height=.2\textwidth,align=b]{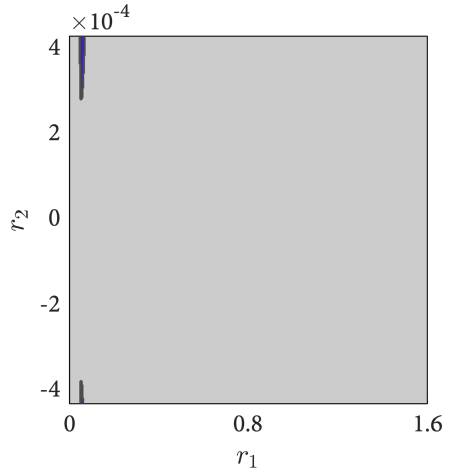} &
		\includegraphics[height=.2\textwidth,trim={6cm 0 0 0},clip,align=b]{figures/nls/nls-stability/stability-surface-reg-leg-k-0-rho-0-1.pdf} \\
	\end{tabular}

	\caption{Stability regions for the Parareal configuration from \cref{tab:nls-parareal-configuration} with classical exponential integrators (top row) and repartitioned exponential integrators (bottom row). Each column corresponds to a different Parareal iteration $k$. The gray region is the stability region \cref{eq:stability-region-non-diffusive}, and color shows the amplification factor $|R(z_1=ir_1,z_2=ir_2)|$ outside the stability region. The $r_1$ and $r_2$ axis limits on the graphs correspond exactly to \cref{eq:nls-r1-r2-ranges}.
	}
	
	\label{fig:parareal-stability-nls}

\end{figure}

In \cref{fig:parareal-stability-nls} we compare the stability regions of the Parareal configuration from \cref{tab:nls-parareal-configuration} equipped with either classical or repartitioned ERK methods. We immediately see that the stability regions associated with classical ERK integrators only encompass a small subset of the rectangular parameter region \cref{eq:nls-r1-r2-ranges} and that the rate of instability worsens significantly as $k$ increases. In contrast, repartitioning greatly expands the stability region of the Parareal method and the remaining instabilities are sufficiently small that they will not affect the quality of the final solution. \Cref{fig:parareal-stability-extra1,fig:parareal-stability-extra2} from \cref{app:additional-stability-plots} contain additional stability plots that reveal a wider range of $r_2$ values. Although repartitioning greatly improves stability, exponential Parareal is only stable when $|r_2| \ll |r_1|$. In other words, the linear term in \cref{eq:semilinear-ode} must contain the majority of the stiffness.
 
 Overall, linear stability analysis is consistent with the convergence diagrams from \cref{fig:nls-experiment-smooth,fig:nls-experiment-oscillatory} which show that Parareal with classical ERK methods grows increasingly unstable as the iteration count $k$ increases. From the linear stability diagrams we also see that the Parareal iteration greatly magnifies the instabilities that are present in the serial ERK4 integrator. For comparison, the serial ERK4 method is stable across the entire range of $(r_1,r_2)$ values shown in \cref{fig:parareal-stability-nls} and its instability rates near the line $r_2 = 0$ are smaller than 1.01; see fig. 3 in \cite{buvoli2022stability}. This implies that repartitioning is important for Parareal even on short timescales where serial exponential integrators do not require it.

\subsection{Linear convergence analysis}
\label{subsec:parareal-linear-convergence}

We now study the convergence rate of the Parareal iteration. We again consider the partitioned Dahlquist equation \cref{eq:partitioned-dahlquist} and determine the subset of parameters that lead to guaranteed rapid convergence. In \cref{subsec:parareal-linear-formulas} we showed that the difference between the Parareal solution and a serial fine integrator solution evolves according to the iteration \cref{eq:iteration_matrix_parareal_error}. Since Parareal fully converges after exactly $\Nprocs{}$ iterations, the matrix $\mathbf{E}$ is nilpotent and the convergence rate cannot be derived from its spectrum. Nevertheless, as originally proposed in \cite{Ruprecht2018}, monotonic convergence is guaranteed if $\|\mathbf{E}\| < 1$ since
	\begin{align}
			\|\mathbf{e}^{k+1}\| \le \| \mathbf{E}\| \|\mathbf{e}^k\| <  \|\mathbf{e}^k\|.
			\label{eq:bounding-error-iteration}
	\end{align} 
	For convergence to occur within a small number of Parareal iterations, we require $\|\mathbf{E}\| \ll 1$; for example if $\|\mathbf{E}\| < 1/10$ it will take 10 iterations to reduce the error by 10 digits.   In \cite{buvoli2021imexparareal} we showed that the $\infty$-norm of $\mathbf{E}$ is
	\begin{align}
		\| \mathbf{E} \|_\infty = \frac{1 - |R_{\mathcal{G}}|^{\Nprocs}}{1 - |R_{\mathcal{G}}|} |R_{\mathcal{G}} - R_{\mathcal{F}}|.
		\label{eq:parareal_iteration_inf_norm}
	\end{align}
	where the values $R_{\mathcal{G}}$ and $R_{\mathcal{F}}$ are the stability functions from \cref{eq:coarse-fine-stability-functions}. The $\infty$-norm is  convenient to use since it is both interpretable and easy to compute.
	
	Using \cref{eq:bounding-error-iteration,eq:parareal_iteration_inf_norm} we define the {\em convergence region} $\mathcal{C}_\infty$ to be the set of all ($z_1$, $z_2$) pairs for which the $\infty$-norm of $\mathbf{E}$ is smaller than one
	\begin{align}
	 	\mathcal{C}_\infty = \left\{ (z_1, z_2) \in \mathbb{C}: \| \mathbf{E}(z_1,z_2) \|_\infty < 1\right\}.	
	 	\label{eq:convergence-region}
	\end{align}
	Since we are only considering non-diffusive equations with $\lambda_1,\lambda_2 \in i\mathbb{R}$, we will study the two-dimensional convergence region
	\begin{align}
		\widehat{\mathcal{C}}_\infty = \left\{ (r_1, r_2) \in \mathbb{R}: \| \mathbf{E}(ir_1,ir_2) \|_\infty < 1\right\}.
		\label{eq:convergence-region-non-diffusive}		
	\end{align}		
	Note that the matrix $\mathbf{E}$ does not depend on the iteration $k$ so a single convergence region pertains to a Parareal configuration with an arbitrary $\NIters$. 
	
\begin{figure}

	\begin{center}
		\renewcommand*{\arraystretch}{1.5}
		\begin{tabular}{c|c|cr}
			{\footnotesize $\mathbf{\NGprop = 1}$} &
			{\footnotesize $\mathbf{\NGprop = 2}$} &
			{\footnotesize $\mathbf{\NGprop = 3}$} \\[1em]
			\includegraphics[height=.26\textwidth]{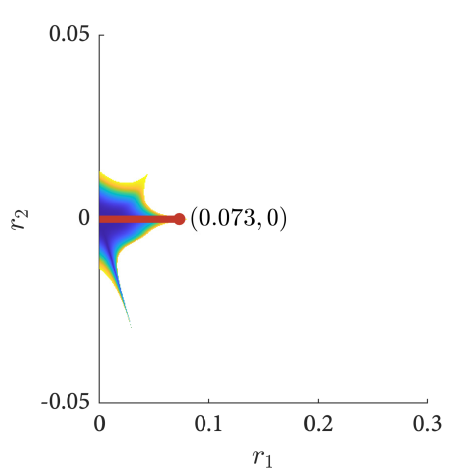} &
			\includegraphics[height=.26\textwidth]{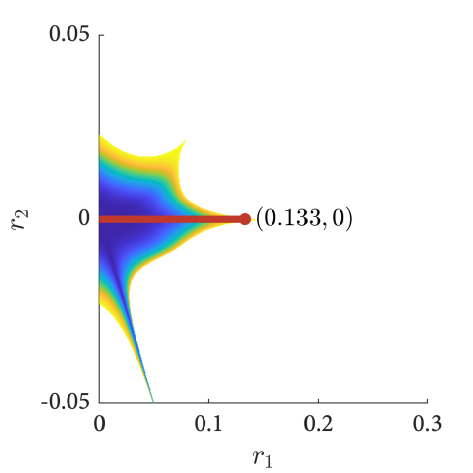} &
			\includegraphics[height=.26\textwidth]{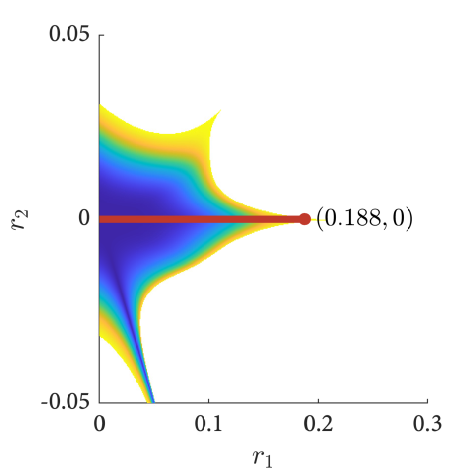} & 
			\includegraphics[height=.26\textwidth,trim={6cm 0 0 0}, clip]{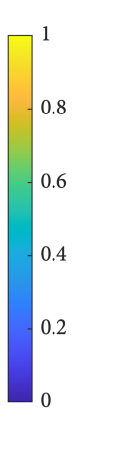} \\	
			\hspace{1em} {\footnotesize Magnified $r_2$ axis} &
			\hspace{1em} {\footnotesize Magnified $r_2$ axis} &
			\hspace{1em} {\footnotesize Magnified $r_2$ axis} \\	[0.5em]			
			\includegraphics[height=.26\textwidth]{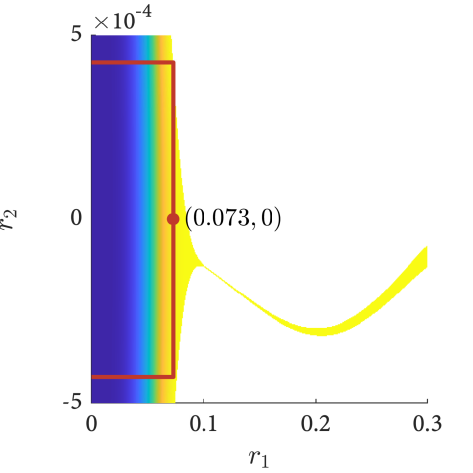} &
			\includegraphics[height=.26\textwidth]{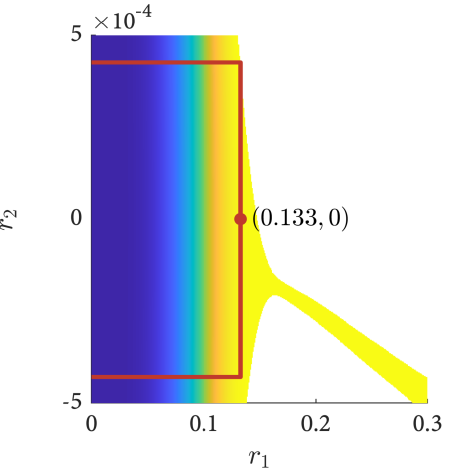} &
			\includegraphics[height=.26\textwidth]{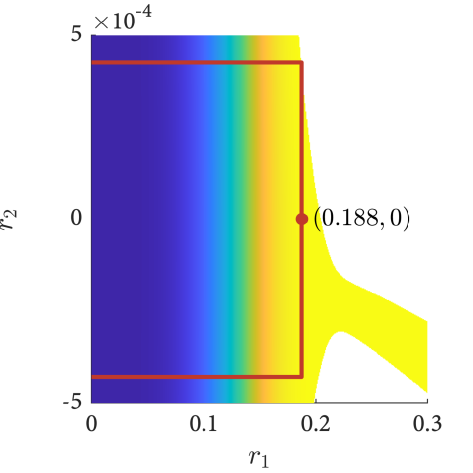}
		\end{tabular}	
	\end{center}

	\caption{Convergence regions \cref{eq:convergence-region-non-diffusive} for the three Parareal configurations from \cref{fig:nls-ng-experiment}. The columns represent different choices for $\NGprop$ and the bottom row shows a magnified $r_2$ axis compared with the top row. Color represents the $\infty$-norm of the matrix $\mathbf{E}$, from \cref{eq:parareal-iteration}, and the red rectangles are the largest rectangular subset of the $r_1$ and $r_2$ range from \cref{eq:nls-r1-r2-ranges} that can be contained inside the convergence region. The $r_1$-coordinate of the labeled point corresponds to the width of the rectangular subset rounded to three digits; it is defined as $r_1^{\max} = \max_\rho$ subject to $\{ r_1 = \rho, z_2 \in [-2h,2h] \} \subseteq \widehat{C}_\infty$ for $h=14/2^{16}$.
	}
	
	\label{fig:nls-parareal-convergence-regions}	
		
\end{figure}

	We now apply linear convergence analysis to understand why high-frequency oscillations cause problems for Parareal and why increasing the number of coarse steps $\NGprop$ improves convergence. We again consider the three Parareal configurations from \cref{fig:nls-ng-experiment} with $\NGprop\in\{1,2,3\}$ and all other parameters from \cref{tab:nls-parareal-configuration}. We are primarily interested to see if the convergence regions of these three Parareal configurations enclose the rectangular region \cref{eq:nls-r1-r2-ranges}.	 In \cref{fig:nls-parareal-convergence-regions}, we present the Parareal convergence regions; the red rectangles in each plot show the largest rectangular subset of \cref{eq:nls-r1-r2-ranges} that can be enclosed inside each convergence region.
		
	The first observation is that the convergence regions near ($r_1=0$, $r_2=0$) grow approximately linearly in size with respect to $\NGprop$. This observation follows directly from \cref{eq:parareal_iteration_inf_norm} since increasing $\NGprop$ makes the coarse propagator more accurate, therefore decreasing the quantity $|R_{\mathcal{G}} - R_{\mathcal{F}}|$ (See \cref{rem:conv-ng} in \cref{app:convergence-scaling}). Overall, linear convergence analysis confirms that increasing $\NGprop$ leads to a Parareal configuration that will resolve a larger number of high-frequency temporal components.
	
	The second observation is that the convergence regions are small and fail to fully enclose \cref{eq:nls-r1-r2-ranges}. More precisely, while all of the $r_2$ range is inside the convergence region, less than twenty percent of the $r_1$ range is included, even when $\NGprop = 3$. However, recall that in \cref{fig:nls-ng-experiment} we saw that the Parareal configuration with $\NGprop=3$ was able to accurately converge to the fine solution; we will explore this fact in more detail in \cref{subsec:parareal-nls-validation}.
	
	Our third and final observation involves convergence rates for small, fixed $(r_1,r_2)$ and follows directly from \cref{rem:conv-ng}. Specifically, if we let $q$ be the order of the coarse integrator, then, for fixed $(r_1,r_2)$, $\|\mathbf{E}\|_\infty = \mathcal{O}(1/\NGprop^q)$. Therefore increasing $\NGprop$ will also increase the Parareal convergence rate.
	
\subsection{Validating convergence results for the nonlinear Schr\"{o}dinger equation}
\label{subsec:parareal-nls-validation}

We now validate how closely the predictions of linear convergence analysis, made using the Dahlquist parameters \cref{eq:nls-z1-z2-ranges}, align with the results from \cref{fig:nls-ng-experiment}. The nonlinear Schr\"{o}dinger equation was spatially discretized using a Fourier spectral discretization that represents the solution as the sum of $M$ Fourier modes, such that
\begin{align}
	u(x,t) = \sum_{n=-M/2,}^{M/2-1} \mathbf{a}_n(t) e^{ikx/4}.
	\label{eq:fourier-spatial-discretization}
\end{align}
For even $M$, the Fourier coefficients $\mathbf{a}_n(t)$ evolve according to \cref{eq:nls-system-discretized,eq:semilinear-ode} with $\mathbf{y} = [\mathbf{a}_0, \ldots, \mathbf{a}_{M/2-1}, \mathbf{a}_{-M/2}, \ldots, \mathbf{a}_{-1}]^T$. Therefore, the differential equation that governs the $n$th coefficient is
\begin{align}
	\dot{\mathbf{a}}_n = \lambda_1(n) \mathbf{a}_n + \left[ N(\mathbf{y}) \right]_{1+\alpha(n)} &\quad \text{where} \quad \lambda_1(n) = in^2 / 16,
\end{align}
$[N(\mathbf{y})]_j$ is the $j$th component of the nonlinearity, and $\alpha(n) = n + M \pmod{M}$; $1+\alpha(n)$ is simply the index of $\mathbf{a}_n$ in the vector $\mathbf{y}$. To conduct linear analysis, we replace the coupled nonlinearity with the decoupled linear term $\lambda_2 \mathbf{a}_n$, with $\lambda_2 \in i[-2, 2]$. It then follows that a Parareal iteration with fine timestep $h$ will monotonically converge to the solution of $\mathbf{a}_n(t)$ only if ($h\lambda_1(n)$, $h\lambda_2$) is inside the convergence region $\widehat{C}_\infty$. In other words, linear convergence analysis predicts that Parareal will only monotonically converge to the subset of Fourier coefficients $\mathbf{a}_n(t)$ for which $n$ satisfies the set inequality
	\begin{align}
		\{ r_1 = r_1(n),~  r_2 \in [-2h,2h] \} \subseteq \widehat{\mathcal{C}}_\infty \quad \text{for} \quad r_1(n) =  h n^2 / 16.
		\label{eq:convergence-conditions-z1-z2}	
	\end{align}
	Ignoring stability considerations, all other $\mathbf{a}_n(t)$ will remain at the accuracy achieved by the coarse integrator until $k \to \Nprocs$. We can simplify the condition \cref{eq:convergence-conditions-z1-z2} by introducing 
	\begin{align}
		r_1^{\text{max}} = \max_{\rho} 
		\quad \text{subject to} \quad 
		\{ r_1 = \rho, r_2 \in [-2h,2h] \}  \subseteq \widehat{\mathcal{C}}_\infty,
		\label{eq:r1-max-def}
	\end{align}
	which is the width of the largest rectangle that includes the entire $r_2$ range and is enclosed by the convergence region. The values of $r_1^{\max}$ for the three Parareal configurations considered in \cref{fig:nls-ng-experiment} are the $r_1$-coordinates of the labeled points in \cref{fig:nls-parareal-convergence-regions}. Using \cref{eq:r1-max-def} it follows immediately that the condition \cref{eq:convergence-conditions-z1-z2} is equivalent to the inequality
\begin{align}
	r_1(n) < r_1^{\text{max}}.
\end{align}

\begin{figure}

	\begin{tabular}{ll}
		\renewcommand*{\arraystretch}{1.9}
		\begin{tabular}{l|ll}
			$\NGprop$ & $r_1^{\text{max}}$ & Convergent $\mathbf{a}_n(t)$ \\	
			\hline
			1 & .073 & $n \in [-73,73]$ \\
			2 & .133 & $n \in [-99,99]$ \\
			3 & .188 & $n \in [-118,118]$
		\end{tabular}	&
		\includegraphics[width=0.5\textwidth,align=c]{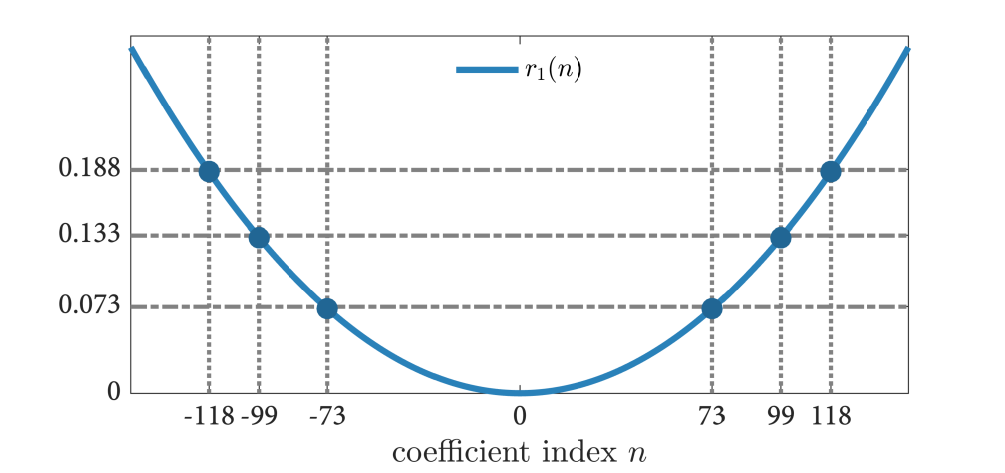}	
	\end{tabular}

	\caption{Monotonically convergent Fourier coefficients, as predicted by linear convergence analysis, for the three Parareal configurations from \cref{fig:nls-ng-experiment} with $\NGprop \in \{1,2,3\}$. The set of convergent mode indices is defined as $\{ n \in \mathbb{Z} : r_1(n) \le r_1^{\text{max}} \}$ for $r_1(n) = hn^2/16$ and $h=14/2^{16}$. The table contains $r_1^{max}$ values (the $x$ coordinates of the labeled points in \cref{fig:nls-parareal-convergence-regions}) along with the indices of the corresponding convergent coefficients. The solid blue line in the plot shows $r_1(n)$, the horizontal dash-dotted gray lines correspond to the three $r_1^{\text{max}}$ values, and the pairs of vertical dashed gray lines are the upper and lower bounds for the predicted convergent coefficients.}

	\label{fig:nls_linear_operator}
	
\end{figure}

\begin{figure}

	\begin{center}
	
		\renewcommand*{\arraystretch}{1.5}
		\begin{tabular}{c|c|c}
			{\footnotesize $\mathbf{\NGprop = 1}$} &
			{\footnotesize $\mathbf{\NGprop = 2}$} &
			{\footnotesize $\mathbf{\NGprop = 3}$} \\[1em]
			{\footnotesize Error in Fourier Coefficient $\mathbf{a}_n(14)$} & 
			{\footnotesize Error in Fourier Coefficient $\mathbf{a}_n(14)$} & 
			{\footnotesize Error in Fourier Coefficient $\mathbf{a}_n(14)$} \\			
			\includegraphics[width=.32\textwidth]{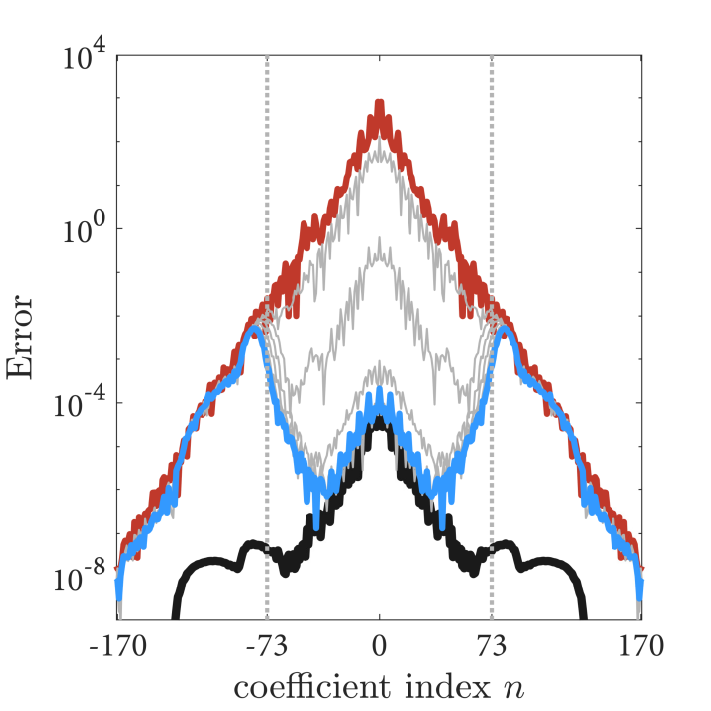} &
			\includegraphics[width=.32\textwidth]{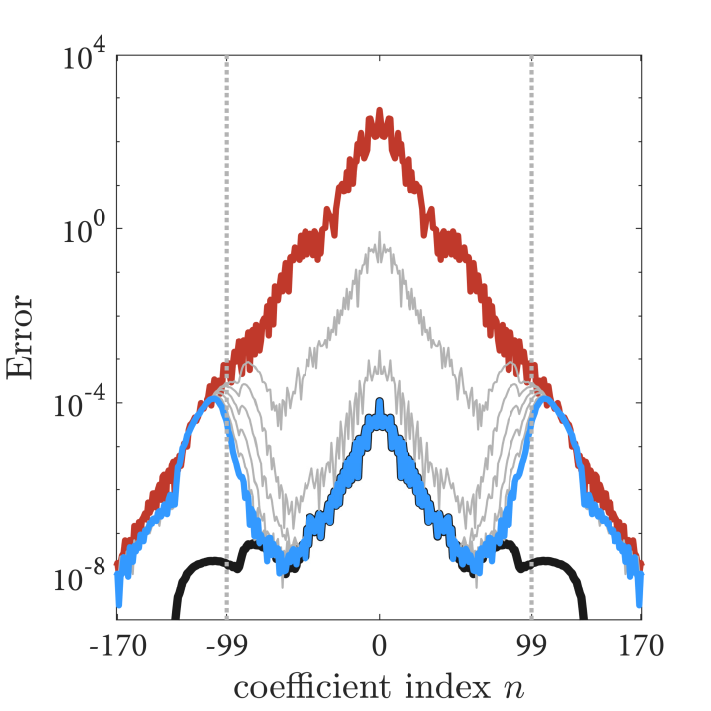} &
			\includegraphics[width=.32\textwidth]{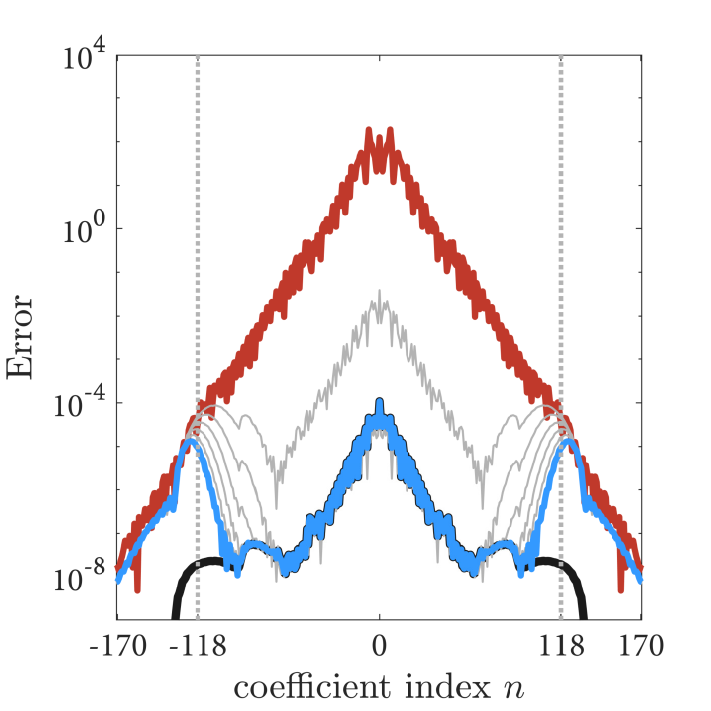} \\	
			\hspace{1em} {\footnotesize Error Norm vs Iteration} &
			\hspace{1em} {\footnotesize Error Norm vs Iteration} &
			\hspace{1em} {\footnotesize Error Norm vs Iteration} \\				
			\includegraphics[width=.32\textwidth]{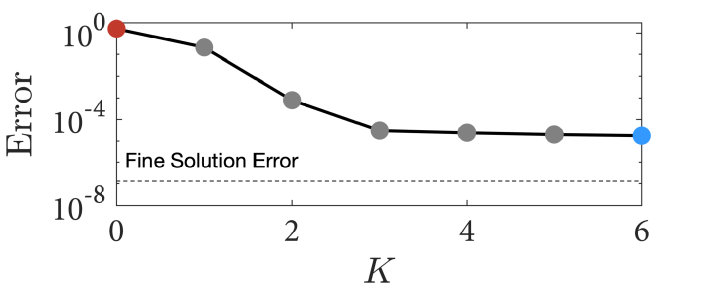} &
			\includegraphics[width=.32\textwidth]{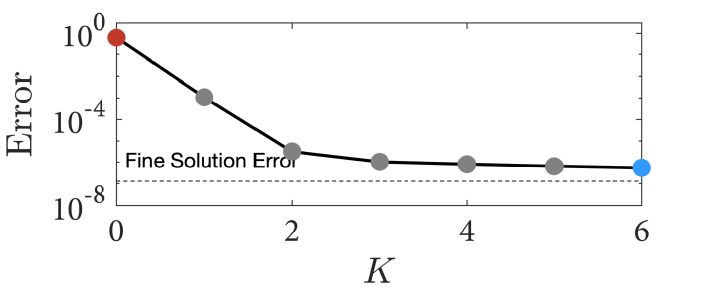} &
			\includegraphics[width=.32\textwidth]{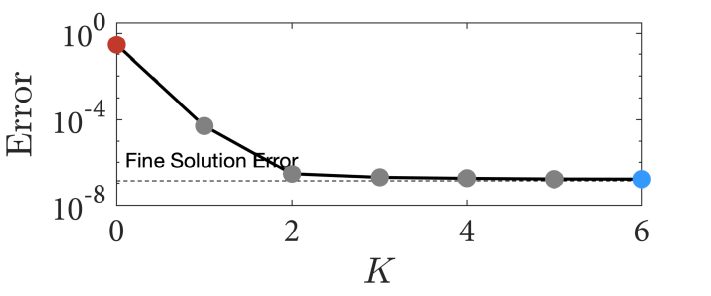}
		\end{tabular}
		
		\vspace{1em}
		
		\begin{tabular}{p{0.96\textwidth}} \hline
			\begin{footnotesize}
				\hfill
				\textcolor{nls_fourier_red}{\hdashrule[0.2ex]{3em}{2.5pt}{}} Coarse ($\NIters = 0$)
				\hfill
				\textcolor{nls_fourier_gray}{\hdashrule[0.2ex]{3em}{1pt}{}} Parareal ($1 \le \NIters \le 5$)
				\hfill
				\textcolor{nls_fourier_blue}{\hdashrule[0.2ex]{3em}{2.5pt}{}} Parareal ($\NIters = 6$)
				\hfill
				\textcolor{nls_fourier_black}{\hdashrule[0.2ex]{3em}{2.5pt}{}} Fine Propagator ($K = \Nprocs$)
				\hfill
			\end{footnotesize}
			\hfill \\[.3em] \hline
		\end{tabular}	
		
	\end{center}

	\caption{Plots describing the same numerical experiment as the one from \cref{fig:nls-ng-experiment}. Each column corresponds to a Parareal configuration with a different value of $\NGprop$ and the colors represent different Parareal iteration numbers. {\bf Top Row:} Error in the Fourier coefficients $\mathbf{a}_n(t)$ of the solution \cref{eq:fourier-spatial-discretization} at $t=14$ for Parareal with $K \in \{0, \ldots, 6\}$. The black line corresponds to the fine propagator $\Fprop$ run in serial. The two vertical dotted lines in each plot are the upper and lower bounds for the convergent spatial modes as predicted by linear analysis and are identical to those shown in \cref{fig:nls_linear_operator}. {\bf Bottom Row:} Error norm between the reference solution and the Parareal method at $t=14$. These plots are identical to the one shown in \cref{fig:nls-ng-experiment}. }
	
	\label{fig:nls-fourier-convergence}
	
\end{figure}

In \cref{fig:nls_linear_operator} we present a table of the $r_1^{\text{max}}$ values for the three parareal configurations from \cref{fig:nls-ng-experiment}, along with the resulting estimates of the monotonically convergent Fourier modes. Then in \cref{fig:nls-fourier-convergence} we validate these estimates by comparing the Fourier coefficients $\mathbf{a}_n(t)$ obtained using the Parareal iteration to those obtained using the serial ERK4 integrator. Overall we see that linear convergence analysis very accurately predicts the convergent Fourier coefficients. Moreover, we see that the accuracy of all Fourier coefficients with an $r_1(n)$ that was outside of the convergence region did not improve substantially beyond what was achieved using the coarse integrator.

\subsection{Convergence regions for additional Parareal configurations}
\label{subsec:convergence-additional-configuration}

\begin{figure}

	\centering
	
	\begin{tabular}{ccccc}
		\hspace{2em} $\Gprop = \text{ERK1}$ & \hspace{2em} $\Gprop = \text{ERK2}$ & \hspace{2em} $\Gprop = \text{ERK3}$ & \hspace{2em} $\Gprop = \text{ERK4}$ \\ 
		\includegraphics[height=.2\textwidth,align=b]{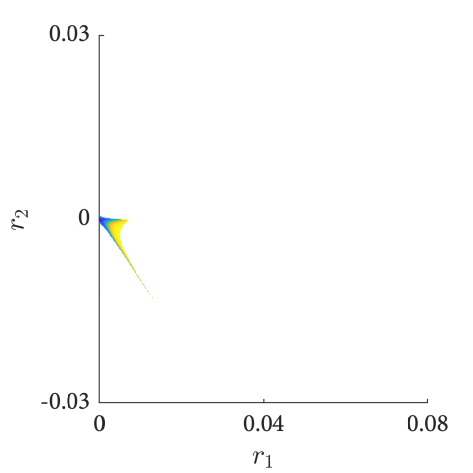} &
		\includegraphics[height=.2\textwidth,align=b]{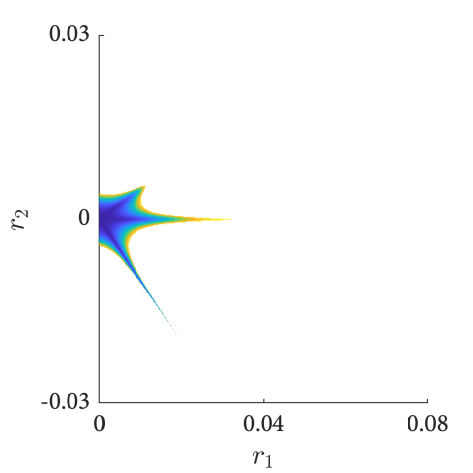} &
		\includegraphics[height=.2\textwidth,align=b]{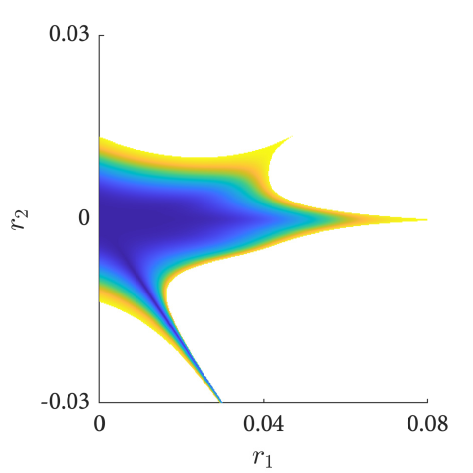} &
		\includegraphics[height=.2\textwidth,align=b]{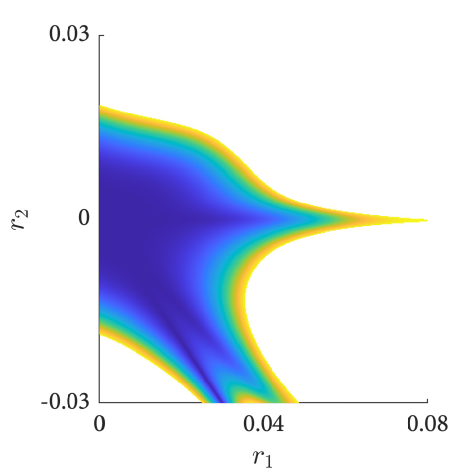} &
		\includegraphics[height=.2\textwidth,trim={6cm 0 0 0},clip,align=b]{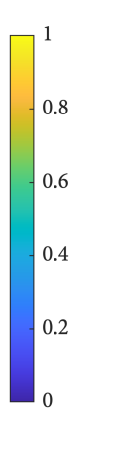} \\	
	\end{tabular}

	\caption{Convergence regions for the Parareal configuration from \cref{tab:nls-parareal-configuration} with different coarse integrators $\Gprop$.}
	
	\label{fig:parareal-convergence-coarse}

\end{figure}

Convergence regions depend on all the Parareal parameters from \cref{tab:parareal-parameters} and on the repartitioning constant $\rho$ from \cref{eq:rho-repartitioning}. Here we investigate the effects of changing the coarse integrator $\Gprop$ and the repartitioning constant $\rho$. \Cref{fig:parareal-convergence-coarse} presents convergence regions for Parareal configurations with ${\Gprop \in \{ \text{ERK1, ERK2, ERK3, ERK4} \}}$ and all other parameters taken from \cref{tab:nls-parareal-configuration}. We see that using a higher-order coarse integrator results in a larger convergence region. This follows directly from \cref{eq:parareal_iteration_inf_norm} since the increased accuracy of a high-order coarse propagator decreases the quantity $|\Gprop - \Fprop|$. Therefore, replacing a low-order coarse propagator with a higher-order one, provides another way to improve convergence for high-frequency temporal modes. Naturally, any convergence gains must be weighed against the decrease in parallel speedup \cref{eq:parareal_speedup} caused by the more expensive high-order coarse propagator.

\begin{figure}

	\centering
	
	\begin{tabular}{ccccc}
		$\rho = 0$ & $\rho = \pi/256$ & $\rho = \pi/64$ & $\rho = \pi/16$ \\ 
		\includegraphics[height=.2\textwidth,align=b]{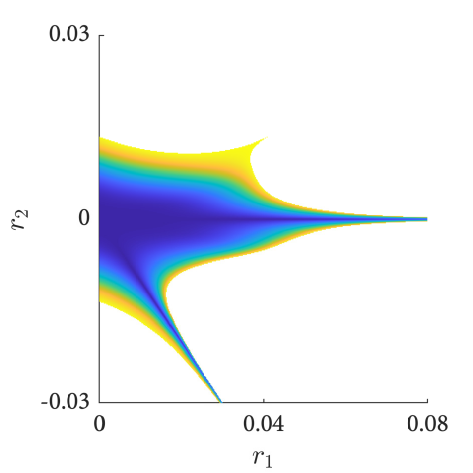} &
		\includegraphics[height=.2\textwidth,align=b]{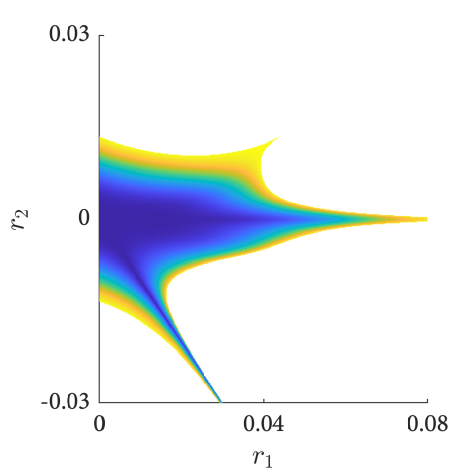} &
		\includegraphics[height=.2\textwidth,align=b]{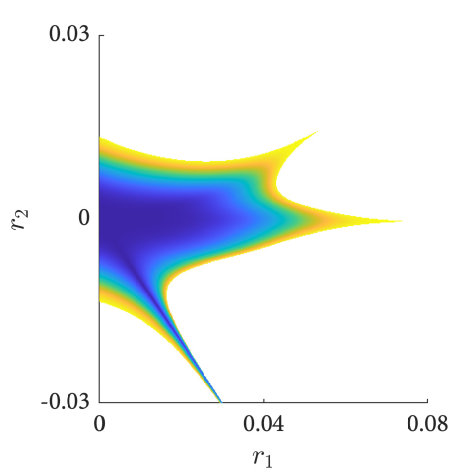} &
		\includegraphics[height=.2\textwidth,align=b]{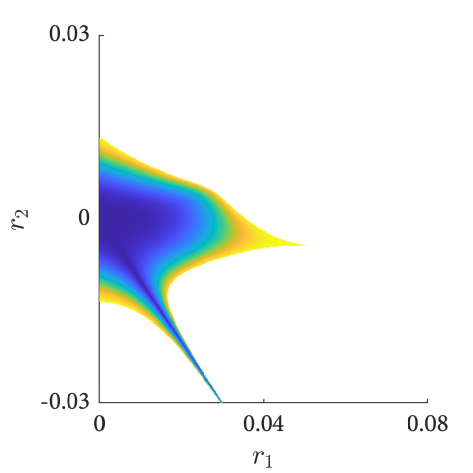} &
		\includegraphics[height=.2\textwidth,trim={6cm 0 0 0},clip,align=b]{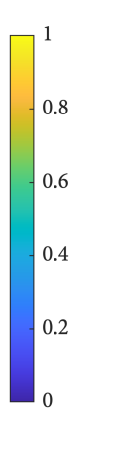} \\
	
	\end{tabular}

	\caption{Convergence regions of the Parareal configuration from \cref{tab:nls-parareal-configuration} with four repartitioning constants $\rho$ from \cref{eq:rho-repartitioning}.
	}
	
	\label{fig:parareal-convergence-rho}

\end{figure}

Next, we briefly discuss how the repartitioning parameter $\rho$ from \cref{eq:rho-repartitioning} affects convergence; recall that $\rho$ is the angle (in radians) that the eigenvalues of the linear operator are rotated into the left-half plane. \Cref{fig:parareal-convergence-rho} contains convergence regions for the Parareal configuration from \cref{tab:nls-parareal-configuration} with repartitioning parameters $\rho \in \{ 0, ~\pi/256,~ \pi/64,~ \pi / 16 \}$. When no repartitioning is applied (i.e. $\rho = 0$) an exponential integrator will exactly solve a Dahlquist equation \cref{eq:partitioned-dahlquist} with $\lambda_2=0$. Therefore, the Parareal convergence region for $\rho = 0$ extends infinitely along the line $r_2 = \text{Im}(h\lambda_2) = 0$. Any amount of repartitioning destroys the exactness of the integrator along this line. In practice this is not important since imposing $\lambda_2=0$ is equivalent to forcing the nonlinearity $N(t,y)$ in \cref{eq:model-ode} to be zero. It should be noted however, that increasing the repartitioning parameter leads to a more subtle contraction of the overall convergence region. Therefore, to maximize convergence for high-frequency information, one should select the smallest repartitioning constant that ensures stability.

\subsection{Implications of convergence analysis for solving non-diffusive PDEs}
\label{subsec:convergence-implications-pde}

As demonstrated by linear analysis, Parareal only converges rapidly to equation components that are not overly oscillatory in time. There are many factors affecting the number of high-frequency temporal modes present in a spatially discretized partial differential equation. High-order dispersive derivatives (e.g. $u_{xxx}$ or $iu_{xx}$) possess continuous spectrums with large imaginary eigenvalues. The extent to which the continuous spectrum causes convergence problems will depend on the choice of spatial discretization. Non-diffusive discretizations, such as Fourier pseudo-spectral methods, will present the greatest challenge since the discretized linear operators will also have purely imaginary eigenvalues. Using fine spatial grids will also increase the total number of high-frequency components. Finally, there is the question of whether high-frequency information is necessary for obtaining accurate solutions. This will depend on the initial condition, the length of the integration window, and the characteristics of the problem. PDEs or initial conditions that cause rapid spectral broadening within the integration window will be the most challenging to solve using Parareal. In summary, limited convergence for high-frequency temporal components will manifest as inaccuracies in high-frequency spatial modes of the solution.

For the generic initial value problem \cref{eq:semilinear-ode}, we can extend the linear analysis developed in this section to estimate the convergence properties of a given Parareal configuration. We start by making the following two assumptions:
\begin{enumerate}
		\item The linear operator $\mathbf{L}$ is diagonalizable and $(\lambda_n, \mathbf{v}_n)$ are the $n$th eigenvalue and eigenvector. This allows us to express the solution as a linear combination of the eigenvectors, $\mathbf{y}(t) = \sum_{n=0}^M a_n(t) \mathbf{v}_j$. 
		\item The linear operator contains the majority of the stiffness such that $\rho(\mathbf{L}) \gg \rho(\frac{\partial N}{\partial y})$.
\end{enumerate}
We then bound the spectrum of the linear and nonlinear operators and define the region
\begin{align}
	\mathcal{R}(h) = \left\{ r_1 \in h [0,~ c_1], r_2 \in h [-c_2,~ c_2] \right\} \quad c_1 = \rho(\mathbf{L}), \quad c_2 = \max_{t} \rho \left(\tfrac{\partial N}{\partial \mathbf{y}}(\mathbf{y}(t))\right).
	\label{eq:convergence-region}
\end{align}
To proceed we must first ensure that a given Parareal configuration is sufficiently stable for all ($r_1$, $r_2$) in $\mathcal{R}(h)$. If this holds true, we then determine the largest rectangular subset of $\mathcal{R}(h)$ that: (i) is enclosed by the Parareal convergence region \cref{eq:convergence-region-non-diffusive} and (ii) contains the entire $r_2$ range. The width of this region is  
\begin{align}
	r_1^{\text{max}}(h) = \max_{\omega} 
	\quad \text{subject to} \quad 
	\{ r_1 = \omega, r_2 \in h[-c_2,c_2] \}  \subseteq \widehat{\mathcal{C}}_\infty.
\end{align}
The value of $r_1^{\text{max}}(h)$ will depend on all Parareal parameters except for $\NIters$. Finally, we estimate that a Parareal iteration with fine stepsize $h$ will monotonically converge to the fine integrator solution of any coefficient $a_n(t)$ where $n$ satisfies 
\begin{align}
	|h\lambda_n| < r_1^{\text{max}}(h).
\end{align}
All the remaining coefficients will retain the accuracy achieved with the coarse integrator and fail to converge for small iteration count $k$. Although these estimates are rooted in linear theory, our results in \cref{subsec:parareal-nls-validation} demonstrate that this approach has the potential to accurately predict convergence for nonlinear problems.
\section{Higher-dimensional numerical experiments}
\label{sec:higher-dimensional-numerical-experiments}

We now consider two-dimensional, non-diffusive equations and demonstrate that exponential Parareal can achieve reduced time-to-solution compared to serial exponential integrators. Specifically, we conduct two additional numerical experiments in which we solve the dispersive Kadomtsev-Petviashvili (KP) equation and the hyperbolic Vlasov-Poisson (VP) equation. Both PDEs are equipped with periodic boundary conditions and discretized in space using a Fourier spectral method. Since analytical solutions are not known, we compute a reference solution using ERK4 with a very small timestep. The error is then defined as
	\begin{align}
		\| \mathbf{y}_\text{ref} - \mathbf{y}_\text{method}	\|_\infty  / \|\mathbf{y}_{\text{ref}}\|_\infty.
	\end{align}
where $\mathbf{y}$ represents the solution in physical space. Below we describe the equations, their initial conditions, and the corresponding numerical parameters. 
	
	The {\bf Kadomtsev-Petviashvili} (KP) equation is
	\begin{align}
		\left(u_t + 6uu_x + u_{xxx}\right)_x + 3 \sigma^2 u_{yy} = 0	
		\label{eq:kp-non-evolution}
	\end{align}
	where $\sigma^2 = -1$ leads to KPI that models thin films with large surface tension, while $\sigma^2 = 1$ leads to KPII that models water waves with small surface tension \cite{biondini2008kadomtsev}. Both equations admit the soliton solution $u(x,y,t) = 2p^2 \text{sech}(p(x-4p^2t))$ where $p$ is a free variable. The stability of a perturbed soliton depends on the sign of $\sigma^2$, with the KPI solution being unstable and the KPII solution being stable \cite{frycz1990bending,infeld1994decay,infeld1995numerical}. 
	
	For any well-localized solution in $x$, the KP equation can be expressed in evolution form as
	\begin{align}
		u_t + 6 uu_x + u_{xxx} + 3 \sigma^2 \partial^{-1} u_{yy} = 0 \quad \text{where} \quad  \partial^{-1} f = \frac{1}{2} \left[\int_{-\infty}^x f(s) ds - \int^{\infty}_x f(s) ds\right].
		\label{eq:kp-evolution}
	\end{align}
	To ensure smoothness in time, the initial condition must satisfy the following equality at $t=t_0$
	\begin{align}
		\int_{-\infty}^\infty u_{yy}(x,y,t)dx = 0.
		\label{eq:kp-constraint}
	\end{align}
	If the initial condition does not satisfy this constraint, then an infinitesimally short but infinitely large change occurs in the solution so that \cref{eq:kp-constraint} is satisfied for all $t>t_0$ \cite{ablowitz1991kadomtsev}; this results in a discontinuous solution (in time) at $t=t_0$.
	
	For our numerical experiment we consider the KPI equation equipped with periodic boundary conditions on the domain $x \in [-8\pi, 8\pi]$, $y \in [0, 8\pi]$. We spatially discretize using $972$ grid points in $x$ and $750$ grid points in $y$, and dealias using the standard 3/2 rule. We integrate the equation in Fourier space where the operator $\partial^{-1}$ is equivalent to the Fourier multiplier $-i/k_x$; Note that when $k_x=0$ this mode is singular. However, for any initial condition that satisfies \cref{eq:kp-constraint} we can simply set this multiplier to zero; for more general initial conditions, numerical regularization must be added \cite{klein2007numerical}.
	
\begin{figure}[h]
		
		\centering
		\includegraphics[width=0.8\textwidth,trim={0 0 0 1cm}, clip]{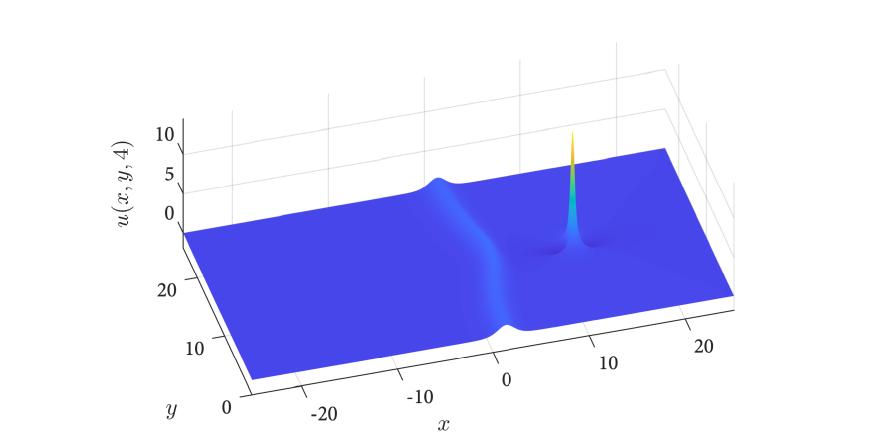}

		\caption{Solution of the KPI equation \cref{eq:kp-evolution} at $t=4$ with initial conditions \cref{eq:kp-initial-condition}.}	
		\label{fig:kp-solution-physical}
		
	\end{figure}

	As in \cite{infeld1994decay} we select our initial condition to be a soliton with a perturbed phase
	\begin{align}
		u(x,y,t=0) = 2\text{sech}^2\left((x + 4 \pi) + \delta \cos(y/4)\right), 	\quad \delta = 1/5,
		\label{eq:kp-initial-condition}
	\end{align}
	and integrate the equation to time $\Tfinal=4$. Our initial condition satisfies \cref{eq:kp-constraint}, therefore, no regularization is needed. As shown in \cref{fig:kp-solution-physical}, the perturbation is unstable and leads to the formation of a two-dimensional soliton. 
	
\begin{figure}[h]	
		\centering
		\includegraphics[width=0.7\textwidth,trim={0 0 0 1cm}, clip]{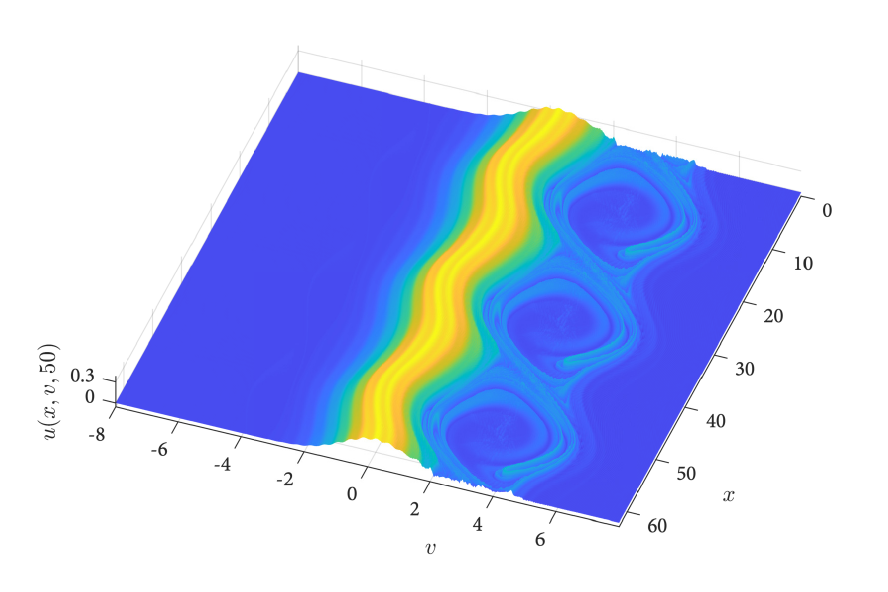}
		
		\caption{Solution of the Vlasov-Poisson equation \cref{eq:valsov-poisson} at $t=50$ for the initial condition \cref{eq:vp-bump-on-tail-initial-condition}.}	
		\label{fig:vp-solution-physical}
	\end{figure}

	The hyperbolic {\bf Vlasov-Poisson} (VP) equation is
	\begin{align}
		f_t + v f_x + E(x,t) f_v = 0, \quad \text{for} \quad
		E_x(x,t) = -1 + \int_{-\infty}^{\infty} f(x,v,t) dv,
		\label{eq:valsov-poisson}
	\end{align}
	and describes the evolution of charged particles in an electric field \cite{glassey1996cauchy}. Our numerical experiment is based on the bump-on-tail experiment from \cite{crouseilles2020exponential}. Specifically, we equip the VP equation with periodic boundary conditions on the domain $x \in [0, 20\pi]$, $v \in [-8, 8]$, and spatially discretize using a 1024 point Fourier discretization in both $x$ and $v$. Our initial condition is
	\begin{align}
		f(x,v,t=0) = \left(\frac{0.9}{\sqrt{2\pi}} e^{-v^2 / 2} + \frac{0.2}{\sqrt{2 \pi}} e^{-2(v - 4.5)^2} \right)  \left( 1 + \frac{4}{100} \cos(0.3 x)\right)
		\label{eq:vp-bump-on-tail-initial-condition}
	\end{align}
	and the solution is integrated to time $\Tfinal=50$. To preserve a diagonal linear operator we solve the equation in physical $v$ space and Fourier $x$ space; see \cref{eq:discrete-vp} in \cref{app:discrete-vp}.  As shown in \cref{fig:vp-solution-physical} the bump-on-tail initial condition excites modes that lead to complex dynamics. 
	
\subsection{Parareal parameter selection and experiment overview}
	
	The Parareal configurations we selected to solve the KP and VP equations are described in \cref{tab:2d-parareal-configurations}. For both equations, we considered multiple configurations that differ only in the number of coarse steps $\NGprop$. We vary this parameter to demonstrate the improved convergence properties associated with larger $\NGprop$ values. For the fine integrator $f$ we always selected ERK4 and set the total number of steps $\Nstepstot$ so that a fully-converged Parareal method produces a highly accurate solution. The remaining parameters were then determined by balancing the achievable parallel speedup with the size of the convergence regions.
	
	The results for the numerical experiments in this section will be summarized in three plots: (i) error versus iteration $\NIters$, (ii) parallel speedup \cref{eq:parareal_speedup} versus iteration $\NIters$, and (iii) error versus run-time (i.e. the computational time required to achieve a specified error). In the error versus run-time plots we also include results for the coarse and fine ERK integrators run in serial. All experiments were performed using 32-core Haswell nodes on the Cray XC40 {\em Cori} at the National Energy Research Scientific Computing Center. For the VP equation we collected timing results by running the Parareal iteration on 64 nodes (2048 compute cores). The full KP experiment requires 256 nodes (8096 compute cores) which exceeded our available computational resources. Therefore, we ran the Parareal iteration in serial on a single node to determine the convergence curves, and then extrapolated the achievable speedup from a smaller experiment where we solved the KP equation on the shortened interval $t\in[0,1]$ using 64 nodes.
		
	Lastly, for each Parareal configuration we compute the convergent spatial modes as predicted by the linear stability analysis from \cref{subsec:convergence-implications-pde}. For simplicity, we assume that any stiffness in the nonlinear term is negligible so that convergence depends exclusively on the eigenvalues of the linear operator (i.e. $c_2 = 0$ in \cref{eq:convergence-region}). 	
		
\begin{table}[h]
	
	\begin{subtable}[t]{\textwidth}
		\caption{{\bf Kadomtsev-Petvaishvili (KP)}}
		\label{subtab:2d-parareal-configurations-kp}	
	
		\vspace{1em}
		\begin{minipage}[t]{0.5\textwidth}
			\centering
			{\em Parareal Parameters} \\[1.35em]
			\renewcommand*{\arraystretch}{1.25}
		\begin{tabular}{ll|ll|ll|ll}
			$\Nprocs$ 		& $8192$ 	& $\NFprop$ 		& $32$	& $f$ 	& ERK3   & $K$  & $1, \ldots, 28$\\ \hline
			$\Nstepstot$		& $2^{18}$ 	& $\NGprop$ & \{1, 2, 3\} 	& $g$	& ERK4  
		\end{tabular}
	
		\end{minipage}
		\begin{minipage}[t]{0.40\textwidth}
			\centering
			{\em $r^{\text{max}}_1$ values} \\[1em]
			\renewcommand*{\arraystretch}{1.25}
			\begin{tabular}{l|lll}
										& $\NGprop = 1$ & $\NGprop = 2$ & $\NGprop = 3$ \\ \hline
				$r^{\text{max}}_1(h)$ 	& 0.0642			& 0.112			& 0.155
			\end{tabular}
		\end{minipage}
	
	\end{subtable}
	
	\vspace{1em}
	\begin{subtable}[t]{\textwidth}
		\caption{{\bf Vlasov-Poisson (VP)}}
		\label{subtab:2d-parareal-configurations-vp}	
		
		\vspace{1em}
		\begin{minipage}[t]{0.5\textwidth}
		\centering
		{\em Parareal Parameters} \\[1.35em]
		\renewcommand*{\arraystretch}{1.25}
		\begin{tabular}{ll|ll|ll|ll}
			$\Nprocs$ 		& $2048$ 	& $\NFprop$ 		& $64$	& $f$ 	& ERK4   & $K$  & $1, \ldots, 12$\\ \hline
			$\Nstepstot$		& $2^{17}$ 	& $\NGprop$ & \{2, 3\} 	& $g$	& ERK3   & 
		\end{tabular}	
	\end{minipage}
	\begin{minipage}[t]{0.40\textwidth}
		\centering
		{\em $r^{\text{max}}_1$ values} \\[1em]
		\renewcommand*{\arraystretch}{1.25}
		\begin{tabular}{l|ll}
								& $\NGprop = 2$ & $\NGprop = 3$ \\ \hline
			$r^{\text{max}}_1(h)$ 	&     0.074 & 0.102
		\end{tabular}
	\end{minipage}
		
	\end{subtable}
	
	\caption{Parareal parameters used to solve the KP equation \cref{eq:kp-evolution} with initial conditions \cref{eq:kp-initial-condition} and the 
	Vlasov-Poisson equation \cref{eq:valsov-poisson} with initial conditions \cref{eq:vp-bump-on-tail-initial-condition}. In the right-most tables we present the associated $r_1^{\text{max}}$ values. All $r_1^{\text{max}}$ values are computed using the stepsize $h = \Tfinal/\Nstepstot$, the repartitioning parameter $\rho = \pi/128$, and by assuming that the stiffness in the nonlinear term is negligible such that $c_2 = 1$ in \cref{eq:convergence-region}.}
	\label{tab:2d-parareal-configurations}
\end{table}


\subsection{Kadomtsev-Petvaishvili -- results and discussion}
	
	The KP equation is challenging to solve because the third-order derivative term leads to a linear operator with large imaginary eigenvalues. As we have seen in \cref{subsec:parareal-nls-validation,subsec:convergence-implications-pde}, these eigenvalues determine the degree of temporal oscillation present in the spatial Fourier coefficients of the solution. Therefore, we must select a Parareal configuration whose convergence region contains at least a large subset of these eigenvalues. Our proposed configurations are described in \cref{subtab:2d-parareal-configurations-kp}. The choice of total steps $\Nstepstot$ ensures that a fully-converged Parareal iteration will produce a solution with an accuracy of $1.2\times 10^{-9}$. Moreover, we selected ERK3 as the coarse integrator because it is stable at sufficiently large stepsizes (see \cref{fig:erk-convergence-kp-vp}) and the convergence region for an ERK3, ERK4 pairing is larger than that of an ERK2, ERK4 pairing (see \cref{fig:parareal-convergence-coarse}).
	
\begin{figure}[h]
		
		\centering	     
		\includegraphics[width=0.48\textwidth,trim={0 7cm 0 0},clip]{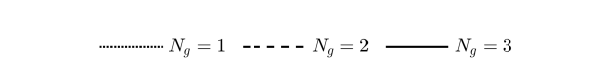}
		
		\begin{subfigure}[b]{0.48\textwidth}
			\centering
			\includegraphics[width=\textwidth]{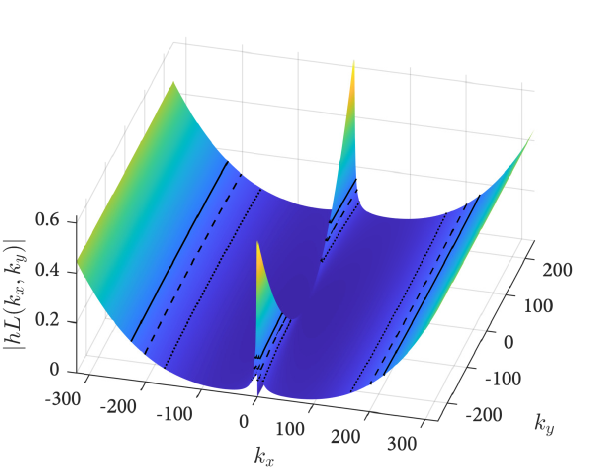}
			\caption{Magnitude of the scaled KP linear operator eigenvalues}
			\label{fig:kp_linop}
		\end{subfigure}
		\hfill
		\begin{subfigure}[b]{0.48\textwidth}
			\centering
			\includegraphics[width=\textwidth]{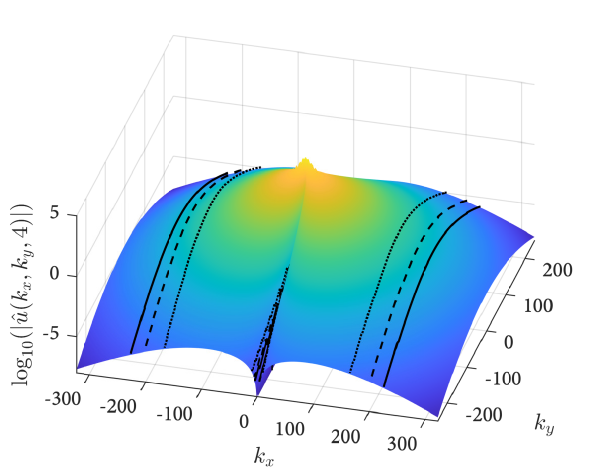}
			\caption{KP Solution at $t=8$ (Fourier transformed in $x$,$y$)}
			\label{fig:kp_fourier}
		\end{subfigure}
		
		\caption{
			Magnitude of the linear operator eigenvalues (\cref{fig:kp_linop}) and final solution in Fourier space (\cref{fig:kp_fourier}) for the KP equation. The axis in both plots are the Fourier wavenumbers in $x$ and $y$. The black contours shows the boundaries of the region $\mathcal{A}$ from \cref{eq:kp-convergence-region} that encloses all convergent spatial modes as predicted by linear convergence analysis. The dotted, dashed, and solid line respectively correspond to the region $\mathcal{A}$ for Parareal configurations with $\NGprop = 1$, $\NGprop = 2$, and $\NGprop = 3$. From linear analysis, we expect that each Parareal iteration will converge monotonically for all spatial modes inside its corresponding region $\mathcal{A}$. 
		}	
		\label{fig:kp-convergence}
		
	\end{figure}

	We can estimate the convergent spatial modes for each choice of $\NGprop$ using the stability analysis developed in \cref{subsec:convergence-implications-pde}. The eigenvalues of the linear operator $\mathbf{L}$ are
	\begin{align}
		\mathbf{L}(k_x,k_y) = 
		\left\{
		\begin{array}{ll}
			-i (\omega_x k_x)^3 & k_x = 0 \\
			-i (\omega_x k_x)^3 + i \frac{(\omega_y k_y)^2}{\omega k_x} & k_x \ne 0	
		\end{array}
		\right.
		\label{eq:kp-convergence-region}
	\end{align}
	where $k_x$ and $k_y$ are integer Fourier wavenumbers and $\omega_x = 2\pi / L_x = 1/8$, $\omega_y = 2\pi / L_y = 1/4$. The convergent spatial modes lie inside the region
	\begin{align}
		\mathcal{A} = \left\{ (k_x, k_y) \in \mathbb{Z}^2 : h|\mathbf{L}(k_x, k_y)| < r_1^{\text{max}}(h) \right\}
		\quad \text{for} \quad h=2^{-16},
		\label{eq:kp-convergence-region}
	\end{align}
	where the values of $r_1^{\text{max}}$ depend on $\NGprop$ and are contained in \cref{subtab:2d-parareal-configurations-kp}.	In \cref{fig:kp-convergence} we overlay the region $\mathcal{A}$ onto the linear operator $\mathbf{L}$ and the Fourier transformed final solution. This allows us to estimate which spatial modes will be accurately computed by each Parareal configuration. Although none of the regions $\mathcal{A}$ enclose the entire ($k_x$, $k_y$) domain, the coefficients for the highest-frequency spatial modes are small and only need to be resolved if we require an extremely accurate solution.
	
\begin{figure}
	
		\hfill
		\begin{subfigure}[t]{0.32\textwidth}
			\centering
			\caption{Error versus Iteration}
			\includegraphics[width=\textwidth,trim={0 0 0.70cm 0}, clip]{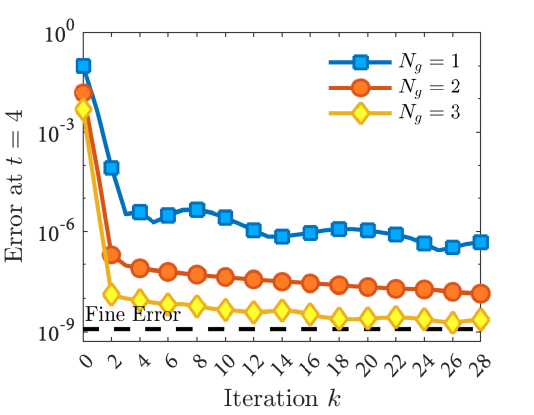}
			\label{fig:kp_error-vs-iteration}
		\end{subfigure}
		\hfill
		\begin{subfigure}[t]{0.32\textwidth}
			\centering
			\caption{Parallel Speedup \cref{eq:parareal_speedup} versus Iteration}
			\includegraphics[width=\textwidth,trim={0 0 0.70cm 0}, clip]{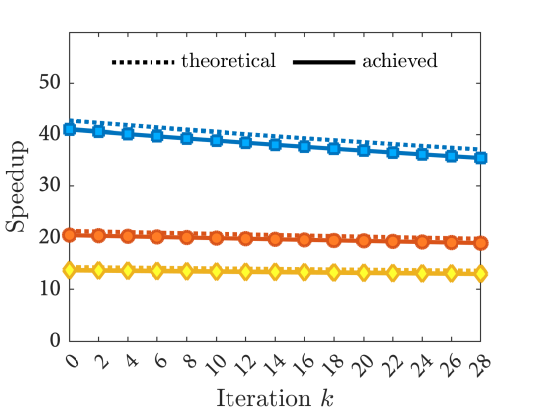}
			\label{fig:kp_parallel-speedup}
		\end{subfigure}
		\hfill
		\begin{subfigure}[t]{0.32\textwidth}
			\centering
			\caption{Error versus Runtime}
			\includegraphics[width=\textwidth,trim={0 0 0.70cm 0}, clip]{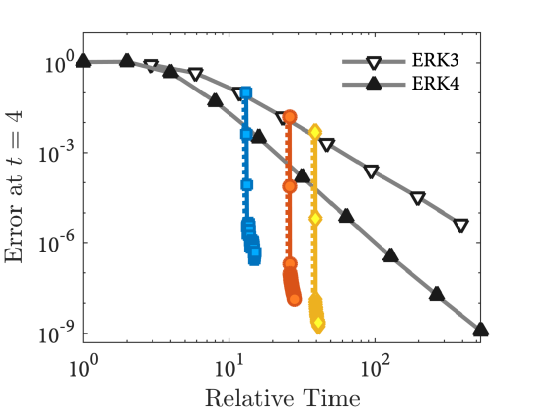}
			\label{fig:kp_error-vs-time}
		\end{subfigure}
		\hfill
		
		\caption{Numerical results for the KP equation using the three Parareal configurations from \cref{tab:2d-parareal-configurations}, that are differentiated with different colored lines. (a) Solution error at the final time $t=4$ as a function of Parareal iteration. (b) Theoretical speedup \cref{eq:parareal_speedup} and achieved speedup from the numerical experiment. (c) Error versus runtime diagram comparing the Parareal configurations to the serial coarse and fine integrators.}		
		\label{fig:kp-ng-experiment}
		
	\end{figure}

	In \cref{fig:kp-ng-experiment} we show convergence, speedup, and error vs runtime results for the exponential Parareal methods applied to the KP equation. We divide our discussion of the results into three parts:
	\begin{enumerate}
		
		\item {\em Convergence}. No exponential Parareal method converged monotonically across the entire $k$ range. Instead they exhibited a rapid reduction in error during the first few iterations before entering a plateau of slow convergence. This behavior is expected, since none of the Parareal convergence regions enclose the full spectrum of the discretized KP linear operator. As predicted by linear analysis, increasing $\NGprop$ improves convergence, and the Parareal configuration with $\NGprop = 3$ is able to obtain a solution that is comparable in accuracy to the serial fine integrator. In contrast, the Parareal configurations with $\NGprop \in \{1,2\}$ do not resolve a sufficient number of high-frequency spatial modes to achieve fine error within 28 iterations; nevertheless, both methods improve the coarse solution by multiple orders of magnitude. These results are analogous to those for the NLS equation shown in \cref{fig:nls-ng-experiment}. 
		
		\item {\em Parallel speedup \cref{eq:parareal_speedup}.} \Cref{fig:kp_parallel-speedup} shows theoretical and achieved parallel speedup in dashed and solid lines, respectively. Unlike our one-dimensional experiments, we now see good agreement between the two types of curves. This is due to the fact that the KP equation is more computationally expensive to integrate over a single timestep than the one-dimensional NLS equation; specifically, the right-hand-side evaluations now require multiple two-dimensional discrete Fourier transforms as opposed to a single one-dimensional transform. In summary, although the penalties incurred due to communication costs are mildly visible (notice that all dashed lines are slightly below the solid lines in \cref{fig:kp_parallel-speedup}), the theoretical speedup estimate \cref{eq:parareal_speedup} provides a realistic measure for real-world speedup of the Parareal iteration. However, we remark that high parallel speedup does not imply convergence or low error; it simply characterizes the runtime of the Parareal iteration. To understand the practical effectiveness of the Parareal iteration, we must investigate error versus runtime. 
		
		\item {\em Error versus runtime.} Perhaps the most important result is the error versus runtime plot, from which we see that all three Parareal configurations are able to compute high-accuracy solutions significantly faster than the serial ERK methods. Moreover, despite their failure to converge to the fine solution within 28 iterations, the Parareal configurations with $\NGprop \in \{1,2\}$ are the fastest methods for obtaining moderately less accurate solutions. We summarize the improvements of the Parareal configurations with $k=28$ over serial ERK4 in the table below:		
			\begin{center}
				\renewcommand*{\arraystretch}{1.25}
				\begin{tabular}{l|ll}
					$\NGprop$ & Error Tol at $k=28$ & Improvement compared to serial ERK4 \\\hline			
					1 & $2.7\times 10^{-7}$  & 9.09x faster  \\ 
					2 & $1.3 \times 10^{-8}$ & 10.15x faster \\ 
					3 & $1.8 \times 10^{-9}$ & 11.71x faster		
				\end{tabular}
			\end{center}
			Lastly we note that \cref{fig:kp_error-vs-time} shows the error versus runtime for the Parareal method using both theoretical and achieved speedup (dashed and solid lines respectively). We see that the losses due to communication only have a very minor effect on performance, since the dashed curves lie just to the left of the solid curves.
	
	\end{enumerate}

\subsection{Vlasov-Poisson -- results and discussion}
	
	The Vlasov-Poisson equation does not contain high-order spatial derivatives, therefore it is possible to select Parareal parameters that simultaneously offer good Parallel speedup and convergence properties. Our proposed Parareal configurations are described in \cref{subtab:2d-parareal-configurations-vp}. The choice of total steps $\Nstepstot$ ensures that a fully-converged Parareal iteration will produce a solution with an accuracy of $6.6\times 10^{-8}$ (see \cref{fig:erk-convergence-kp-vp}).
	
	We again apply linear analysis to estimate the convergent spatial modes for each choice of $\NGprop$. The eigenvalues of the linear operator $\mathbf{L}$ are
	\begin{align}
		\mathbf{L}(k_x,v) = i \omega_x k_x v 
		\quad \text{for} \quad
		k_x \in \mathbb{Z}, ~
		v \in \mathcal{V},
		\quad \text{and} \quad
		\mathcal{V} = \{ -8 + 16j/N_v \}_{j=0}^{N_v-1}	
		\label{eq:kp-convergence-region}
	\end{align}
	where $v$ represents a discrete grid point on the domain $[-8,8]$, $k_x$ is the Fourier wavenumber in $x$, and $\omega_x = 2\pi / L_x = 1/10$. The convergent spatial modes lie inside the region
	\begin{align}
		\mathcal{B} = \left\{ (k_x, v) \in \mathbb{Z} \times \mathcal{V} : h|\mathbf{L}(k_x, v)| < r_1^{\text{max}}(h) \right\} 
		\quad \text{for} \quad
		h=50/2^{17},
		\label{eq:vp-convergence-region}
	\end{align}
	where the values of $r_1^{\text{max}}$ are contained in \cref{subtab:2d-parareal-configurations-vp}. Because we are only Fourier transforming in the $x$ direction, it is important that our Parareal configuration accurately computes all the components in the $v$ domain since we cannot assume spectral decay in the physical $v$ direction.
	
\begin{figure}[h]
	
		\centering	     
		
		\begin{subfigure}[b]{0.48\textwidth}
			\centering
			\includegraphics[width=\textwidth]{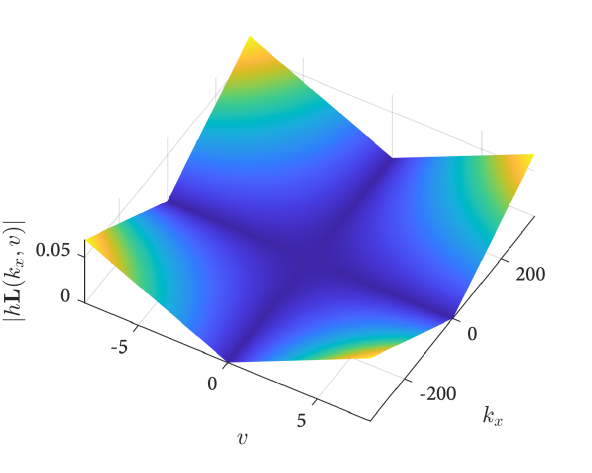}
			\caption{Magnitude of the VP linear operator eigenvalues}
			\label{fig:vp_linop}
		\end{subfigure}
		\hfill
		\begin{subfigure}[b]{0.48\textwidth}
			\centering
			\includegraphics[width=\textwidth]{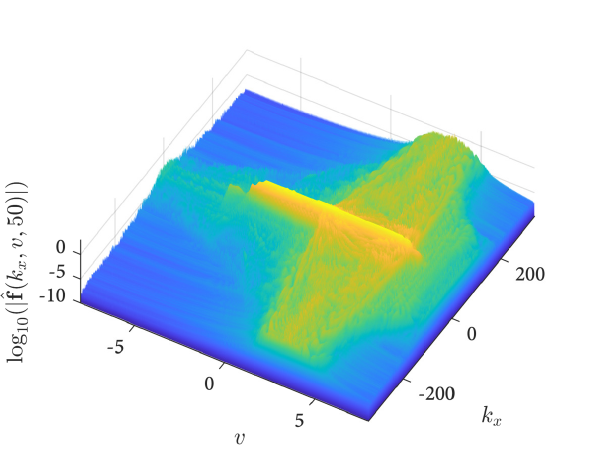}
			\caption{VP Solution at $t=50$ (Fourier transformed in $x$ only)}
			\label{fig:vp_fourier}
		\end{subfigure}
		
		\caption{Magnitude of the linear operator eigenvalues (\cref{fig:vp_linop}) and the transformed final solution (\cref{fig:vp_fourier}) for the Vlasov-Poisson equation. The axes in both plots represent the Fourier wavenumber in $x$ and the spatial variable $v$. The region $\mathcal{B}$ from \cref{eq:vp-convergence-region} for both the Parareal configurations from \cref{tab:2d-parareal-configurations} encloses the entire discrete $(k_x$, $v)$ domain; therefore no black contours are shown on the plots.}
		
		\label{fig:vp-convergence-region}
	
	\end{figure}

	 In \cref{fig:vp-convergence-region} we overlay the convergence region $\mathcal{B}$ onto the linear operator $\mathbf{L}$ and the transformed final solution. The convergence regions for Parareal configurations with both $\NGprop=1$ and $\NGprop=2$ enclose the entire discrete ($k_x$, $v$). We note that the largest diagonal element of the scaled linear operator $h\mathbf{L}$ is 0.0693, therefore the convergence region for the Parareal configuration with $\NGprop = 1$ just barely encloses the eigenvalues since $r_1^{\text{max}}(h) = 0.0740$.

\begin{figure}
	
		\hfill
		\begin{subfigure}[t]{0.32\textwidth}
			\centering
			\caption{Error versus Iteration}
			\includegraphics[width=\textwidth,trim={0 0 0.70cm 0}, clip]{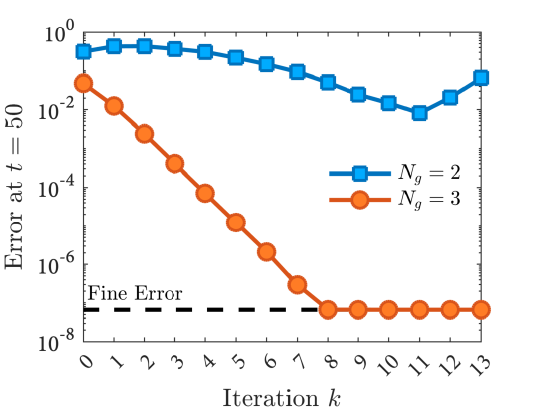}
			\label{fig:vp_error-vs-iteration}
		\end{subfigure}
		\hfill
		\begin{subfigure}[t]{0.32\textwidth}
			\centering
			\caption{Parallel Speedup \cref{eq:parareal_speedup}}
			\includegraphics[width=\textwidth,trim={0 0 0.70cm 0}, clip]{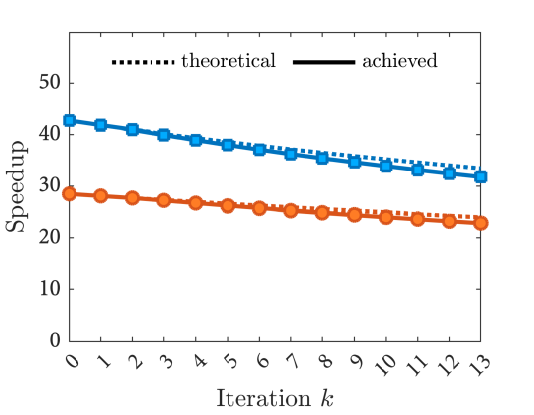}
			\label{fig:vp_parallel-speedup}
		\end{subfigure}
		\hfill
		\begin{subfigure}[t]{0.32\textwidth}
			\centering
			\caption{Error versus Runtime}
			\includegraphics[width=\textwidth,trim={0 0 0.70cm 0}, clip]{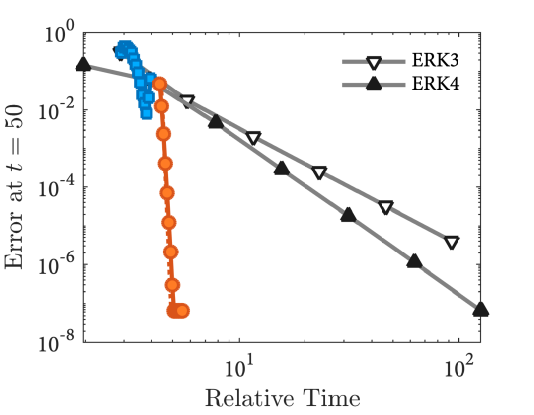}
			\label{fig:vp_error-vs-time}
		\end{subfigure}
		\hfill
		
		\caption{Numerical results for the VP equation using the two Parareal configurations from \cref{tab:2d-parareal-configurations}, that are differentiated with different colored lines. (a) Solution error at the final time $t=50$ as a function of Parareal iteration. (b) Theoretical speedup \cref{eq:parareal_speedup} and achieved speedup from the numerical experiment. (c) Error versus runtime diagram comparing the Parareal configurations to the serial coarse and fine integrators.}
		
		\label{fig:vp-ng-experiment}
		
	\end{figure}

	In \cref{fig:vp-ng-experiment} we show convergence, speedup, and error versus runtime plots for the exponential Parareal method applied to the VP equation. We again divide our discussion of the results into three parts:
	\begin{enumerate}
		
		\item {\em Convergence}. The Parareal method with $\NGprop=2$ failed to converge, while the method with $\NGprop=3$ displayed monotonic convergence and achieved the fine error tolerance after eight iterations. The failure of convergence for $\NGprop=2$ is likely due to several reasons. First, linear analysis is not guaranteed to provide an accurate prediction for all nonlinear equations. Moreover, linear analysis predicts that the Parareal with $\NGprop=2$ is only just barely  convergent, so a larger safety margin may be required to properly predict convergence on nonlinear problems. Lastly, it is also possible that the divergent iteration is due to instabilities. Specifically, our assumption that the nonlinear term is completely non-stiff may be inaccurate due to the presence of the term $-f_v$ in the nonlinearity. Fortunately, modestly increasing $\NGprop$ resolves these issues and leads to a stable, monotonically convergent Parareal iteration.
		
		\item {\em Parallel Speedup \cref{eq:parareal_speedup}}. As with the KP equation, we see very good agreement between the theoretical and achieved parallel speedup. We again see very minor penalties due to communication (notice that all dashed lines are slightly below the solid lines), however the differences are even smaller than those for the KP equation. This follows from the fact that the cost per timestep is more expensive for the VP equation than for the KP equation.

		\item {\em Error versus Runtime.} The Parareal configuration with $\NGprop = 3$ was able to obtain a solution with an error of $6.6 \times 10^{-8}$ twenty four times faster than the serial ERK4 method.	
			\begin{center}
				\renewcommand*{\arraystretch}{1.25}
				\begin{tabular}{l|ll}
					$\NGprop$ & Error Tol at $k=8$ & Improvement compared to serial ERK4 \\\hline			
					3 & $6.6 \times 10^{-8}$ & 24x faster		
				\end{tabular}
			\end{center}
			The improvement in time-to-solution of exponential Parareal over the serial ERK4 method is substantially greater on the VP equation than it was for the KP equation (24x vs 10x). This difference is made possible by the lack of highly-oscillatory temporal components in the VP equation. Specifically, this allowed us to run the coarse integrator at a significantly larger stepsize relative to the fine integrator. For comparison, the coarse integrator of the Parareal configurations with $\NGprop=3$ for the KP and VP equations, were respectively run with a stepsize that was 10.6 and 21.3 times larger than the fine integrator. Overall this experiment demonstrates the potential for very significant reduction in time-to-solution when applying exponential Parareal to accurately solve hyperbolic equations.
			 
	\end{enumerate}

\section{Conclusions and future work}
\label{sec:conclusion}

In this paper we applied exponential integrators within the Parareal iteration and presented linear analysis that can be used to study the stability and convergence properties of the resulting methods on non-diffusive equations. We then demonstrated that exponential Parareal methods can achieve significantly reduced time-to-solution compared to serial exponential integrators on non-diffusive partial differential equations. 

We draw two main conclusions from this work. First we showed that repartitioning is essential for obtaining a Parareal configuration that is stable on stiff non-diffusive equations. Second, through linear analysis we were able to better understand the convergence characteristics of the Parareal iteration in the absence of diffusion. Specifically we saw that the Parareal iteration achieves fine integrator accuracy for low-frequency (i.e. non-stiff) oscillatory modes and coarse integrator accuracy for high-frequency (i.e. stiff) oscillatory modes. When solving non-diffusive partial differential equations this phenomenon makes it impossible to guarantee rapid convergence for high-frequency spatial modes. Therefore, exponential Parareal is best suited for non-diffusive equations and initial conditions that do not cause rapid spectral broadening.

To the best of the authors' knowledge, this is the first paper to investigate the usage of ETD-RK methods within the Parareal iteration. Our initial results look promising as we have demonstrated the ability to achieve reduced time-to-solution using exponential Parareal on both hyperbolic and dispersive equations. Nevertheless, there are still many avenues that require further exploration. In particular all of the numerical experiments presented in this paper involve diagonal linear operators that greatly simplify the computation of the exponential $\varphi$-functions. In future work we plan to study exponential Parareal integrators in the more general setting with non-diagonal linear operators and examine the resulting effects on computational performance.

\section*{Acknowledgements}
The work of Buvoli was funded by the National Science Foundation, Computational Mathematics Program DMS-2012875. The work of Minion was supported by the U.S. Department of Energy, Office of Science, Office of Advanced Scientific Computing Research, Applied Mathematics program under contract number DE-AC02005CH11231. 

\bibliographystyle{siam}
\bibliography{references_nourl,pint}

\appendix
\section{Method coefficients}
\label{app:coefficients}

This appendix contains the Butcher tableaux described in \cref{subsec:erk} for the exponential Runge-Kutta integrators that are used in this paper. We use the abbreviations $\varphi_{i} = \varphi_i(h\mathbf{L})$ and $\varphi_{i,j} = \varphi_i(c_jh\mathbf{L})$ for the $\varphi$-functions.

\begin{itemize}
	\item {\bf ERK1:} first-order exponential Euler method \cref{eq:exponential-euler}
	\begin{center}
		\renewcommand*{\arraystretch}{1.5}
		\begin{tabular}{l|llll}
			$0$ 			&   \\ \hline
							& $\varphi_1$
		\end{tabular}
	\end{center}
	\item {\bf ERK2:} second-order method from \cite{cox2002ETDRK4}
	\begin{center}
	\renewcommand*{\arraystretch}{1.5}
		\begin{tabular}{l|llll}
			0 				& \\
			$1$ 			& $\varphi_{1}$  \\ \hline
							& $\varphi_1 - \varphi_2$ 		& $\varphi_2$ 
		\end{tabular}
	\end{center}
	\item {\bf ERK3:} third-order method from \cite{cox2002ETDRK4}	
	\begin{center}
		\renewcommand*{\arraystretch}{1.5}
		\begin{tabular}{l|llll}
			0 				& \\
			$\tfrac{1}{2}$ 	& $\tfrac{1}{2} \varphi_{1,2}$  \\
			$1$				& $-\varphi_{1}$ 					& $2\varphi_1$ & \\ \hline
							& $\varphi_1 - 3 \varphi_2 + 4 \varphi_3$ 		& $4\varphi_2 - 8 \varphi_3$ 	& $- \varphi_2 + 4 \varphi_3$
		\end{tabular}
	\end{center}
	\item {\bf ERK4:} fourth-order method from \cite{krogstad2005IF}
	\begin{center}
		\renewcommand*{\arraystretch}{1.5}
		\begin{tabular}{l|llll}
			0 				& \\
			$\tfrac{1}{2}$ 	& $\tfrac{1}{2} \varphi_{1,2}$  \\
			$\tfrac{1}{2}$ 	& $\tfrac{1}{2} \varphi_{1,2} - \varphi_{2,2}$	& $\varphi_{2,2}$				& \\
			$1$				& $\varphi_{1} - 2 \varphi_2$ 					& 0 							& $2\varphi_2$ \\ \hline
							& $\varphi_1 - 3 \varphi_2 + 4 \varphi_3$ 		& $2\varphi_2 - 4 \varphi_3$ 	& $2\varphi_2 - 4 \varphi_3$ & $-\varphi_2 + 4 \varphi_3$
		\end{tabular}
	\end{center}
\end{itemize}

\newpage
\section{Nonlinear Schr\"{o}dinger Serial ERK Results}
\label{app:nls-serial-convergence}

\begin{figure}

	\begin{minipage}[t]{0.48 \linewidth}
		\centering
		{\bf Error versus Stepsize} -- Initial Condition \cref{eq:nls-low-frequency-ic}
		
		\includegraphics[width=\linewidth]{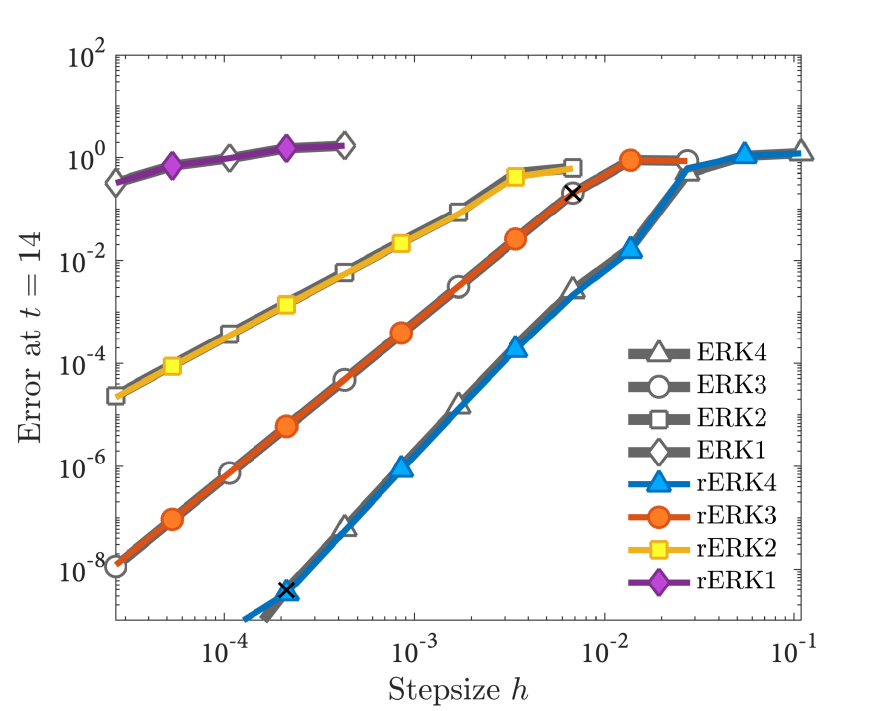}
	\end{minipage}
	\begin{minipage}[t]{0.48 \linewidth}
		\centering
		{\bf Error versus Computational Time} -- Initial Condition \cref{eq:nls-low-frequency-ic}

		\includegraphics[width=\linewidth]{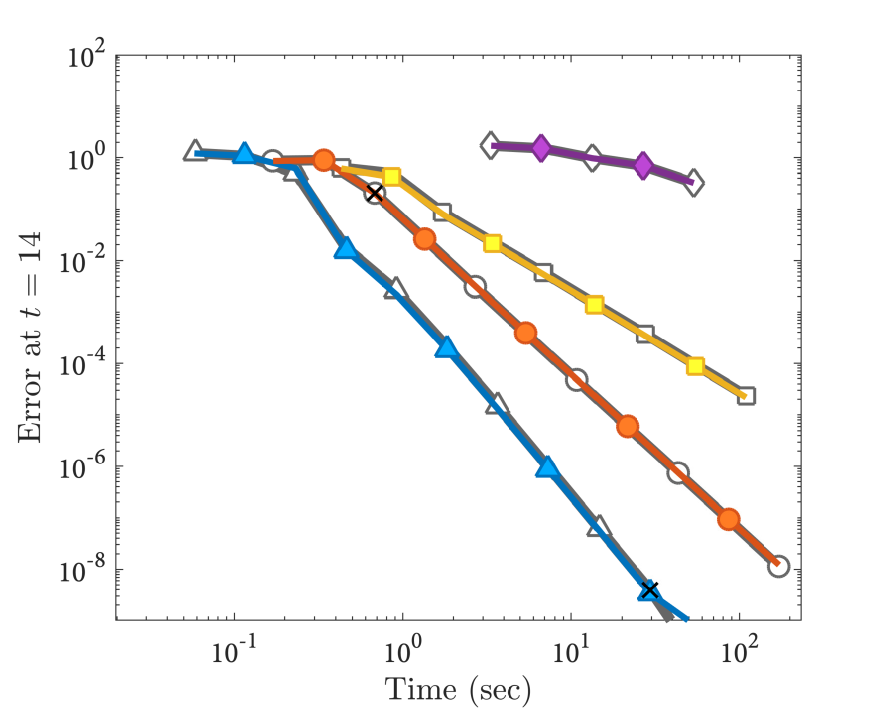}
	\end{minipage}
	\hfill
	
	\caption{Convergence diagram (left) and precision diagram (right) for the exponential Runge-Kutta methods listed in \cref{app:coefficients} run on the nonlinear Schr\"{o}dinger equation \cref{eq:nls} with initial condition \cref{eq:nls-low-frequency-ic}. Colored lines correspond to repartitioned integrators (rERK) and gray lines correspond to unmodified exponential integrators (ERK) -- repartitioning has no effect on this problem.  The black crosses on the ERK3 and ERK4 method respectively correspond to the stepsizes of the coarse and fine integrators for the Parareal method described in \cref{tab:nls-parareal-configuration}.}
	\label{fig:nls-serial-convergence-smooth} 
		
\end{figure}

\begin{figure}

	\begin{minipage}[t]{0.48 \linewidth}
		\centering
		{\bf Error versus Stepsize} -- Initial Condition \cref{eq:nls-high-frequency-ic}
		
		\includegraphics[width=\linewidth]{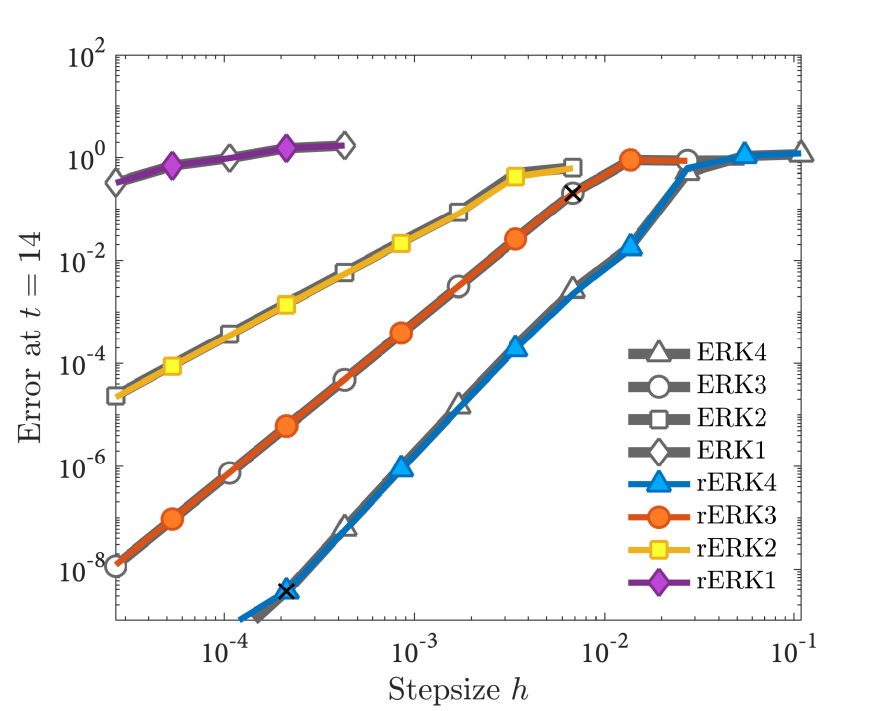}
	\end{minipage}
	\begin{minipage}[t]{0.48 \linewidth}
		\centering
		{\bf Error versus Computational Time} -- Initial Condition \cref{eq:nls-high-frequency-ic}

		\includegraphics[width=\linewidth]{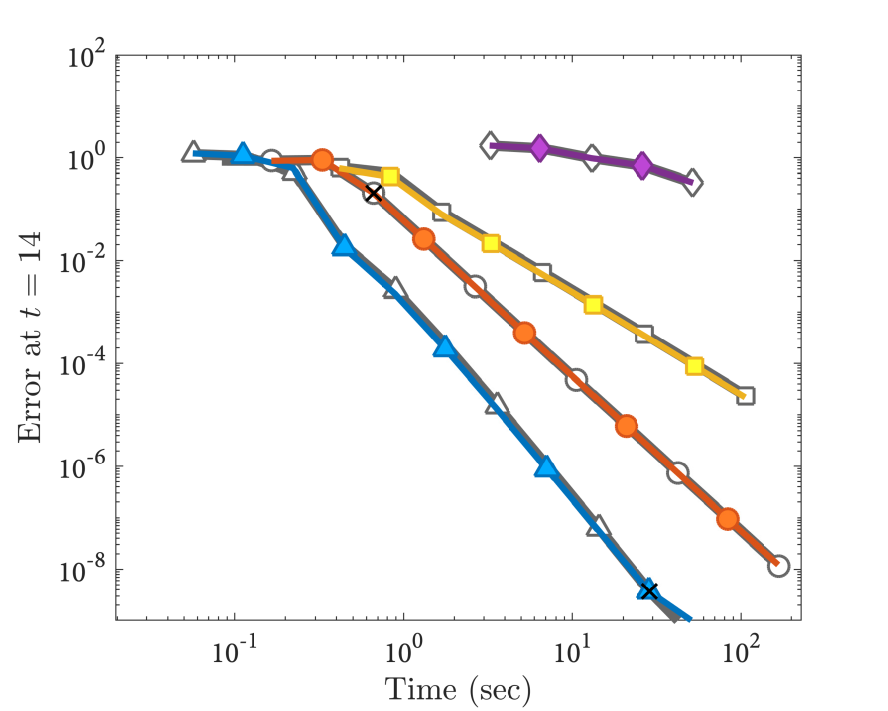}
	\end{minipage}
	\hfill
	
	\caption{Identical to \cref{fig:nls-serial-convergence-smooth} except we are now considering the initial condition \cref{eq:nls-high-frequency-ic}.}
	\label{fig:nls-serial-convergence-oscillatory} 
		
\end{figure}

\begin{figure}

	\begin{minipage}[t]{0.48 \linewidth}
		\centering
		{\bf Error vs Stepsize} -- Initial Condition \cref{eq:nls-full-spectrum-ic}
		
		\includegraphics[width=\linewidth]{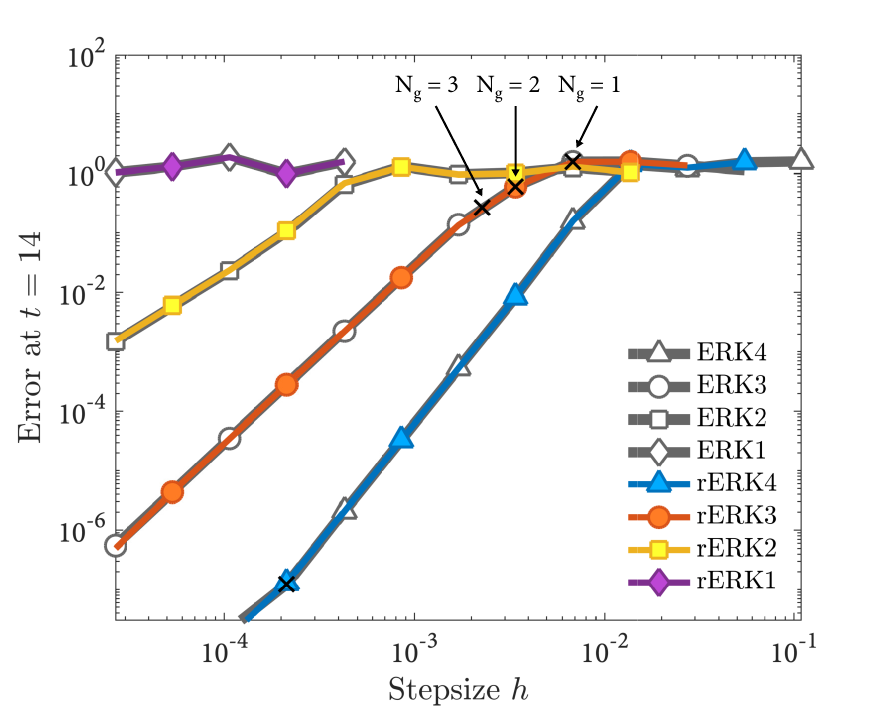}
	\end{minipage}
	\begin{minipage}[t]{0.48 \linewidth}
		\centering
		{\bf Error vs Computational Time} -- Initial Condition \cref{eq:nls-full-spectrum-ic}

		\includegraphics[width=\linewidth]{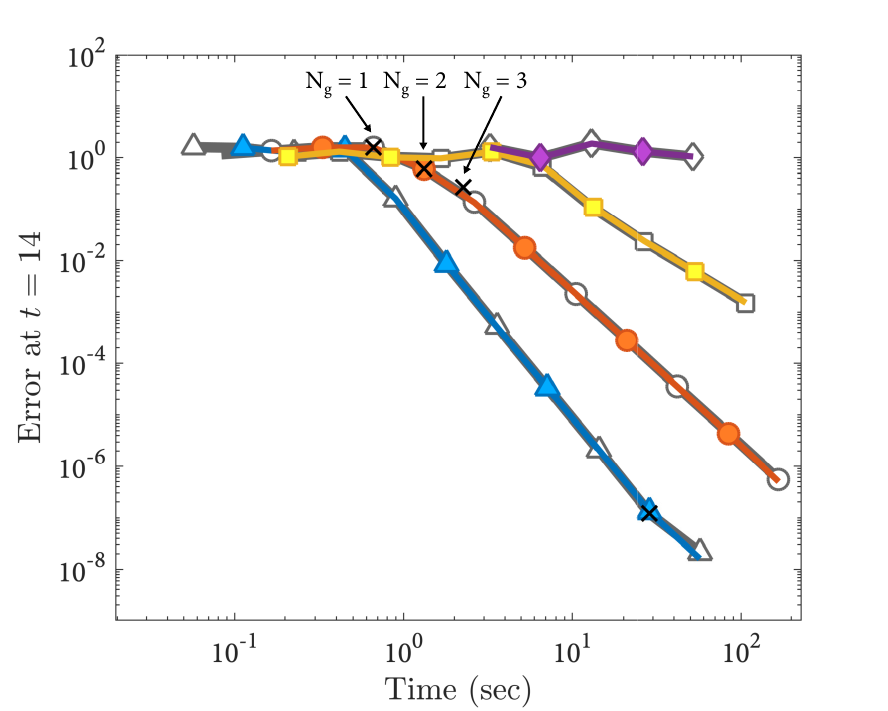}
	\end{minipage}

	\caption{Identical to \cref{fig:nls-serial-convergence-smooth} except we are now considering the initial condition \cref{eq:nls-full-spectrum-ic}.}
	\label{fig:nls-serial-convergence-full-spectrum} 
		
\end{figure}

We solve the nonlinear Schr\"{o}dinger equation \cref{eq:nls} using the serial ERK methods from from \cref{app:coefficients} using $2^p$  timesteps where $p=7, 8, \ldots, 19$. In \cref{fig:nls-serial-convergence-smooth,fig:nls-serial-convergence-oscillatory,fig:nls-serial-convergence-full-spectrum} we show accuracy and convergence diagrams for the initial conditions \cref{eq:nls-low-frequency-ic,eq:nls-high-frequency-ic,eq:nls-full-spectrum-ic}, respectively.

\section{Parareal with $\NIters=0,\ldots,160$ for NLS with initial condition \cref{eq:nls-high-frequency-ic}}
\label{app:nls-experiment-oscillatory-extended-k}

\begin{figure}
	
	\includegraphics[width=\linewidth]{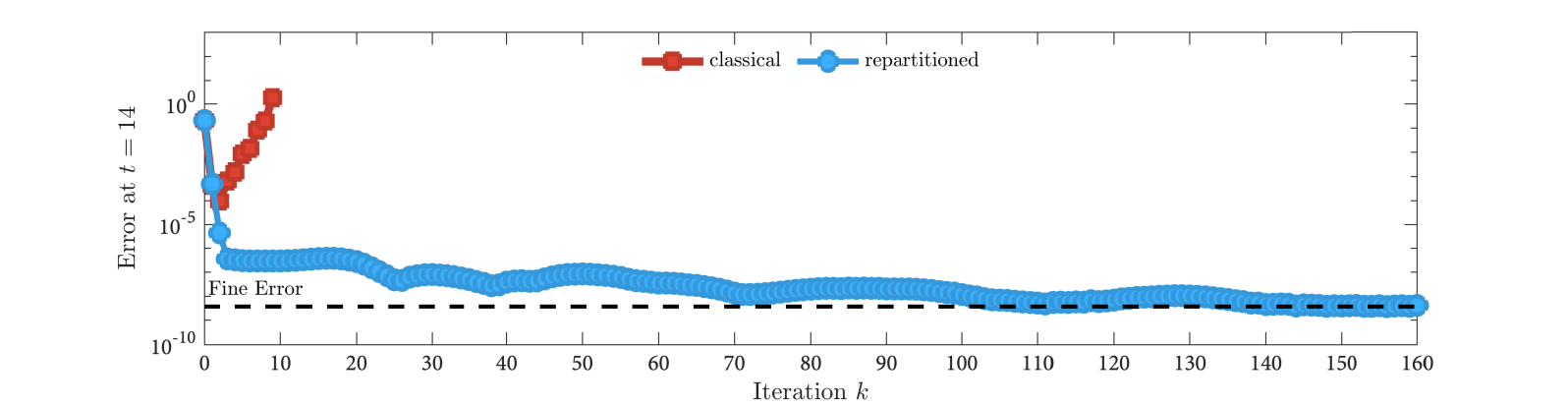}	
	
	\caption{Error of the Parareal iteration from \cref{tab:nls-parareal-configuration} with $K=0,\ldots,160$ applied to the NLS equation \cref{eq:nls} with initial conditions \cref{eq:nls-high-frequency-ic}. Parareal with classical ERK becomes unstable (i.e. computation produces NaN values) for $k>9$.}
	\label{fig:nls-experiment-oscillatory-extended-k}	
\end{figure}

In \cref{fig:nls-experiment-oscillatory-extended-k} we show error versus iteration for a larger number of Parareal iterations than shown in \cref{fig:nls-experiment-oscillatory}. Specifically, we solve the nonlinear Schr\"{o}dinger equation \cref{eq:nls} with initial conditions \cref{eq:nls-high-frequency-ic} using the Parareal configuration from \cref{tab:nls-parareal-configuration}, except with $\NIters=0,\ldots,160$.  We see that Parareal with classical ERK becomes completely unstable for $k>9$, while Parareal with repartioned ERK converges after $150$ iterations achieving a maximum theoretical speedup of $10.16$.

\section{NLS solution for initial condition \cref{eq:nls-full-spectrum-ic}}
\label{app:nls-full-spectrum-solution}

\begin{figure}
	
	\begin{minipage}[t]{0.48 \linewidth}
		\centering
		{\bf NLS Solution} -- Initial Condition \cref{eq:nls-low-frequency-ic}
		
		\includegraphics[width=\linewidth]{figures/nls/nls-solution-smooth}	
	\end{minipage}
	\begin{minipage}[t]{0.48 \linewidth}
		\centering
		{\bf NLS Solution} -- Initial Condition \cref{eq:nls-full-spectrum-ic}

		\includegraphics[width=\linewidth]{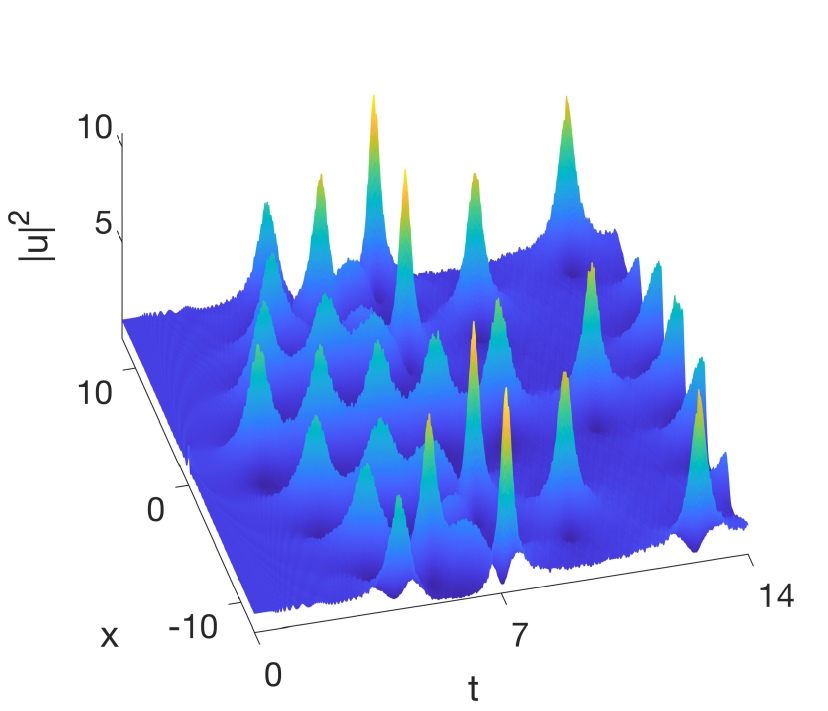}
	\end{minipage}
	
	\caption{Solutions of the nonlinear Schr\"{o}dinger equation \cref{eq:nls} for the initial conditions \cref{eq:nls-low-frequency-ic} and \cref{eq:nls-full-spectrum-ic}.}
	\label{fig:nls-solution-visualization}
\end{figure}

Figure \ref{fig:nls-solution-visualization} shows the two different NLS solutions arising from the initial conditions \cref{eq:nls-low-frequency-ic} and \cref{eq:nls-full-spectrum-ic}.

\section{Additional Stability Plots}
\label{app:additional-stability-plots}

\begin{figure}

	\centering
	
	Parareal Stability Regions and Instability Factors
	\rule{\textwidth}{0.5px} \\[0.5em]
	
	\begin{tabular}{cccccc}
	& $\NIters=0$ & $\NIters = 2$ & $\NIters=4$ & $\NIters=6$ \\ 
		\rotatebox{90}{\hspace{2em} Classical ERK} & 
		\includegraphics[height=.2\textwidth,align=b]{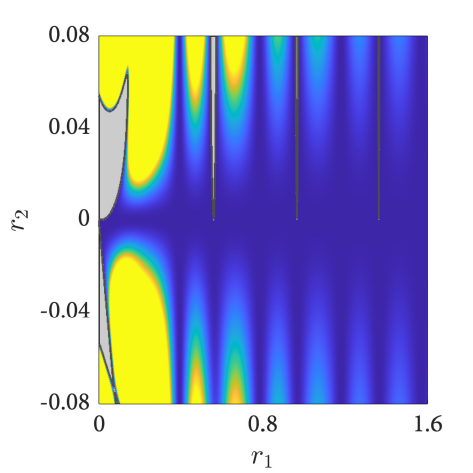} &
		\includegraphics[height=.2\textwidth,align=b]{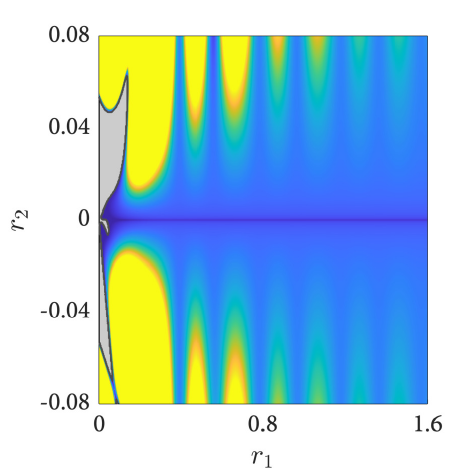} &
		\includegraphics[height=.2\textwidth,align=b]{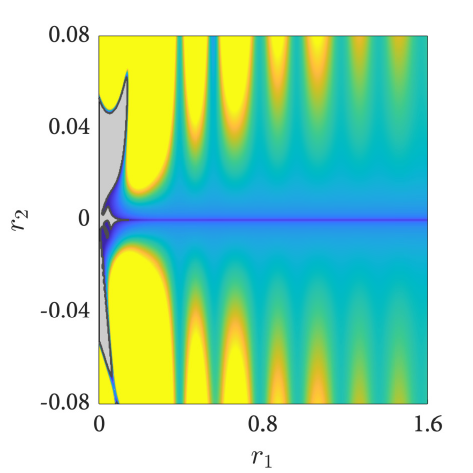} &
		\includegraphics[height=.2\textwidth,align=b]{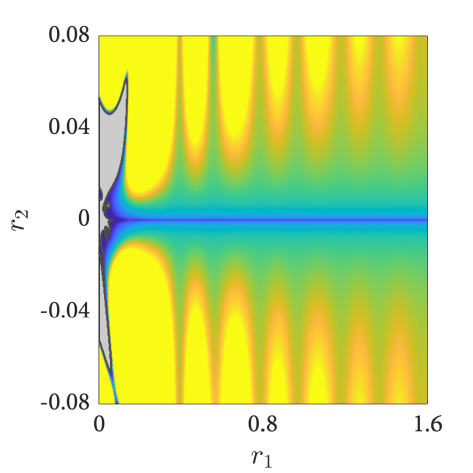} &
		\includegraphics[height=.2\textwidth,trim={6cm 0 0 0},clip,align=b]{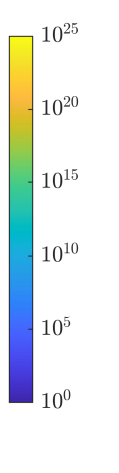} \\
		
		\rotatebox{90}{\hspace{1em} Repartitioned ERK} &
		\includegraphics[height=.2\textwidth,align=b]{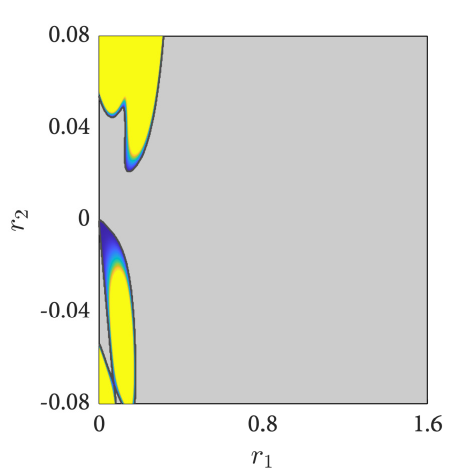} &
		\includegraphics[height=.2\textwidth,align=b]{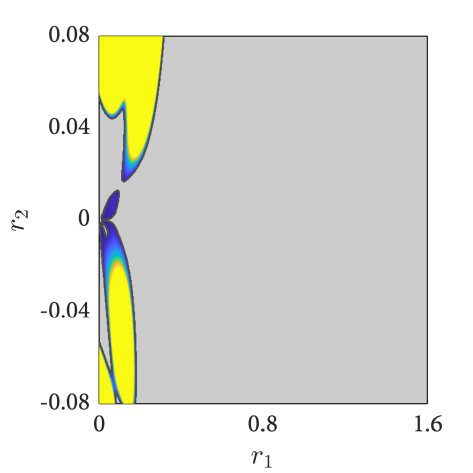} &
		\includegraphics[height=.2\textwidth,align=b]{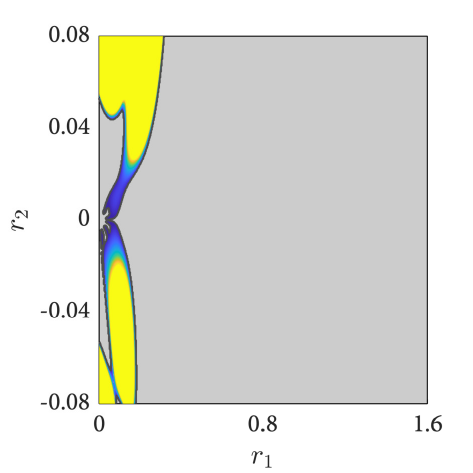} &
		\includegraphics[height=.2\textwidth,align=b]{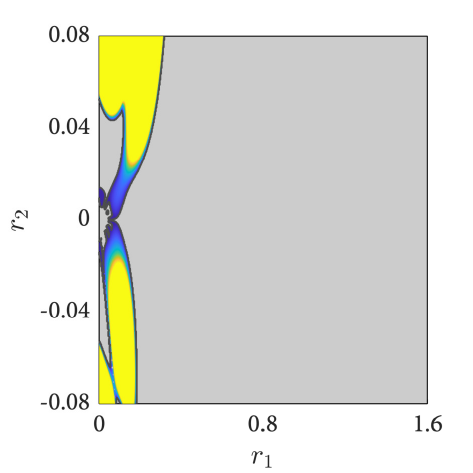} &
		\includegraphics[height=.2\textwidth,trim={6cm 0 0 0},clip,align=b]{figures/nls/nls-stability/stability-surface-demag-leg-k-0-rho-0-1}
	\end{tabular}

	\caption{Additional stability regions for the Parareal configuration from \cref{tab:nls-parareal-configuration} with classical exponential integrators (top column) and repartitioned exponential integrators (bottom column). The gray region is the stability region \cref{eq:stability-region}, and color shows the amplification factor $|R(z_1=ir_1,z_2=ir_2)|$ outside the stability region where the method is unstable. These figures show a wider range of  $r_2$ values than \cref{fig:parareal-stability-nls}.
	}
	
	\label{fig:parareal-stability-extra1}

\end{figure}

\begin{figure}

	\centering
	
	Parareal Stability Regions and Instability Factors
	\rule{\textwidth}{0.5px} \\[0.5em]
	
	\begin{tabular}{cccccc}
	& $\NIters=0$ & $\NIters = 2$ & $\NIters=4$ & $\NIters=6$ \\ 
		\rotatebox{90}{\hspace{2em} Classical ERK} & 
		\includegraphics[height=.2\textwidth,align=b]{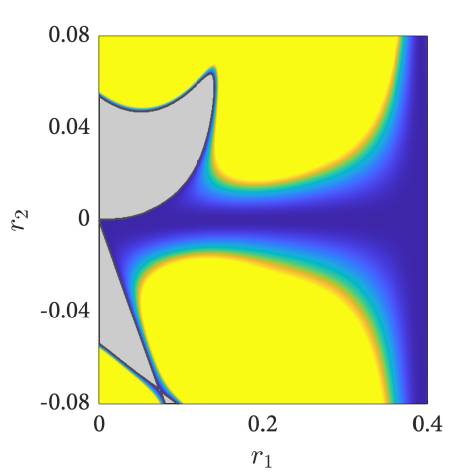} &
		\includegraphics[height=.2\textwidth,align=b]{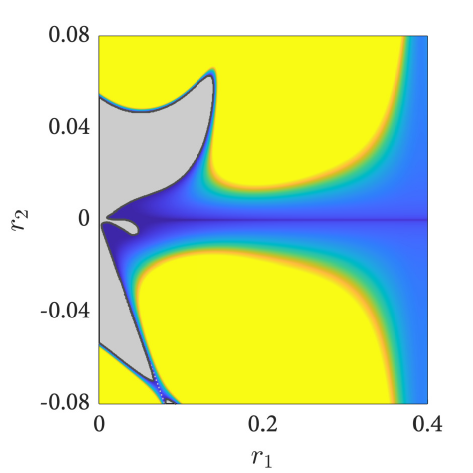} &
		\includegraphics[height=.2\textwidth,align=b]{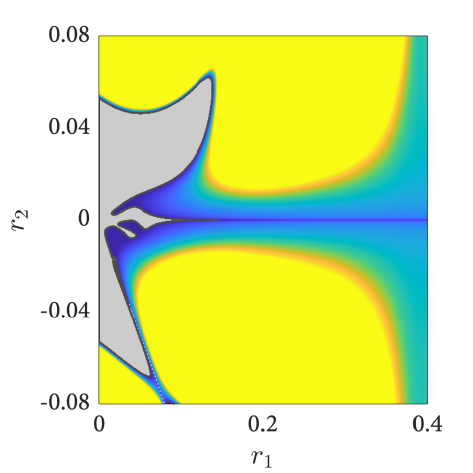} &
		\includegraphics[height=.2\textwidth,align=b]{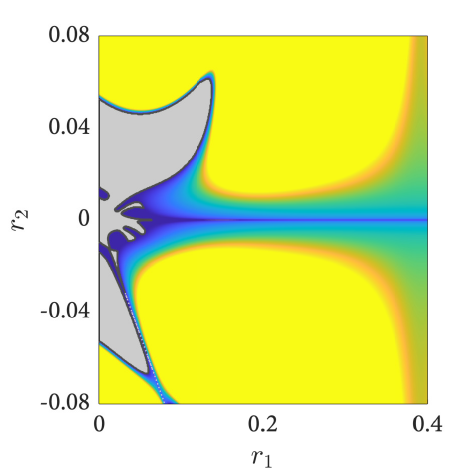} &
		\includegraphics[height=.2\textwidth,trim={6cm 0 0 0},clip,align=b]{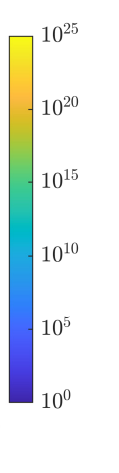} \\
		
		\rotatebox{90}{\hspace{1em} Repartitioned ERK} &
		\includegraphics[height=.2\textwidth,align=b]{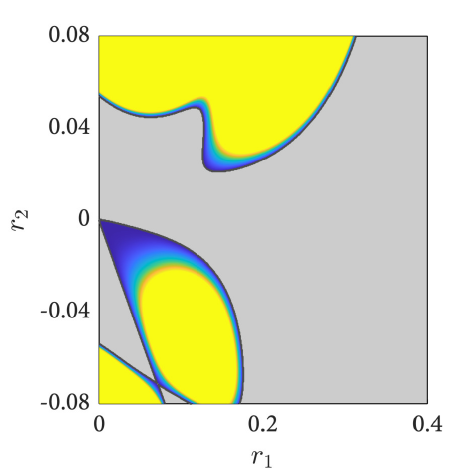} &
		\includegraphics[height=.2\textwidth,align=b]{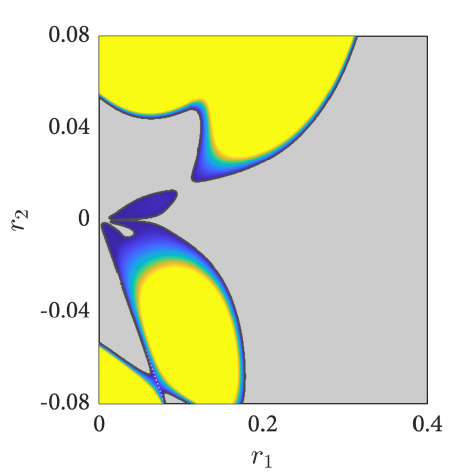} &
		\includegraphics[height=.2\textwidth,align=b]{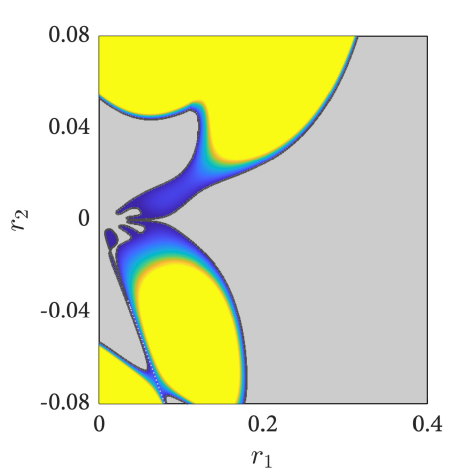} &
		\includegraphics[height=.2\textwidth,align=b]{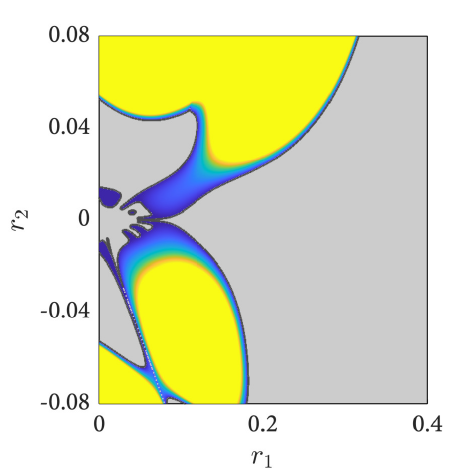} &
		\includegraphics[height=.2\textwidth,trim={6cm 0 0 0},clip,align=b]{figures/nls/nls-stability/stability-surface-demag-r2-mag-r1-leg-k-0-rho-0-1}
	\end{tabular}

	\caption{Additional stability regions for the Parareal configuration from \cref{tab:nls-parareal-configuration} with classical exponential integrators (top column) and repartitioned exponential integrators (bottom column). The grey region is the stability region \cref{eq:stability-region}, and color shows the amplification factor $|R(z_1=ir_1,z_2=ir_2)|$ outside the stability region where the method is unstable. These figures show a narrrower range of $r_1$ and a wider range of $r_2$ than \cref{fig:parareal-stability-nls}.
	}
	
	\label{fig:parareal-stability-extra2}

\end{figure}

\Cref{fig:parareal-stability-extra1,fig:parareal-stability-extra2} contain additional stability plots for the Parareal configuration  \cref{tab:nls-parareal-configuration} that show different $(r_1,r_2)$ ranges than \cref{fig:parareal-stability-nls}.

\section{Remark regarding convergence region scaling}
\label{app:convergence-scaling}

	\begin{remark}{Stability regions grow approximately linearly in $\NGprop$ for small $(|r_1|,|r_2|)$.}
		\label{rem:conv-ng}
		To show this, we first assume that we are in a regime where $(|r_1|,|r_2|)$ is small so that the coarse and the fine integrator both exhibit asymptotic error properties. If we construct the coarse propagator $\Gprop$ using $\NGprop$ steps of a $q$th order integrator, then
		\begin{align}
			|R_{\mathcal{G}}| &= \left|e^{i(r_1+r_2)}\right| + \mathcal{O}\left(\tfrac{|r_1|+|r_2|}{\NGprop}\right)^q < 1 + C_1\left(\tfrac{|r_1|+|r_2|}{\NGprop}\right)^q, \\ 
			|R_{\mathcal{G}} - R_{\mathcal{F}}| &= \mathcal{O}\left(\tfrac{|r_1| + |r_2|}{\NGprop}\right)^q < C_2\left(\tfrac{|r_1| + |r_2|}{\NGprop}\right)^q.
		\end{align} 
		The norm of the error matrix can be bounded above by
		\begin{align}
			\| \mathbf{E} \|_\infty = \frac{1 - |R_{\mathcal{G}}|^{\Nprocs}}{1 - |R_{\mathcal{G}}|} |R_{\mathcal{G}} - R_{\mathcal{F}}| <
			\Nprocs C_2\left(\tfrac{|r_1| + |r_2|}{\NGprop}\right)^q + \mathcal{O}\left(\tfrac{|r_1|+|r_2|}{\NGprop}\right)^q.
		\end{align}
		Convergence is guaranteed if $\| \mathbf{E} \|_\infty < 1$, which is equivalent to
		\begin{align}
			|r_1| + |r_2| < \frac{1}{\NGprop C_2 \Nprocs}	+ \mathcal{O}\left(\tfrac{|r_1|+|r_2|}{\NGprop}\right)^q.
		\end{align}
		Ignoring the higher order terms, the size of the region $|r_1| + |R_2| < (C_2 \NGprop \Nprocs)^{-1}$ grows linearly in $\NGprop$.		
	\end{remark}

\section{Spatially discretized Vlasov-Poisson equation}
\label{app:discrete-vp}

For notational simplicity we represent the discrete VP solution as the matrix $\mathbf{f}$ where $\mathbf{f}_{jk}$ approximates the continuous solution at the grid point ($v_j$, $x_k$). Next we define the scaled Fourier wavenumber vectors
	\begin{align}
		\mathbf{k}^v = 	\frac{\pi}{L_v} [0, \ldots, N_v/2-1, -N_v/2, \ldots, -1], &&
		\mathbf{k}^x = \frac{2\pi}{L_x} [0, \ldots, N_x/2-1, -N_x/2, \ldots, -1],		
	\end{align}
	for $L_v = 16$, $L_x = 20\pi$, and the $\mathbb{R}^{N_v,Nx}$ matrices $\mathbf{K}^{v}_{jk} = \mathbf{k}^v_j$, $\mathbf{K}^{x}_{jk} = \mathbf{k}^x_k$. If $\mathcal{F}_v(\cdot)$ and $\mathcal{F}_x(\cdot)$ represent the discrete Fourier transform in $v$ and $x$, then the transformed variable $\widehat{\mathbf{f}} = \mathcal{F}_x(\mathbf{f})$ satisfies
\begin{align}
	\frac{d}{dt}\widehat{\mathbf{f}}_{jk} = v_j i \mathbf{k}^{x}_k \widehat{\mathbf{f}}_{jk} + \mathcal{F}_x\left(\mathbf{E} \text{.*} \frac{\partial \mathbf{f}}{\partial v}\right)_{jk} &&	
	\frac{\partial \mathbf{f}}{\partial v} = \mathcal{F}^{-1}_x\left(\mathcal{F}^{-1}_v\left( i \mathbf{K}^v \text{.*} \mathcal{F}_v\left(\widehat{\mathbf{f}}\right)\right)\right)
	\label{eq:discrete-vp}
\end{align}
where .* denotes the Hadamard product, and the discrete electric field $\mathbf{E}_{jk} = \mathbf{e}_k \approx E(x_k, t)$ is 
\begin{align}
	\mathbf{e} = \mathcal{F}^{-1}_x \Bigg( \mathbf{k^{-x}} \text{.*} \underbrace{ \left( \mathbf{b} +  \Delta v \sum_{j=1}^{N_v} \widehat{\mathbf{f}}_{jk} \right) }_{\mathcal{F}_x(-1 + \int_0^{20\pi} f(x,v,t)dv)}\Bigg),
	&& \text{for} &&
	\begin{aligned}
		\mathbf{b} &= \mathcal{F}_x(-[1, \ldots, 1]^T) \in \mathbb{R}^{N_x}, \\	
		\mathbf{k}^{-x}_{k} &= 
		\begin{cases}
 			0 & j = 1, \\
 			1/\mathbf{k}^x_k & j > 1
 		\end{cases} \in \mathbb{R}^{N_x}, \\
 		\Delta v &= 16 / N_v.
	\end{aligned}
\end{align}
Note that the integral term $\int_0^{20\pi} f(x,v,t)dv$ is treated using the trapezoidal rule, which convergences exponentially on periodic domains.

\section{ERK convergence diagrams for KP and Vlassov-Poisson}
\label{app:erk-convergence-kp-vp}

\begin{figure}

	\begin{minipage}[t]{0.48 \linewidth}
		\centering
		{\bf Error versus Stepsize} -- KP equation \cref{eq:kp-evolution}
		
		\includegraphics[width=\linewidth]{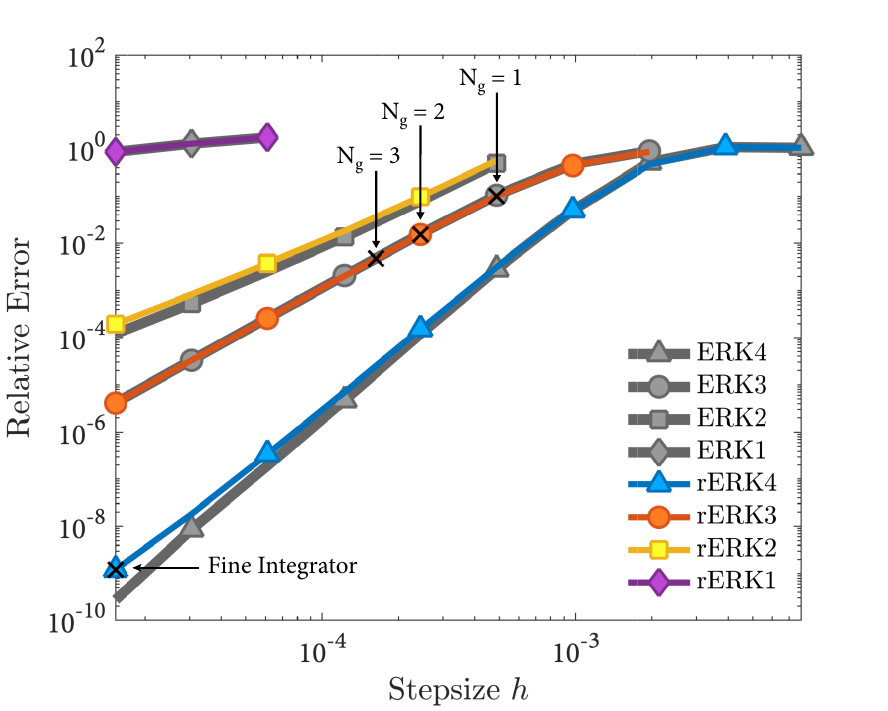}
	\end{minipage}
	\begin{minipage}[t]{0.48 \linewidth}
		\centering
		{\bf Error versus Stepsize} -- VP Equation \cref{eq:valsov-poisson}
		
		\includegraphics[width=\linewidth]{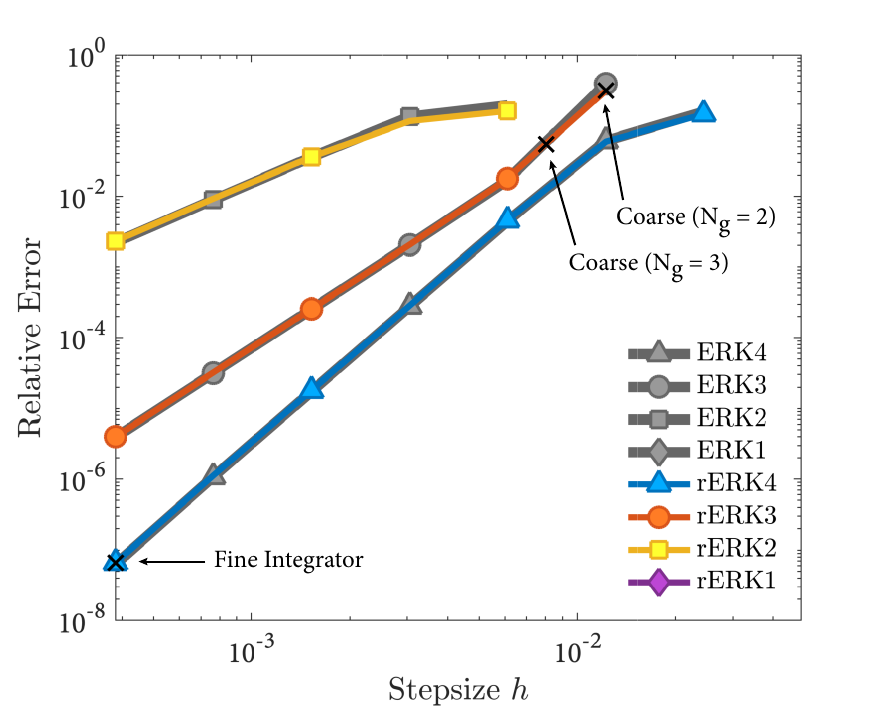}
	\end{minipage}
	\hfill
	
	\caption{Serial ERK convergence diagrams for the KP and VP equations. We can related these plots to the Parareal configurations described in \cref{tab:2d-parareal-configurations}. Specifically, the labeled black crosses at large timesteps correspond to the stepsizes of the coarse integrator, while the black cross on the ERK4 method at the smallest stepsize corresponds to the stepsize of the fine integrator. }
	\label{fig:erk-convergence-kp-vp} 
		
\end{figure}

In \cref{fig:erk-convergence-kp-vp} we show convergence diagrams for the serial ERK methods applied to the KP and VP equations. The plots also contain black crosses that indicate the step-sizes of the coarse and fine integrators for the Parareal configurations described in \cref{sec:higher-dimensional-numerical-experiments}.

\end{document}